\documentclass[11pt]{article}

\usepackage{anysize}
\usepackage{paralist}

\usepackage{float}
\usepackage{fancybox}

\usepackage{pstricks}       

\usepackage{amsmath,amssymb} %
\usepackage{amsthm}          
\usepackage{enumerate} %
\usepackage{url}             
\usepackage{rotating}        

\usepackage{bbm}
\usepackage{ae,aecompl}      
\usepackage[mathscr]{eucal}  

\usepackage{graphicx}        
\usepackage{subfigure}       
\usepackage{psfrag}          
\usepackage{booktabs}        

\usepackage{sidecap}
\usepackage[footnotesize,sc]{caption2}  

\usepackage{tabularx}         %
\usepackage{xspace}         
\usepackage{makeidx}          
\usepackage{index}            

\usepackage{constructions}

\usepackage{color}       
\definecolor{darkgreen}{rgb}{0,.6,0}

\let\orgbfseries\bfseries
\renewcommand{\bfseries}{\orgbfseries\sffamily}

\DeclareMathOperator{\verticesOP}{vert}

\DeclareMathOperator{\affOP}{aff}
\DeclareMathOperator{\convOP}{conv}
\DeclareMathOperator{\coneOP}{cone}
\DeclareMathOperator{\relintOP}{relint}

\newcommand{\DerivativeComplex}[1]{\ensuremath{\mathcal{D}({#1})}}
\newcommand{\ComplexOf}[1]{\ensuremath{\mathcal{C}({#1})}}
\newcommand\Schlegel{\mathop{\textsc{Schlegel}}}
\newcommand\SchlegelCap{\mathop{\textsc{SchlegelCap}}}

\newcommand{\LiftedPrism}[2]{\mathop{\textsc{LiftedPrism}}(#1,#2)}
\newcommand{\LiftedPrismOverTwoBalls}[4]{\mathop{\textsc{LiftedPrism}}((#1,#2),(#3,#4))}
\newcommand\Glue{\mathop{\textsc{Join}}}
\newcommand\conv{\mathop{\operatorname{conv}}}
\newcommand{\relint}[1]{\ensuremath{\relintOP(#1)}}
\newcommand{\simplicialBarycentricSubdiv}[1]{\mathop{\rm sd}{({#1})}}
\newcommand\R{{\mathbb R}}

\newcommand{\vertices}[1]{\ensuremath{\verticesOP({#1})}} 
\newcommand{\facets}[1]{\ensuremath{\operatorname{fac}({#1})}} 
\newcommand{\support}[1]{\ensuremath{|#1|}}

\newcommand{\setdef}[2]{\ensuremath{\left\{{#1}\,:\,{#2}\right\}}}
\newcommand{\eps}{\varepsilon}
\newcommand{\immersed}{\looparrowright}

\newcommand{\boundaryOP}{\ensuremath{\partial}}
\newcommand{\lifted}[2]{\operatorname{lift}(#1,#2)}

\newcommand{\unitvector}[1]{\ensuremath{\boldsymbol{e}_{#1}}}

\newcommand{\PileOfCubes}[2]{\ensuremath{P_{#1}(#2)}}

\newcommand{\faces}[2]{\ensuremath{\mathcal{F}_{#1}({#2})}} 

\theoremstyle{plain}
\newtheorem{theorem}{Theorem}[section]
\newtheorem{lemma}[theorem]{Lemma}
\newtheorem{corollary}[theorem]{Corollary}
\newtheorem{proposition}[theorem]{Proposition}

\theoremstyle{definition}
\newtheorem{definition}[theorem]{Definition}

\newtheorem*{example}{Example}
\newtheorem*{remark*}{Remark}

\theoremstyle{remark}
\newtheorem*{notation}{Notation}

  \newenvironment{Figure}[1][]{%
    \begin{figure}[H]\centering
    \begin{psfrags}}{%
   \end{psfrags}\end{figure}}
  \newenvironment{FigureHere}[1][]{%
    \begin{figure}[H]\centering
    \begin{psfrags}}{%
   \end{psfrags}\end{figure}}

\newcommand{\transformsToArrow}[2]{
    \raisebox{#1}{\includegraphics[width=#2]{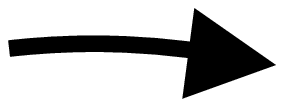}}}

{\end{list}}

{\end{list}}

{\end{list}}

{\end{list}}

\newenvironment{thm_enumerate_i}{%
\begin{list}{\rm(\roman{enumi})}%
{\usecounter{enumi}%
\setlength{\topsep}{2mm}%
\setlength{\partopsep}{0mm}%
\setlength{\parskip}{0mm}%
\setlength{\parsep}{0mm}%
\setlength{\itemsep}{2mm}%
\settowidth{\labelwidth}{(viii)}%
\setlength{\leftmargin}{0mm}%
\addtolength{\leftmargin}{\labelwidth}%
\addtolength{\leftmargin}{\labelsep}%
\setlength{\itemindent}{0mm}}}%
{\end{list}}

\newenvironment{inputoutput}{
\begin{list}{\bf ???}{%
\settowidth{\labelwidth}{Output: }
\settowidth{\leftmargin}{Output:\ }%
\addtolength{\leftmargin}{\labelsep}%
\setlength{\topsep}{1mm}%
\setlength{\partopsep}{0mm}%
\setlength{\parskip}{0mm}%
\setlength{\parsep}{0mm}%
\setlength{\itemsep}{2ex}%
}}%
{\end{list}}

\newenvironment{symboldesc}[1][$A,B,C$]{%
\begin{list}{\bf ???}{%
\settowidth{\labelwidth}{#1}%
\settowidth{\leftmargin}{\qquad #1}%
\setlength{\topsep}{1mm}%
\setlength{\partopsep}{0mm}%
\setlength{\parskip}{0mm}%
\setlength{\parsep}{0mm}%
\setlength{\itemsep}{2ex}%
}}%
{\end{list}}

\newenvironment{compactsymboldesc}[1][$A,B,C$]{%
\begin{list}{\bf ???}{%
\settowidth{\labelwidth}{#1}
\settowidth{\leftmargin}{\qquad #1}%
\setlength{\topsep}{1mm}%
\setlength{\partopsep}{0mm}%
\setlength{\parskip}{0mm}%
\setlength{\parsep}{0mm}%
\setlength{\itemsep}{0ex}%
}}%
{\end{list}}

\newenvironment{dense_itemize}{%
\begin{list}{$\triangleright$}%
{\setlength{\topsep}{1mm}%
\setlength{\partopsep}{0mm}%
\setlength{\parskip}{0mm}%
\setlength{\parsep}{0mm}%
\setlength{\itemsep}{0mm}%
\setlength{\labelwidth}{4mm}%
\setlength{\leftmargin}{0mm}%
\addtolength{\leftmargin}{\labelwidth}%
\addtolength{\leftmargin}{\labelsep}%
\setlength{\itemindent}{0mm}}}%
{\end{list}}

\newenvironment{steps}[1][1]{%
\begin{list}{\bf(\arabic{enumi})}%
{\usecounter{enumi}%
\setlength{\topsep}{2mm}%
\setlength{\partopsep}{2mm}%
\setlength{\parskip}{0mm}%
\setlength{\parsep}{0mm}%
\setlength{\itemsep}{2mm}%
\settowidth{\labelwidth}{(iiiv)}%
\setlength{\leftmargin}{0mm}%
\addtolength{\leftmargin}{\labelwidth}%
\addtolength{\leftmargin}{\labelsep}%
\setlength{\itemindent}{0mm}}
\setcounter{enumi}{#1}\addtocounter{enumi}{-1}}%
{\end{list}}

\newenvironment{compactsteps}[1][1]{%
\begin{list}{\bf(\arabic{enumi})}%
{\usecounter{enumi}%
\setlength{\topsep}{2mm}%
\setlength{\partopsep}{2mm}%
\setlength{\parskip}{0mm}%
\setlength{\parsep}{0mm}%
\setlength{\itemsep}{0mm}%
\settowidth{\labelwidth}{(iiiv)}%
\setlength{\leftmargin}{0mm}%
\addtolength{\leftmargin}{\labelwidth}%
\addtolength{\leftmargin}{\labelsep}%
\setlength{\itemindent}{0mm}}
\setcounter{enumi}{#1}\addtocounter{enumi}{-1}}%
{\end{list}}

\iftrue
\newcounter{precondno}
\newcounter{globalprecondno}
\newenvironment{list_of_preconds}{%
\begin{list}{\quad(P\arabic{precondno})}%
{\usecounter{precondno}%
\setlength{\topsep}{1mm}%
\setlength{\partopsep}{0mm}%
\setlength{\parskip}{0mm}%
\setlength{\parsep}{0mm}%
\setlength{\itemsep}{0mm}%
\settowidth{\labelwidth}{\quad (P12)}%
\setlength{\leftmargin}{0mm}%
\addtolength{\leftmargin}{\labelwidth}%
\addtolength{\leftmargin}{\labelsep}%
\setlength{\itemindent}{0mm}}
\let\orgitem\item
\renewcommand{\item}{\stepcounter{globalprecondno}\orgitem}
\setcounter{precondno}{\value{globalprecondno}}
}%
{\end{list}}
\fi

\newcommand{\subfig}[2]{
\begin{tabular}{@{}c@{}}{{#2}}\\{\footnotesize #1}\end{tabular}}
       
\renewcommand{\index}[1]{}

\author{Alexander Schwartz%
        \thanks{Partially supported by the DFG Research Center 
                ``Mathematics for Key Technologies'' (FZT 86) in Berlin
 and by the German Israeli Foundation (G.I.F.).}
 \qquad\setcounter{footnote}{6}
        G\"unter M.~Ziegler%
        \thanks{Partially 
          supported by Deutsche Forschungs-Gemeinschaft, via the
  DFG Research Center ``Mathematics in the Key Technologies'' (FZT86),
  the Research Group ``Algorithms, Structure, Randomness'' (Project ZI 475/3),
  and a Leibniz grant (ZI 475/4),
     and by the German Israeli Foundation (G.I.F.)}\\[2ex]  
        \begin{minipage}{10cm}\center\footnotesize
            TU Berlin\\
            Inst.\ Mathematics, MA 6-2\\
            D-10623 Berlin, Germany\\
            \url{{schwartz,ziegler}@math.tu-berlin.de}
       \end{minipage}\\}
\title{\textbf{Construction techniques for cubical complexes,\\
    odd cubical 4-polytopes,\\ and prescribed dual manifolds}}
  \date{January 2, 2004} 

\begin{document}
\maketitle

\begin{abstract}\noindent
  We provide a number of new construction techniques for cubical
  complexes and cubical polytopes, and thus for cubifications 
  (hexahedral mesh generation).
  As an application we obtain an instance 
  of a cubical $4$-polytope that has a non-orientable dual manifold (a
  Klein bottle). This confirms an existence conjecture of Hetyei (1995).
  
  More systematically, we prove that every normal crossing codimension
  one immersion of a compact $2$-manifold into~$\R^3$ is PL-equivalent
  to a dual manifold immersion of a cubical $4$-polytope.  As an
  instance we obtain a cubical $4$-polytope with a cubation of Boy's
  surface as a dual manifold immersion, and with an odd number of
  facets.  Our explicit example has $17\,718$ vertices and $16\,533$ facets.
  Thus we get a parity changing
  operation for $3$-dimensional cubical complexes (hexa meshes); 
  this solves problems of Eppstein, Thurston, and others. 
\end{abstract}

\begin{center}
  \parbox{.8\textwidth}{
   \begin{compactdesc}
     \item[\textbf{Keywords:}] Cubical complexes, cubical polytopes, 
       regular subdivisions, normal crossing codimension one
       PL immersions, construction techniques, cubical meshes, Boy's surface\\
 
       \item[\textbf{MSC 2000 Subject Classification:}]  
                    52B12, 52B11, 52B05, 57Q05\\
   \end{compactdesc}}
\end{center}
\section{Introduction}

A $d$-polytope is \emph{cubical} if all its proper faces are
combinatorial cubes, that is, if each $k$-face of the polytope,
$k\in\{0,\ldots,d-1\}$ is combinatorially equivalent to the
$k$-dimensional standard cube.

It has been observed by Stanley, MacPherson, and others
(cf.~\cite{BabsonChan3} \cite{Jock}) that every cubical $d$-polytope
$P$ determines a PL immersion of an abstract cubical $(d-2)$-manifold
into the polytope boundary $\boundaryOP{P}\cong S^{d-1}$.  The
immersed manifold is orientable if and only if the $2$-skeleton of the
cubical $d$-polytope ($d\ge3$) is ``edge orientable'' in the sense of
Hetyei, who conjectured that there are cubical $4$-polytopes that are
not edge-orientable~\cite[Conj.~2]{hetyei95:_stanl}.

In the more general setting of cubical PL $(d-1)$-spheres, Babson and
Chan \cite{BabsonChan3} have observed that \emph{every} type of normal
crossing PL immersion of a $(d-2)$-manifold into a $(d-1)$-sphere
appears among the dual manifolds of some cubical PL $(d-1)$-sphere.

No similarly general result is available for cubical polytopes.  The
reason for this may be traced/blamed to a lack of flexible
construction techniques for cubical polytopes, and more generally, for
cubical complexes (such as the ``hexahedral meshes'' that are of great
interest in CAD and in Numerical Analysis).
  
In this paper, we develop a number of new and improved construction
techniques for cubical polytopes.  We try to demonstrate that it
always pays off to carry along convex lifting functions of high
symmetry.  The most complicated and subtle element of our
constructions is the ``generalized regular Hexhoop'' of
Section~\ref{subsec:generalized_hexhoop}, which yields a cubification
of a $d$-polytope with a hyperplane of symmetry, where a (suitable)
lifting function may be specified on the boundary.
  
Our work is extended by the first author in \cite{Schwartz2}, where
additional construction techniques for \emph{cubifications}
(i.\,e.~cubical subdivisions of $d$-polytopes with prescribed boundary
subdivisions) are discussed.

Using the constructions developed here, we achieve the following
constructions and results:
\begin{compactitem}[$\bullet$]\itemsep=1mm
\item A rather simple construction yields a cubical $4$-polytope (with
  $72$ vertices and $62$ facets) for which the immersed dual
  $2$-manifold is not orientable: One of its components is a Klein
  bottle. Apparently this is the first example of a cubical polytope
  with a non-orientable dual manifold.  Its existence confirms a
  conjecture of Hetyei (Section~\ref{sec:klein}).
\item More generally, all PL-types of normal crossing immersions of
  $2$-manifolds appear as dual manifolds in the boundary complexes of
  cubical $4$-polytopes (Section~\ref{sec:immersions}).  In the case
  of non-orientable $2$-manifolds of odd genus, this yields cubical
  $4$-polytopes with an odd number of facets.  From this, we also
  obtain a complete characterization of the lattice of $f$-vectors of
  cubical $4$-polytopes (Section~\ref{sec:consequences}).
\item In particular, we construct an explicit example with $17\,718$
  vertices and $16\,533$ facets of a cubical $4$-polytope which has a
  cubation of Boy's surface (projective plane with exactly one triple
  point) as a dual manifold immersion (Section~\ref{sec:boy}).
\item Via Schlegel diagrams, this implies that every $3$-cube has a
  cubical subdivision into an even number of cubes that does not
  subdivide the boundary complex.  Thus for every cubification of a
  $3$-dimensional domain there is also a cubification of the opposite
  parity (Section~\ref{sec:hexameshing}). This answers questions by
  Bern, Eppstein, Erickson, and Thurston \cite{BernEppsteinErickson}
  \cite{Eppstein1999} \cite{Thurston1993}.
\end{compactitem}


\section{Basics}

For the following we assume that the readers are familiar with the
basic combinatorics and geometry of convex polytopes. In particular,
we will be dealing with cubical polytopes (see Gr\"unbaum
\cite[Sect.~4.6]{Gr1-2}), polytopal (e.\,g.\ cubical) complexes,
regular subdivisions (see Ziegler \cite[Sect.~5.1]{Z35}), and Schlegel
diagrams \cite[Sect.~3.3]{Gr1-2} \cite[Sect.~5.2]{Z35}.  For cell
complexes, barycentric subdivision and related notions we refer to
Munkres~\cite{Munk}.  Suitable references for the basic concepts about
PL manifolds, embeddings and (normal crossing) immersions include
Hudson \cite{Hudson} and Rourke \& Sanderson \cite{RourkeSanderson}.
\subsection{Almost cubical polytopes}

All proper faces of a cubical $d$-polytope have to be combinatorial
cubes. We define an \emph{almost cubical}\index{almost
  cubical}\index{polytope!almost cubical} $d$-polytope as a pair
$(P,F)$, where $F$ is a specified facet of~$P$ such that all facets
of~$P$ other than $F$ are required to be combinatorial cubes. Thus,
$F$ need not be a cube, but it will be cubical.

By $\ComplexOf{P}$ we denote the polytopal complex given by a polytope
$P$ and all its faces.  By $\ComplexOf{\boundaryOP P}$ we denote the
\emph{boundary complex}\index{polytope!boundary
  complex}\index{boundary complex} of $P$, consisting of all proper
faces of~$P$. If $P$ is a cubical polytope, then
$\ComplexOf{\boundaryOP P}$ is a cubical complex.  If $(P,F)$ is
almost cubical, then the \emph{Schlegel
  complex}\index{polytope!Schlegel complex} $\ComplexOf{\boundaryOP
  P}{\setminus}\{F\}$ is a cubical complex that is combinatorially
isomorphic to the Schlegel diagram $\Schlegel(P,F)$ of $P$ based
on~$F$.

\subsection{Cubifications}\label{subsec:cubifications}

A \emph{cubification}\index{cubification} of a cubical PL
$(d-1)$-sphere $\mathcal{S}^{d-1}$ is a cubical $d$-ball
$\mathcal{B}^d$ with boundary $\mathcal{S}^{d-1}$.  A double counting
argument shows that every cubical $(d-1)$-sphere that admits a
cubification has an even number of facets.  Whether this condition is
sufficient is a challenging open problem, even for $d=3$
(compare~\cite{BernEppsteinErickson}, \cite{Eppstein1999}).

\subsection{Dual manifolds}\label{subsec:dual_mflds}\label{subsec:NCP}

For every (pure) cubical $d$-dimensional complex $\mathcal{C}$, $d>1$,
the \emph{derivative complex} is an abstract cubical cell
$(d-1)$-dimensional complex~$\DerivativeComplex{\mathcal{C}}$ whose
vertices may be identified with the edge midpoints of the complex,
while the facets ``separate the opposite facets of a facet
of~$\mathcal{C}$,'' that is, they correspond to pairs $(F,[e])$, where
$F$ is a facet of~$\mathcal{C}$ and $[e]$ denotes a ``parallel class''
of edges of $F$. This is a cell complex with $f_1(\mathcal{C})$
vertices and $(d-1)f_{d-1}(\mathcal{C})$ cubical facets of dimension
$d-1$, $d-1$ of them for each facet of $\mathcal{C}$.  Hence the
derivative complex $\DerivativeComplex{\mathcal{C}}$ is pure
$(d-1)$-dimensional.  See Babson \& Chan \cite[Sect.~4]{BabsonChan3}.
 
  \begin{figure}[H]\centering

    \subfigure[$\mathcal{C}$]{\includegraphics[width=.25\textwidth]{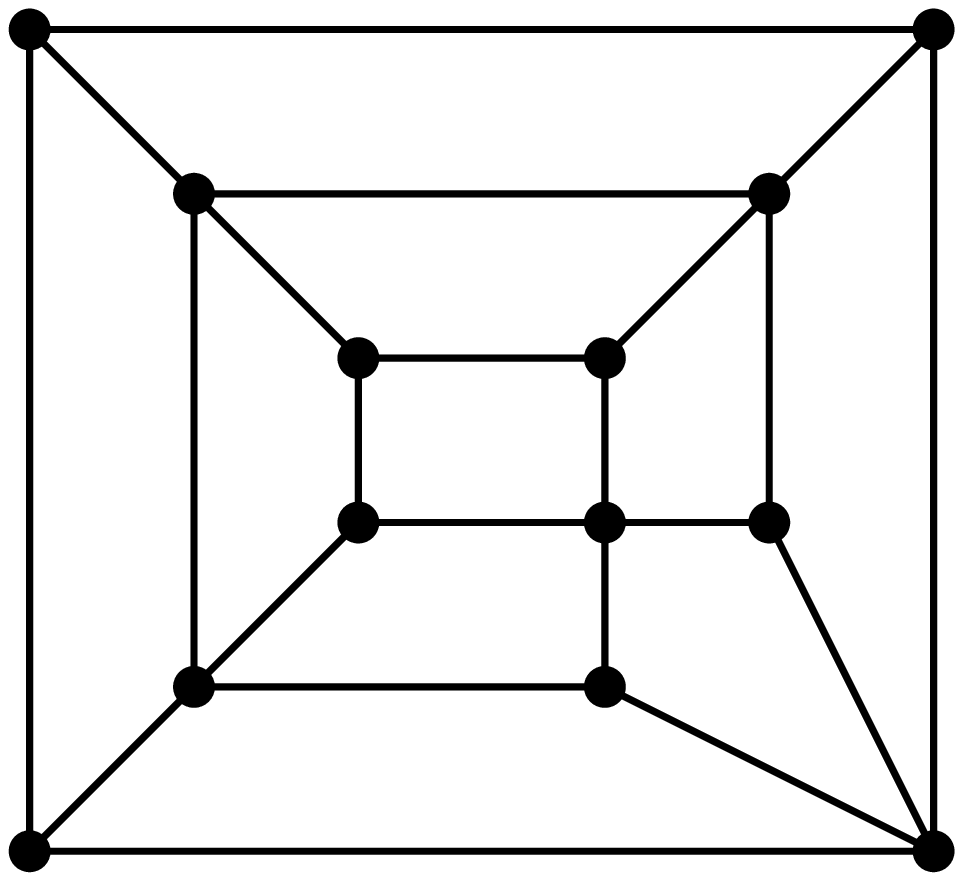}}\qquad
      \subfigure[$\mathcal{C}$ and $\DerivativeComplex{\mathcal{C}}$]{%
        \includegraphics[width=.25\textwidth]{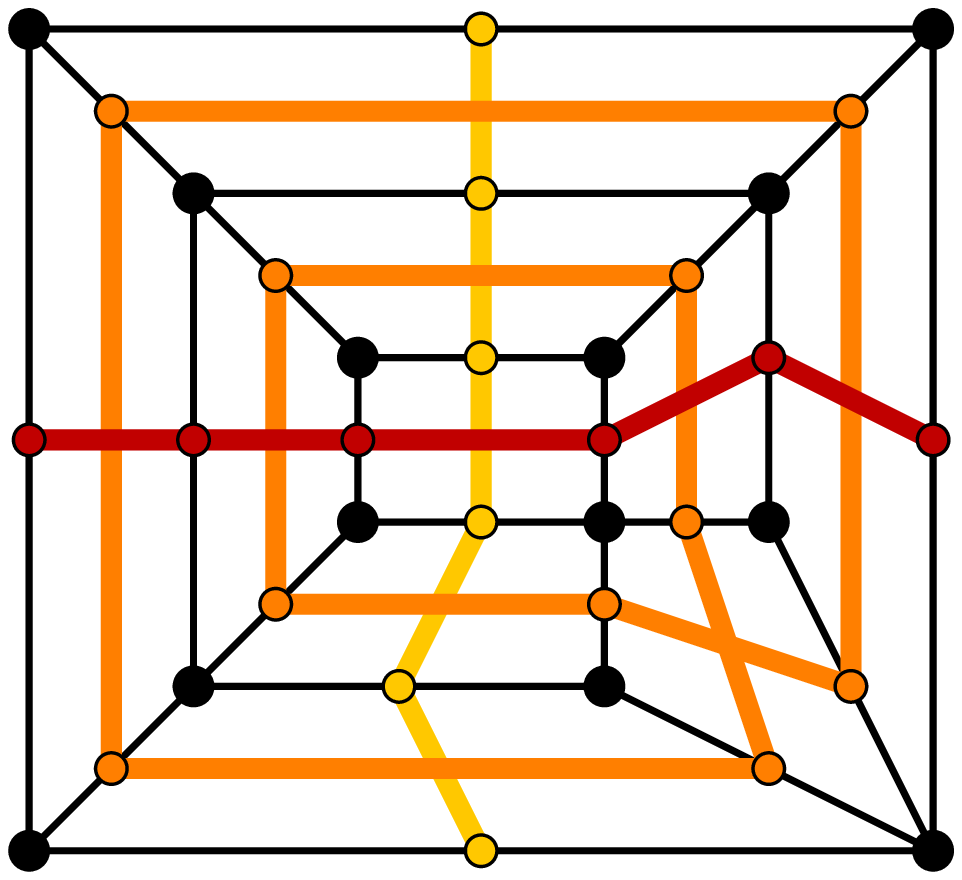}}\qquad
      \subfigure[$\DerivativeComplex{\mathcal{C}}$]{%
         \includegraphics[width=.25\textwidth]{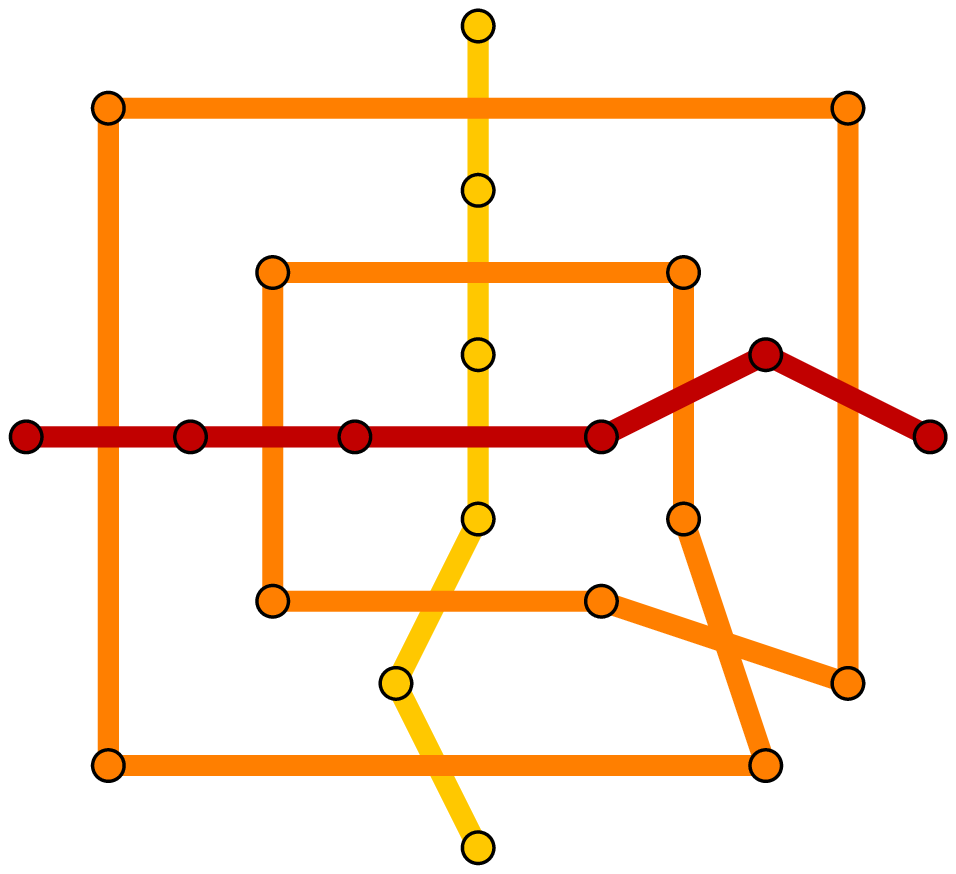}}
\vspace{-4mm}
    \caption[Derivative complex]{The derivative complex of a cubical $2$-complex $\mathcal{C}$.}
    \label{fig:derivative_complex}
  \end{figure}
        
  In the case of cubical PL spheres (for instance boundary complexes
  of cubical polytopes), or cubical PL balls, the derivative complex
  is a (not necessarily connected) manifold, and we call each
  connected component of the derivative
  complex~$\DerivativeComplex{P}$ of a cubical complex~$\mathcal{C}$ a
  \emph{dual manifold} of~$\mathcal{C}$.  If the cubical
  complex~$\mathcal{C}$ is a sphere then the dual manifolds of
  $\mathcal{C}$ are manifolds without boundary.
  If $\mathcal{C}$ is a ball, then some (possibly
  all) dual manifolds have non-empty boundary components, namely the
  dual manifolds of $\boundaryOP\mathcal{C}$.\\

  The derivative complex, and thus each dual manifold, comes with a
  canonical immersion into the boundary of~$P$. More precisely, the
  barycentric subdivision\index{barycentric subdivision!simplicial} of
  $\DerivativeComplex{P}$ has a simplicial map to the barycentric
  subdivision of the boundary complex $\boundaryOP P$, which is a
  codimension one normal crossing immersion into the simplicial sphere
  $\simplicialBarycentricSubdiv{\ComplexOf{\boundaryOP P}}$.
  (\emph{Normal crossing} means that each multiple-intersection point
  is of degree $k\leq d$ and there is a neighborhood of each multiple
  intersection point that is PL isomorphic to (a neighborhood of) a
  point which is contained in $k$ pairwise perpendicular
  hyperplanes.)

   Restricted to a dual manifold, this immersion may be an embedding
   or not.
 \begin{Figure}[ht]\centering
    {\includegraphics[width=.27\textwidth]{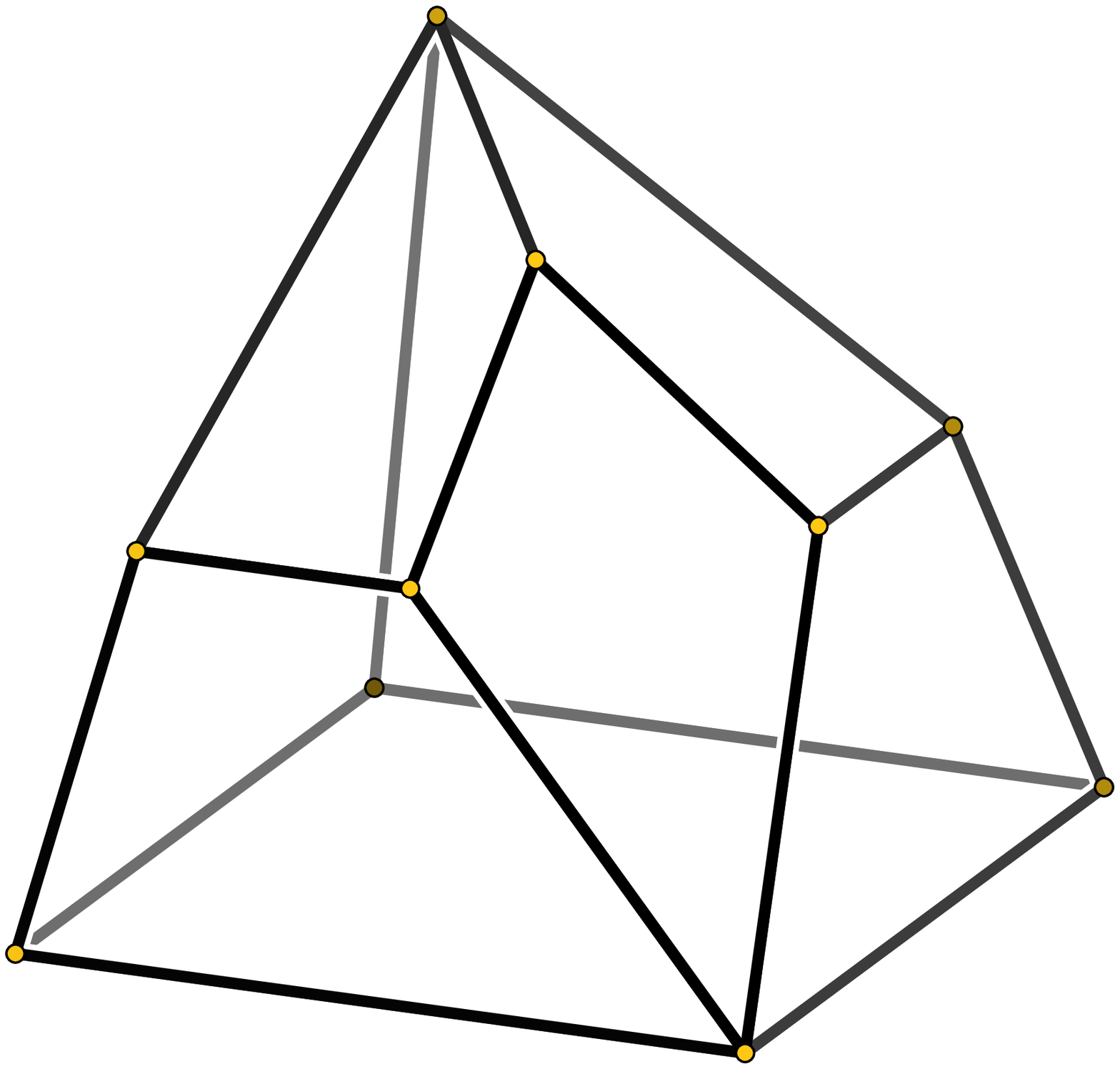}}\qquad\qquad
    {\includegraphics[width=.27\textwidth]{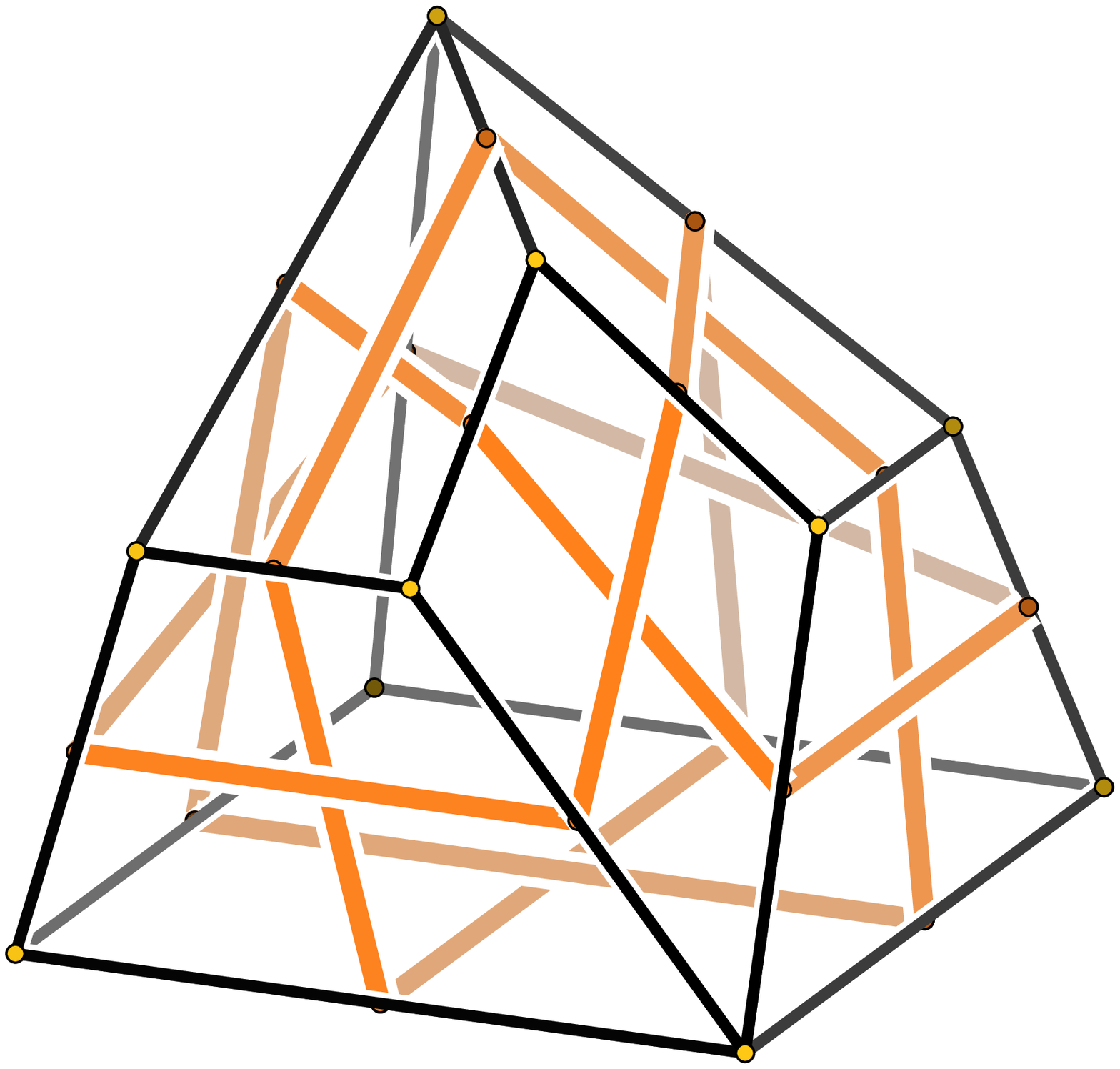}}
  \caption[The cubical octahedron]{
    The cubical octahedron $O_8$ (the only combinatorial type of a
    cubical $3$-polytope with 8 facets), and its single immersed dual
    manifold.}
  \label{fig:cubical_octahedron}
\end{Figure}
  
   In the case of cubical $3$-polytopes, the derivative complex may
   consist of one or many $1$-spheres. For example, for the $3$-cube
   it consists of three $1$-spheres, while for the ``cubical
   octahedron'' $O_8$ displayed in Figure~\ref{fig:cubical_octahedron}
   the dual manifold is a single immersed~$S^1$ (with $8$ double
   points).

In the case of $4$-polytopes, the dual manifolds are surfaces (compact
$2$-manifolds without boundary).
As an example, we here display a Schlegel diagram of a ``neighborly
cubical'' $4$-polytope (with the graph of the $5$-cube), with
$f$-vector $(32,80,96,48)$.

\begin{Figure}\centering
  \includegraphics[width=.35\textwidth]{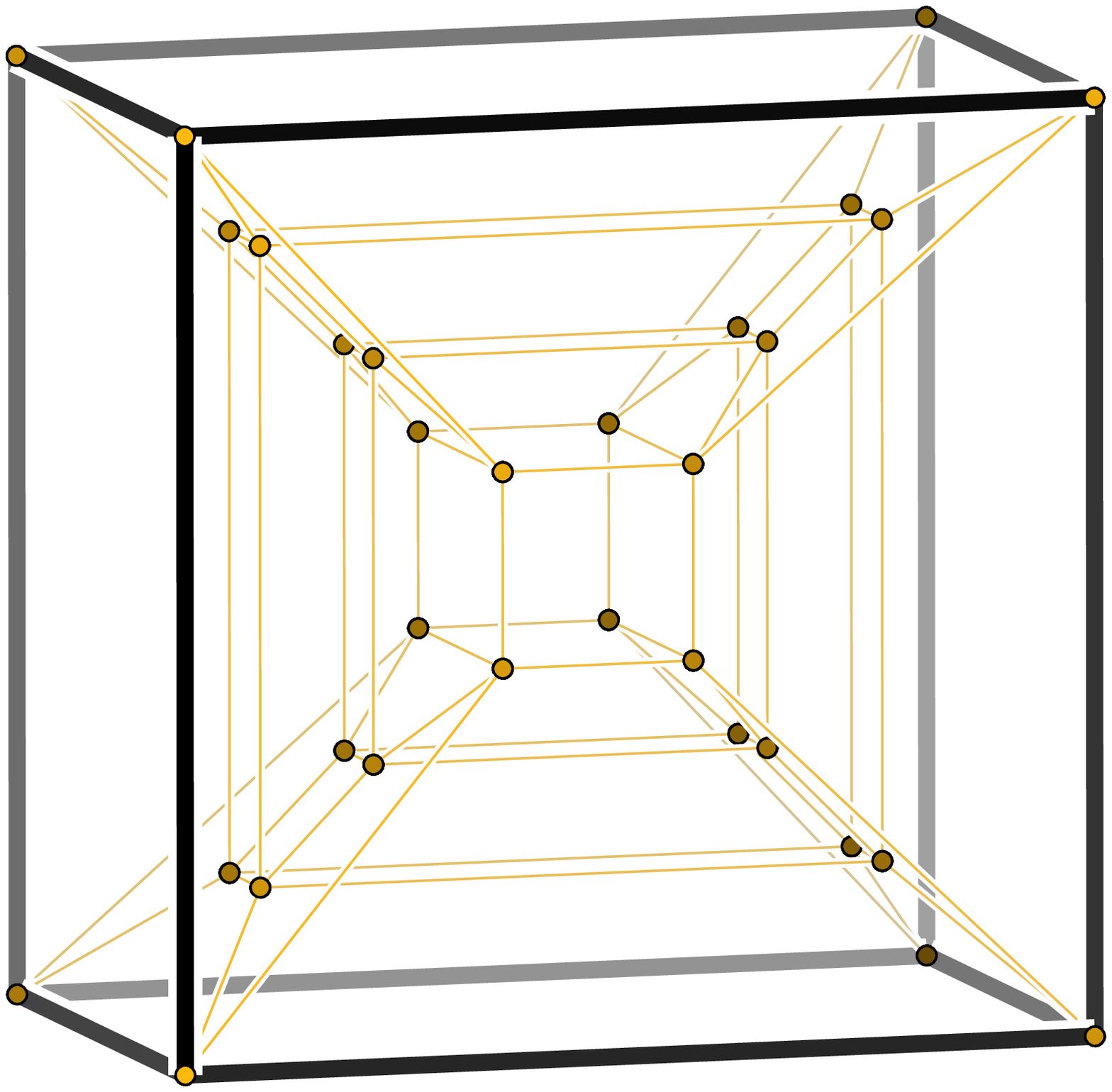}\qquad
  \includegraphics[width=.35\textwidth]{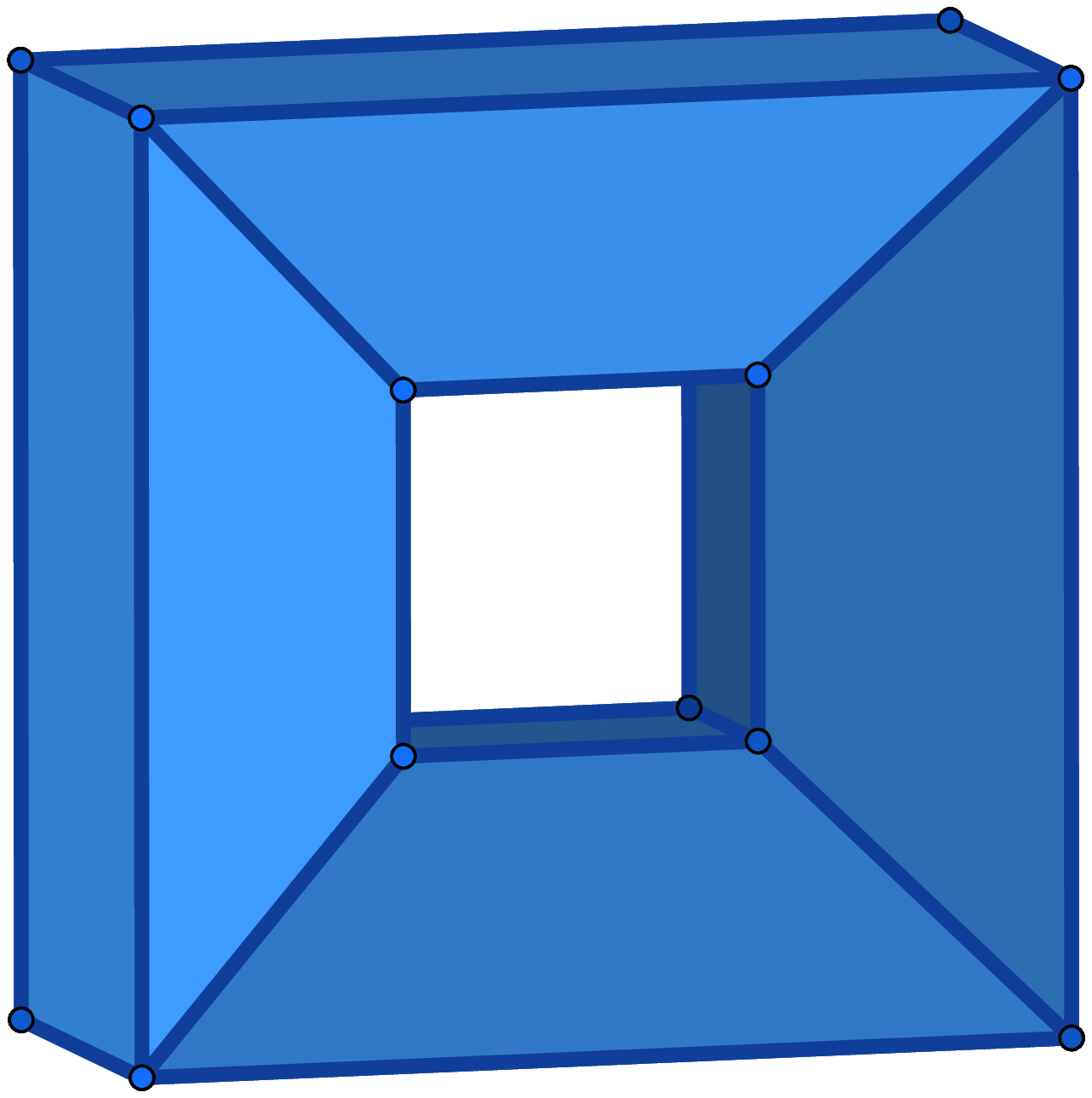}
  \label{fig:ncp} 
  \caption[The neighborly cubical $4$-polytope $C^5_4$]{%
    A Schlegel diagram of the ``neighborly cubical'' $4$-polytope
    $C^5_4$ with the graph of the $5$-cube, and its dual torus. All
    other dual manifolds are embedded $2$-spheres.}
\end{Figure}

  According to Joswig \& Ziegler \cite{JoswigZiegler} this may be
  constructed as
  \[
      C^5_4\ \ :=\ \ \conv((Q\times 2Q)\ \cup\ (2Q\times Q)), \qquad\text{where }Q=[-1,+1]^2.
  \] 
  Here the dual manifolds are four embedded cubical $2$-spheres $S^2$
  with $f$-vector $(16,28,14)$ --- of two different combinatorial
  types --- and one embedded torus $T$ with $f$-vector $(16,32,16)$.

\subsection{Orientability}
  
  Let $P$ be a cubical $d$-polytope $(d\ge3$). The immersed dual
  manifolds in its boundary cross the edges of
  the polytope transversally.
  
  Thus we find that orientability of the dual manifolds is equivalent
  to the possibility to give \emph{consistent} edge orientations to
  the edges of the $P$, that is, in each $2$-face of~$P$ opposite edges
  should get parallel (rather than antiparallel) orientations; compare Hetyei
  \cite{hetyei95:_stanl}. Figure~\ref{fig:edge_orientation} shows such
  an edge orientation for a cubical $3$-polytope (whose derivative
  complex consists of three circles, so it has $8$ consistent edge
  orientations in total).

\begin{Figure}[ht]\centering
    \includegraphics[width=40mm]{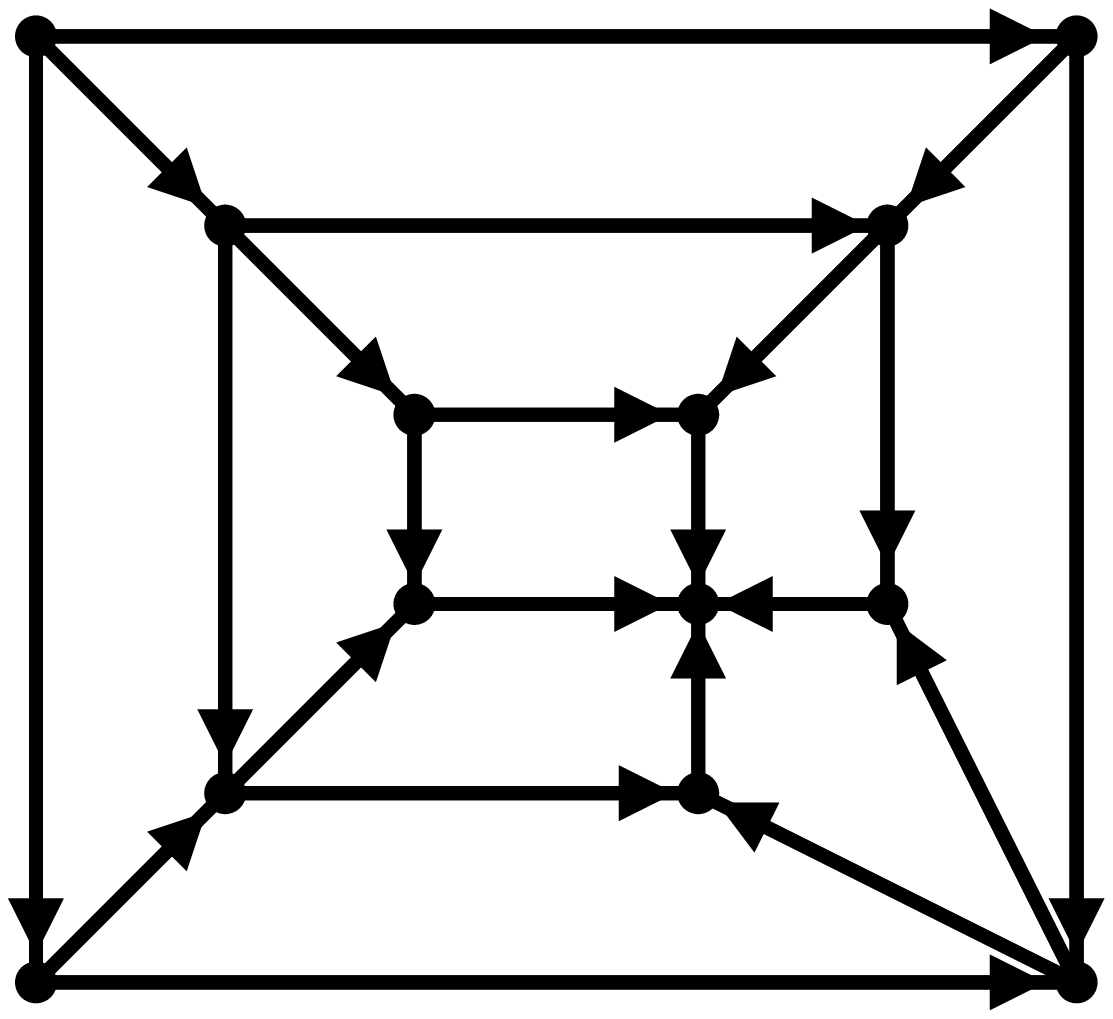}\qquad
    \includegraphics[width=40mm]{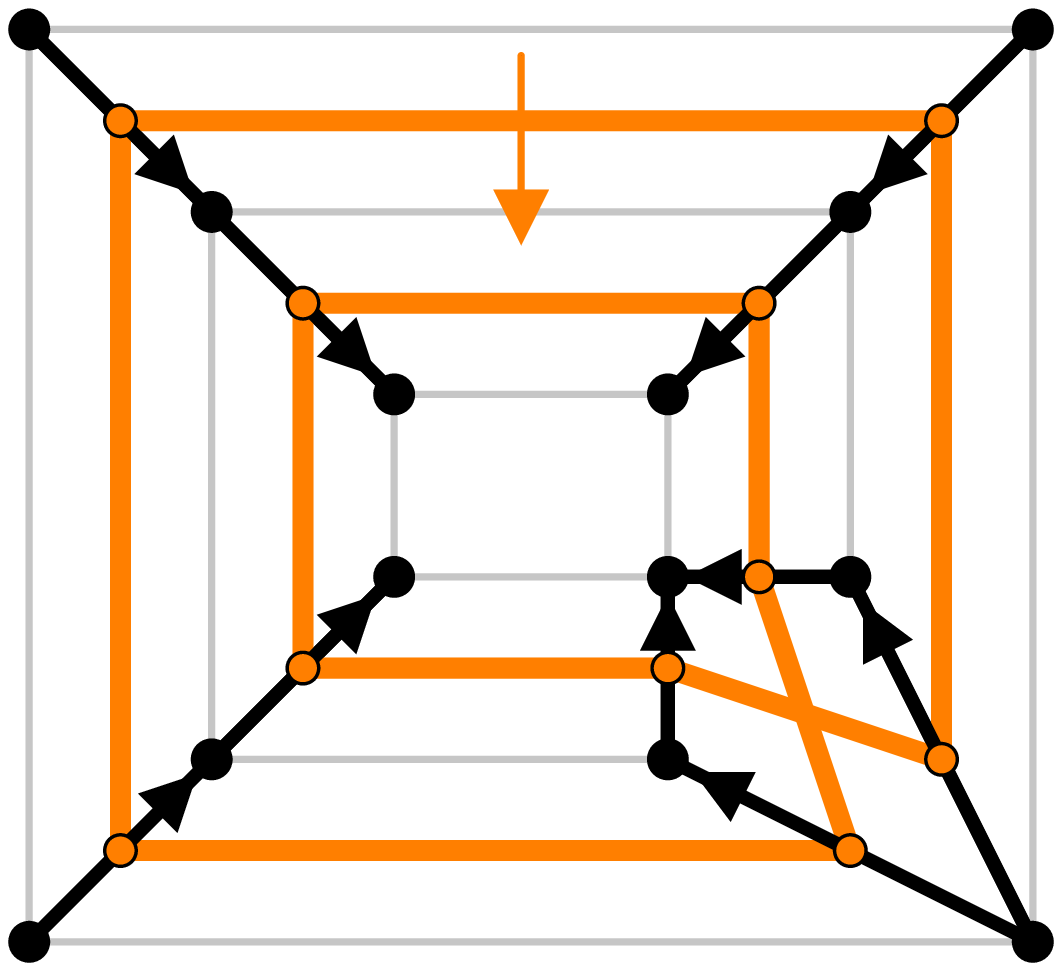}
    \caption{Edge-orientation of (the Schlegel diagram of) a 
             cubical $3$-polytope. 
             The edges marked on the right
             must be oriented consistently.}
     \label{fig:edge_orientation}
\end{Figure}

One can attempt to obtain such edge orientations by moving to from
edge to edge across $2$-faces.  The obstruction to this arises if on a
path moving from edge to edge across quadrilateral $2$-faces we return
to an already visited edge, with reversed orientation, that is, if we
close a \emph{cubical M\"obius strip} with parallel inner edges, as
displayed in the figure.  (Such an immersion is not necessarily 
embedded, that is, some $2$-face may be used twice for the M\"obius
strip.)

   \begin{Figure}[ht]\centering
    \begin{psfrags}
     \psfrag{->}{\rput(1.5mm,0){$\longrightarrow$}}
     \psfrag{v_0}{$\boldsymbol{v}_0$}
     \psfrag{v_1}{$\boldsymbol{v}_1$}
     \psfrag{v_2}{$\boldsymbol{v}_2$}
     \psfrag{v_{l-1}}{$\boldsymbol{v}_{\ell-1}$}
     \psfrag{v_l=w_0}{$\boldsymbol{v}_{\ell}=\boldsymbol{w}_0$}
     \psfrag{w_0}{$\boldsymbol{w}_0$}
     \psfrag{w_1}{$\boldsymbol{w}_1$}
     \psfrag{w_2}{$\boldsymbol{w}_2$}
     \psfrag{w_{l-1}}{$\boldsymbol{w}_{\ell-1}$}
     \psfrag{w_l=v_0}{$\boldsymbol{w}_{\ell}=\boldsymbol{v}_0$}
     \includegraphics[width=80mm]{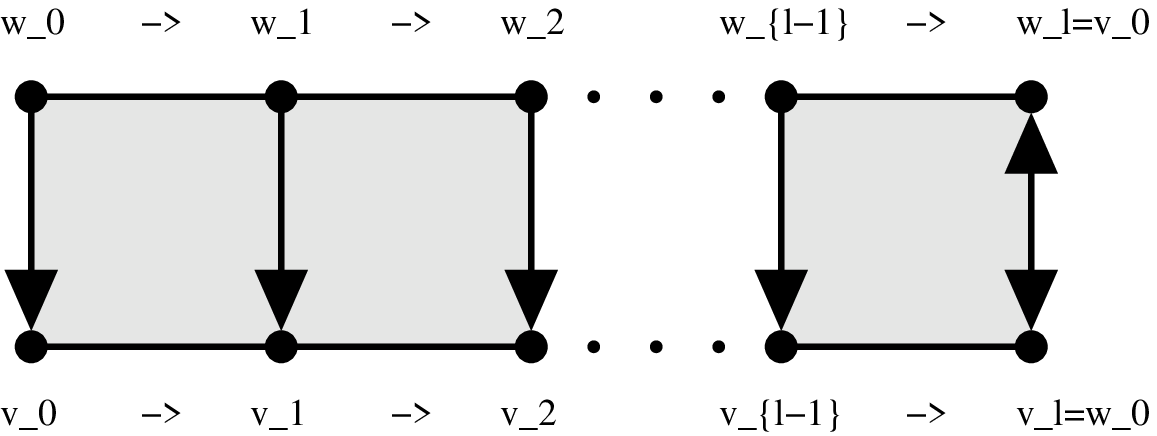}
     \end{psfrags}
     \caption{A cubical M\"obius strip with parallel inner edges.}
     \label{fig:cubical_Moebius_strip}
   \end{Figure}

\begin{proposition}\label{prop:NEO}
   For every cubical $d$-polytope $(d\ge3)$, the 
   following are equivalent:
   \begin{compactitem}[$\bullet$]
        \item All dual manifolds of~$P$ are orientable.

        \item The $2$-skeleton of~$P$ has a consistent edge orientation.

        \item The $2$-skeleton of~$P$ contains no immersion of
           a {cubical M\"obius strip} with parallel inner edges.
    \end{compactitem}
\end{proposition}

\subsection{From PL immersions to cubical PL spheres}\label{subsect:bc}

The emphasis in this paper is on cubical convex $d$-polytopes.  In the
more general setting of cubical PL ($d-1$)-spheres, one has more
flexible tools available. In this setting, Babson \&
Chan~\cite{BabsonChan3} proved that ``all PL codimension~$1$ normal
crossing immersions appear.''  The following sketch is meant to
explain the Babson-Chan theorem geometrically (it is presented in a
combinatorial framework and terminology in~\cite{BabsonChan3}), and to
briefly indicate which parts of their construction are available in
the polytope world.

  \begin{construction}{Babson-Chan~\cite{BabsonChan3}}
    
    \begin{inputoutput}
       \item[Input:] A normal crossing immersion
         $j:\mathcal{M}^{d-2}\rightarrow \mathcal{S}^{d-1}$ of a
         triangulated PL manifold $\mathcal{M}^{d-2}$ of dimension
         $d-2$ into a PL simplicial ($d-1$)-sphere.
        \item[Output:]  A cubical PL ($d-1$)-sphere with a dual
          manifold immersion PL-equivalent to $j$.
    \end{inputoutput}

  \begin{steps}
   \item
     Perform a barycentric subdivision\index{barycentric subdivision!simplicial}
      on $\mathcal{M}^{d-2}$ and $\mathcal{S}^{d-1}$.
     \begin{Figure}
        \includegraphics[width=4cm]{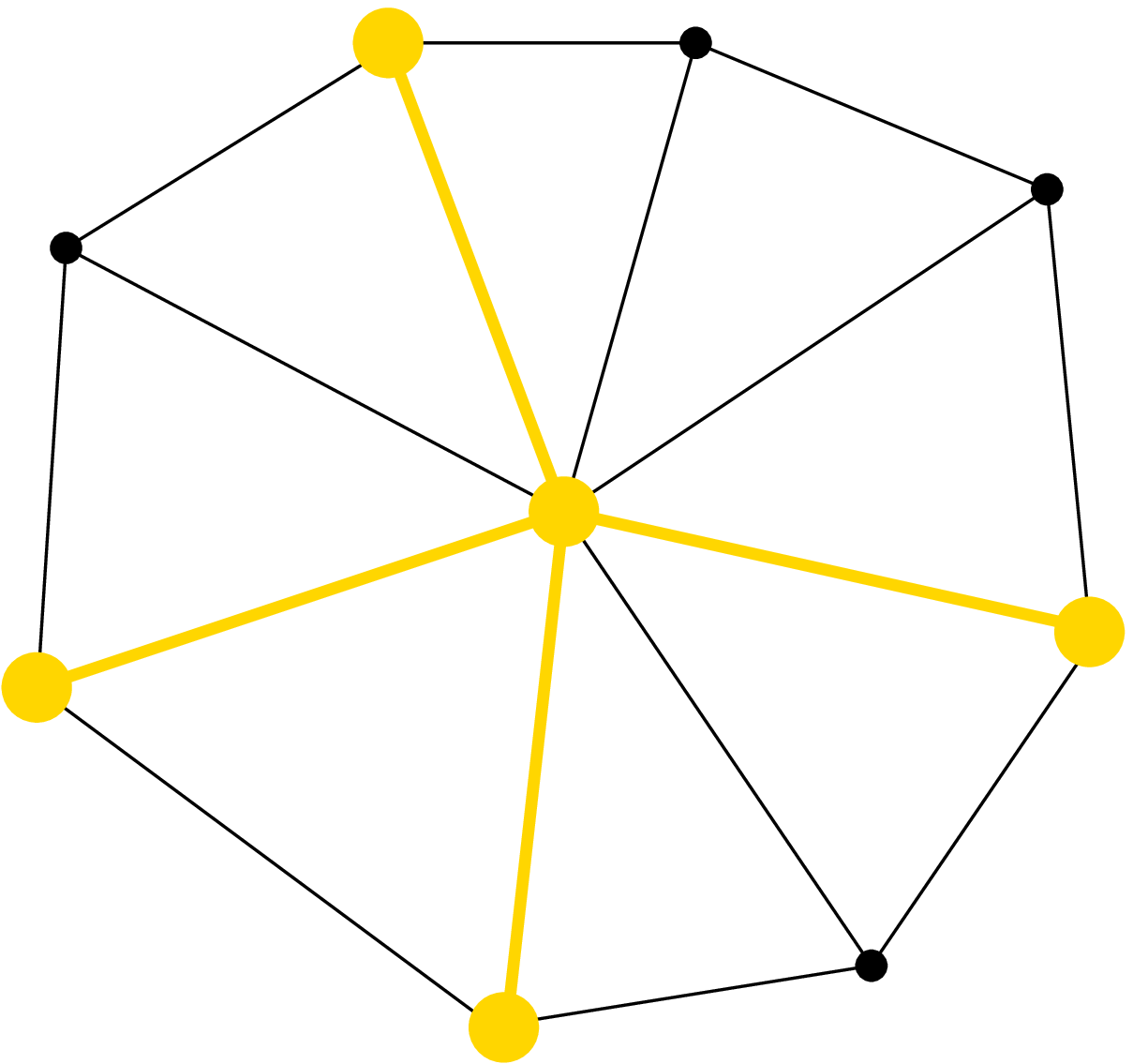}
          \quad\transformsToArrow{18mm}{10mm}\qquad
        \includegraphics[width=4cm]{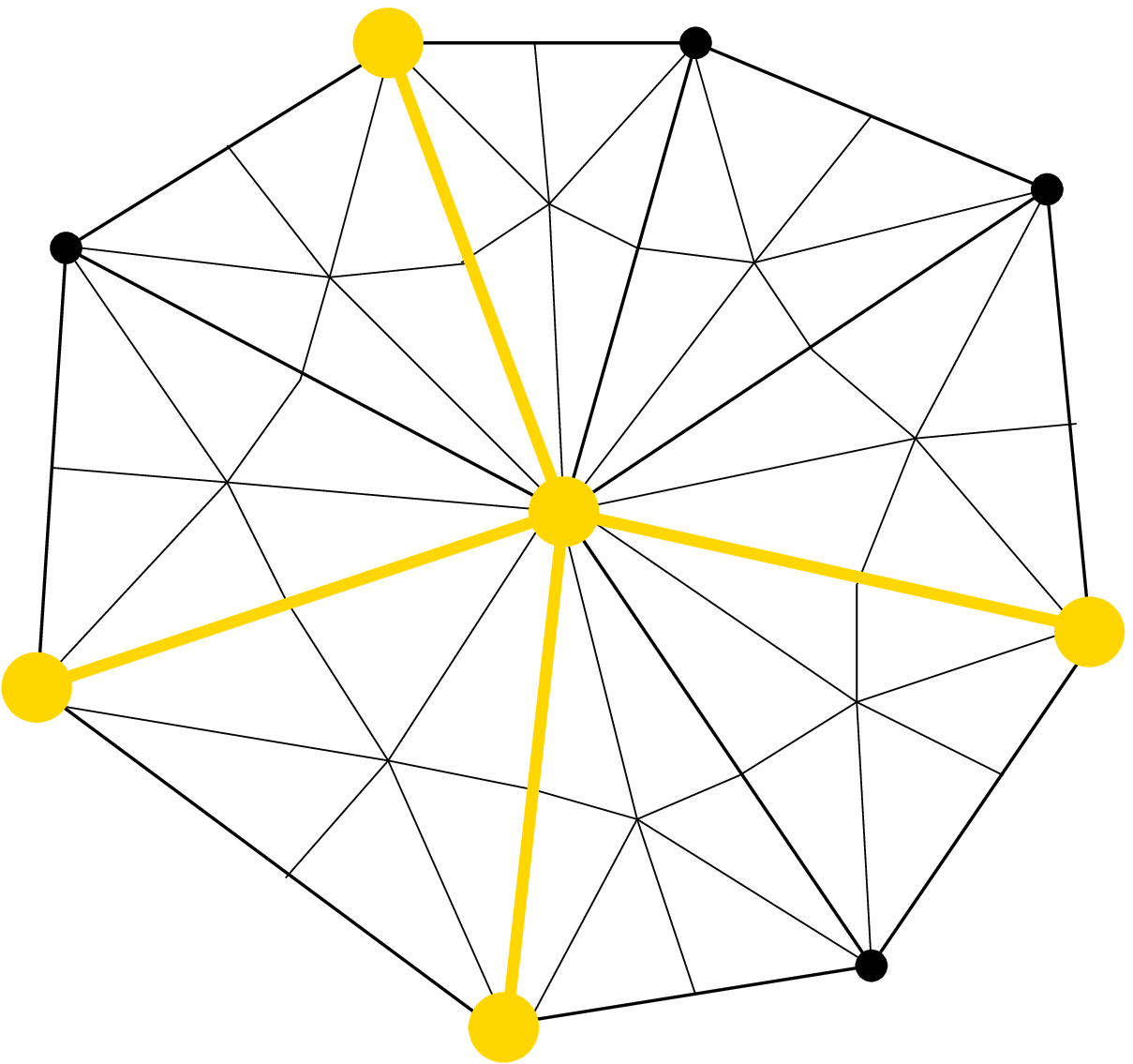}
        \setlength{\captionmargin}{2cm}
        \caption[Step 1 of the Babson-Chan construction]{
          Step 1. Performing a
          barycentric subdivision.  (We illustrate the impact of the
          construction on $2$-ball, which might be part of the
          boundary of a $2$-sphere. The immersion which is shown in
          bold has a single double-intersection point.)}
     \end{Figure}

     (Here each $i$-simplex is replaced by $(i+1)!$ new $i$-simplices,
     which is an even number for $i>0$.  This step is done only to
     ensure parity conditions on the $f$-vector, especially that the
     number of facets of the final cubical sphere is congruent to the
     Euler characteristic of~$\mathcal{M}^{d-2}$.  Barycentric subdivisions are
     easily performed in the polytopal category as well, see Ewald \&
     Shephard \cite{EwSh}.)

  \item
     Perform a ``cubical barycentric subdivision''\index{barycentric subdivision!cubical}
     on $\mathcal{M}^{d-2}$ and $\mathcal{S}^{d-1}$.
     \begin{Figure}
        \includegraphics[width=4cm]{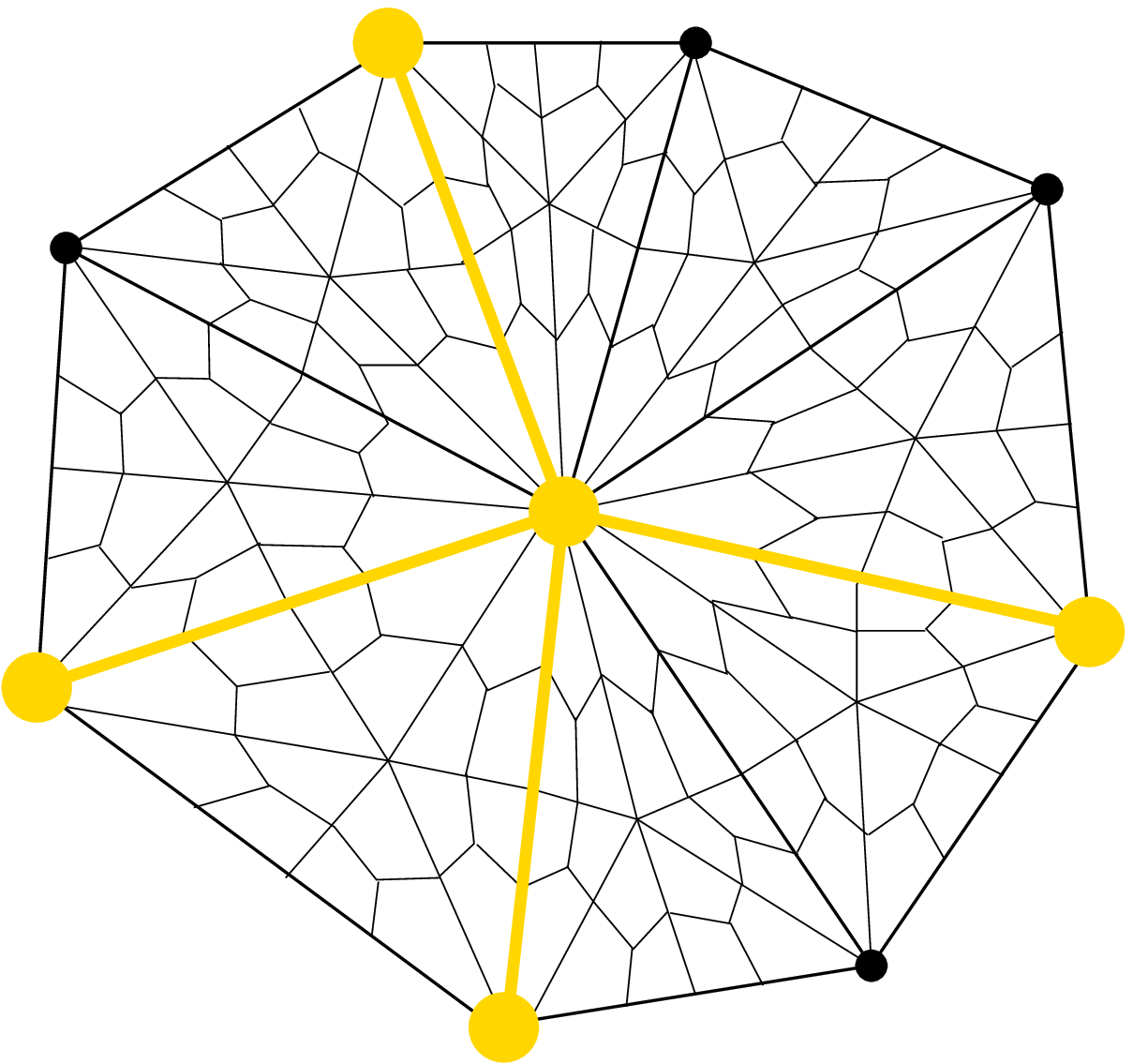}
        \setlength{\captionmargin}{2cm}
        \caption[Step 2 of the Babson-Chan construction]{
        Step 2. Performing a cubical barycentric subdivision.}
     \end{Figure}

     (This is the standard tool for passage from a simplicial complex
     to a PL-homeomorphic cubical complex; here every $i$-simplex is
     subdivided into $i+1$ different $i$-cubes. Such cubations can be
     performed in the polytopal category according to
     Shephard~\cite{shephard66:_approx}: If the starting triangulation
     of $\mathcal{S}^{d-1}$ was polytopal, the resulting cubation will be
     polytopal as well.)
  \item
  ``Thicken'' the cubical ($d-1$)-sphere along the immersed ($d-2$)-manifold,
  to obtain the cubical ($d-1$)-sphere 
   $BC(\mathcal{S}^{d-1},j(\mathcal{M}^{d-2}))$.

      (In this step, every ($d-1-i$)-cube in the $i$-fold multiple
      point locus results in a new $(d-1)$-cube. The original immersed
      manifold, in its cubified subdivided version, now appears as a
      dual manifold in the newly resulting $(d-1)$-cubes.  This last
      step is the one that seems hard to perform for polytopes in any
      non-trivial instance.)\label{pref:last-BC-step-not-polytopal}
\end{steps}\vskip-5mm
     \begin{Figure} 
         \includegraphics[width=5cm]{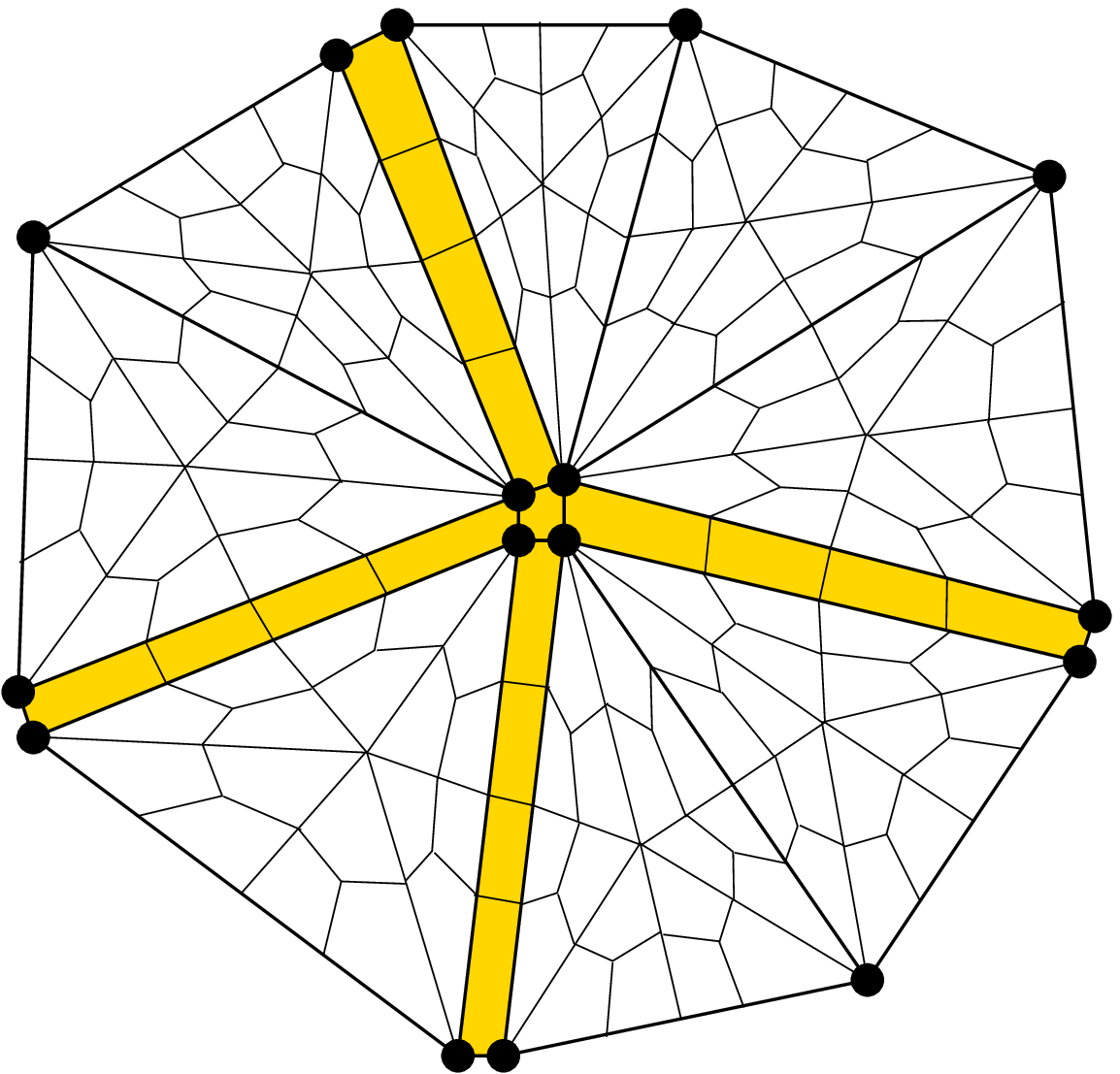}   
        \setlength{\captionmargin}{2cm}
        \caption[The outcome of the Babson-Chan construction]{
          The outcome of the Babson-Chan construction: A cubical
          sphere with a dual manifold immersion that is PL-equivalent
          to the input immersion~$j$.}
     \end{Figure}
\vskip-5mm
\end{construction}



\section{Lifting polytopal subdivisions}\label{seclifting}

\subsection{Regular balls}\label{sec:RCBalls}

In the following, the primary object we deal with is a 
\emph{regular ball}: a regular polytopal subdivision~$\mathcal{B}$ of a convex
polytope $P=\support{\mathcal{B}}$.  

\begin{definition}[regular subdivision, lifting function]
A polytopal subdivision~$\mathcal{B}$ 
is \emph{regular} (also known as
\emph{coherent} or \emph{projective}) if it admits a 
\emph{lifting function}, that is, 
a concave function $f: \support{P}\rightarrow\R$
whose domains of linearity are the facets of the subdivision.  (A
function $g:D\rightarrow\R$ is \emph{concave} if for all
$\boldsymbol{x},\boldsymbol{y}\in D$ and $0 < \lambda < 1$ we have
$g(\lambda\boldsymbol{x} + (1-\lambda)\boldsymbol{y}) \geq \lambda
g(\boldsymbol{x})+ (1-\lambda) g(\boldsymbol{y}).$)
\end{definition}

In this definition, subdivisions of the boundary are allowed, that is,
we do not necessarily require that the faces of
$P=\support{\mathcal{B}}$ are themselves faces in~$\mathcal{B}$.

In the sequel we focus on regular \emph{cubical} balls. Only in some cases we
consider regular non-cubical balls.

\begin{example}
  If~$(P,F)$ is an almost cubical polytope, then the Schlegel diagram
  based on $F$, which we denote by $\Schlegel (P,F)$, is a regular
  cubical ball (without subdivision of the boundary).
\end{example}

\begin{lemma}\label{lem:convexification}
  If~$\mathcal{B}$ is a regular cubical $d$-ball, then there is a
  regular cubical ball~$\mathcal{B'}$ without subdivision of the
  boundary, combinatorially isomorphic to~$\mathcal{B}$.
\end{lemma}

\begin{proof}
  Using a positive lifting function
  $f:\support{\mathcal{B}}\rightarrow\R$, the $d$-ball $\mathcal{B}$
  may be lifted to $\widetilde{\mathcal{B}}$ in~$\R^{d+1}$, by mapping
  each $\boldsymbol{x}\in\support{\mathcal{B}}$ to
  $(\boldsymbol{x},f(\boldsymbol{x}))\in\R^{d+1}$.

Viewed from $\boldsymbol p:=\lambda \boldsymbol e_{d+1}$ for
sufficiently large~$\lambda$, this lifted ball will appear to be
\emph{strictly convex}, that is, its boundary is a convex polytope
(rather that a boundary subdivision of a convex polytope).
Thus one may look at the polytopal complex that consists of the cones
spanned by faces of $\widetilde{\mathcal{B}}$ with apex~$\boldsymbol
p$.  This polytopal complex is regular, since it appears convex when
viewed from~$\boldsymbol p$, which yields a lifting function for the
restriction of $\widetilde{\mathcal{B}}$ to the hyperplane given by
$x_{d+1}=0$, which may be taken to be $\mathcal{B}'$.
\end{proof}
\begin{Figure}
    \psfrag{proj}{}
    \psfrag{base}{\put(-15,-40){$\mathcal{B}$}}
    \psfrag{lift}{}
    \includegraphics[width=.4\textwidth]{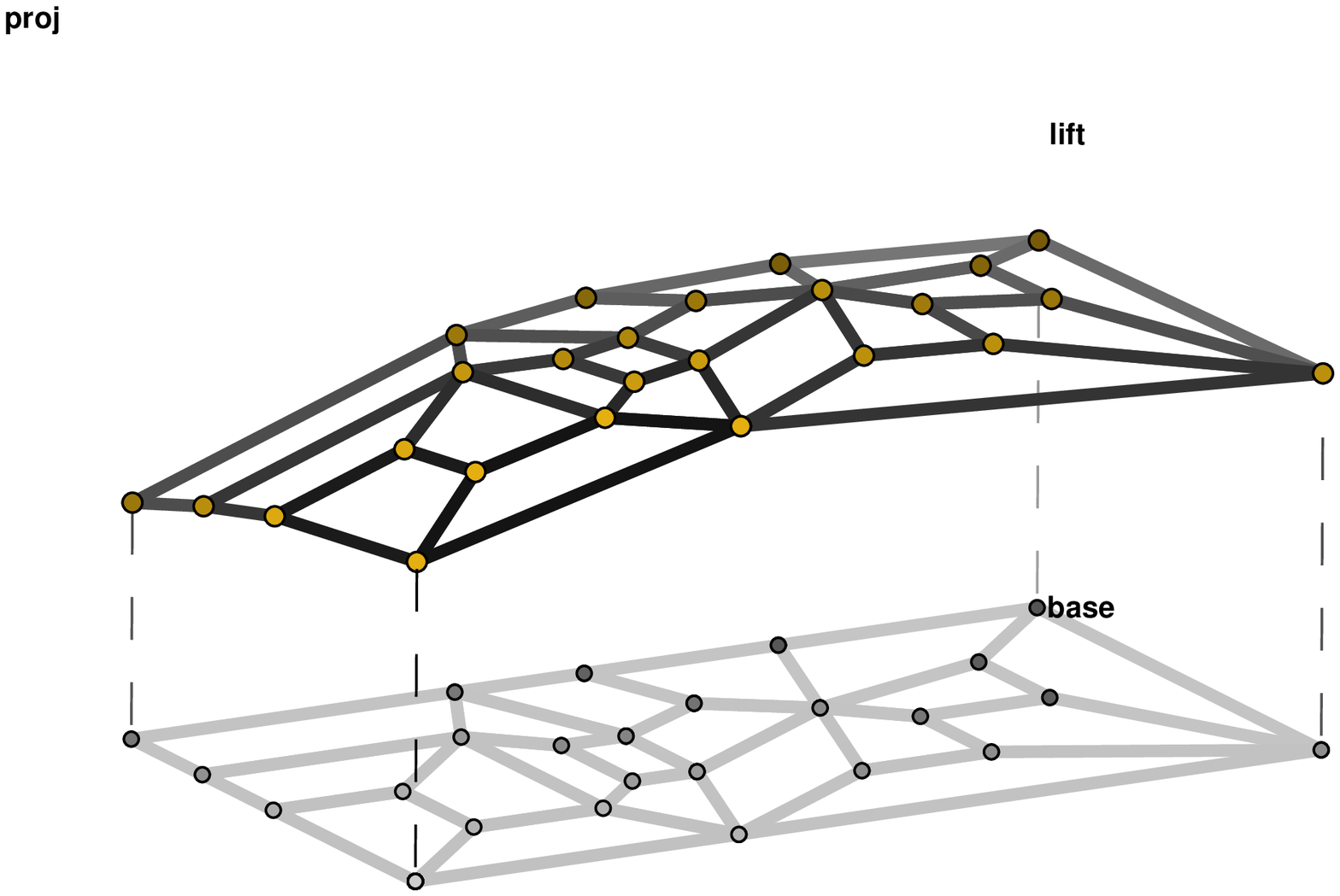}
         \quad\transformsToArrow{25mm}{10mm}
         \quad
    \psfrag{proj}{\rput(8mm,17mm){$\mathcal{B}'$}}
    \psfrag{base}{}
    \psfrag{lift}{}
    \psfrag{p}{\rput(2mm,3mm){$\boldsymbol{p}$}}
     \raisebox{1.1cm}{\includegraphics[width=.4\textwidth]{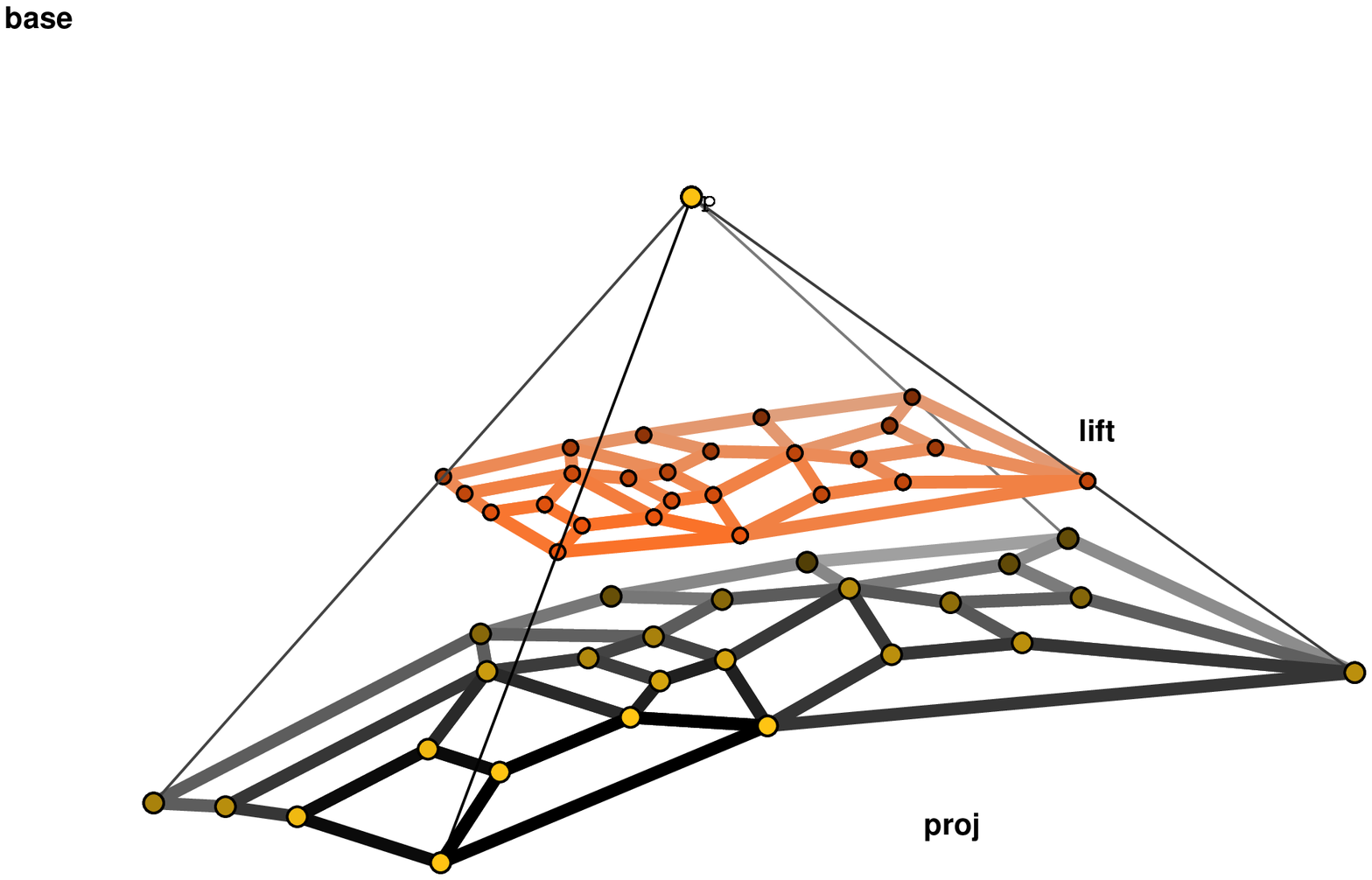}}
  \caption[Convexification of a regular cubical ball]{%
    Illustration of the `convexification' of a regular ball (Lemma \ref{lem:convexification}).}
  \label{fig:convexification}
\end{Figure}

\subsection{Lifted balls}\label{sec:lifted_balls}

When constructing cubical complexes we often deal with regular cubical
balls which are equipped with a lifting function. A \emph{lifted
  $d$-ball} is a pair $(\mathcal{B},h)$ consisting of a regular
$d$-ball~$\mathcal{B}$ and a lifting function $h$ of~$\mathcal{B}$.
The \emph{lifted boundary} of a lifted ball $(\mathcal{B},h)$ is the
pair $(\boundaryOP\mathcal{B},h|_{\boundaryOP\mathcal{B}})$.

If $(\mathcal{B},h)$ is a lifted $d$-ball in $\R^{d'}$ then
$\lifted{\mathcal{B}}{h}$ denotes the copy of~$\mathcal{B}$ in
$\R^{d'+1}$ with vertices
$(\boldsymbol{v},h(\boldsymbol{v}))\in\R^{d'+1}$,
$\boldsymbol{v}\in\vertices{\mathcal{B}}$. (In the sequel we sometimes
do not distinguish between these two interpretations of a lifted
ball.)  We rely on Figure~\ref{fig:lifted_ball} for the illustration
of this correspondence.

\begin{Figure}
    \psfrag{L}{\rput(2mm,2mm){$\mathcal{B}'=\lifted{\mathcal{B}}{h}$}}
    \psfrag{B}{\rput(0mm,0mm){$\mathcal{B}$}}

    \includegraphics[width=.4\textwidth]{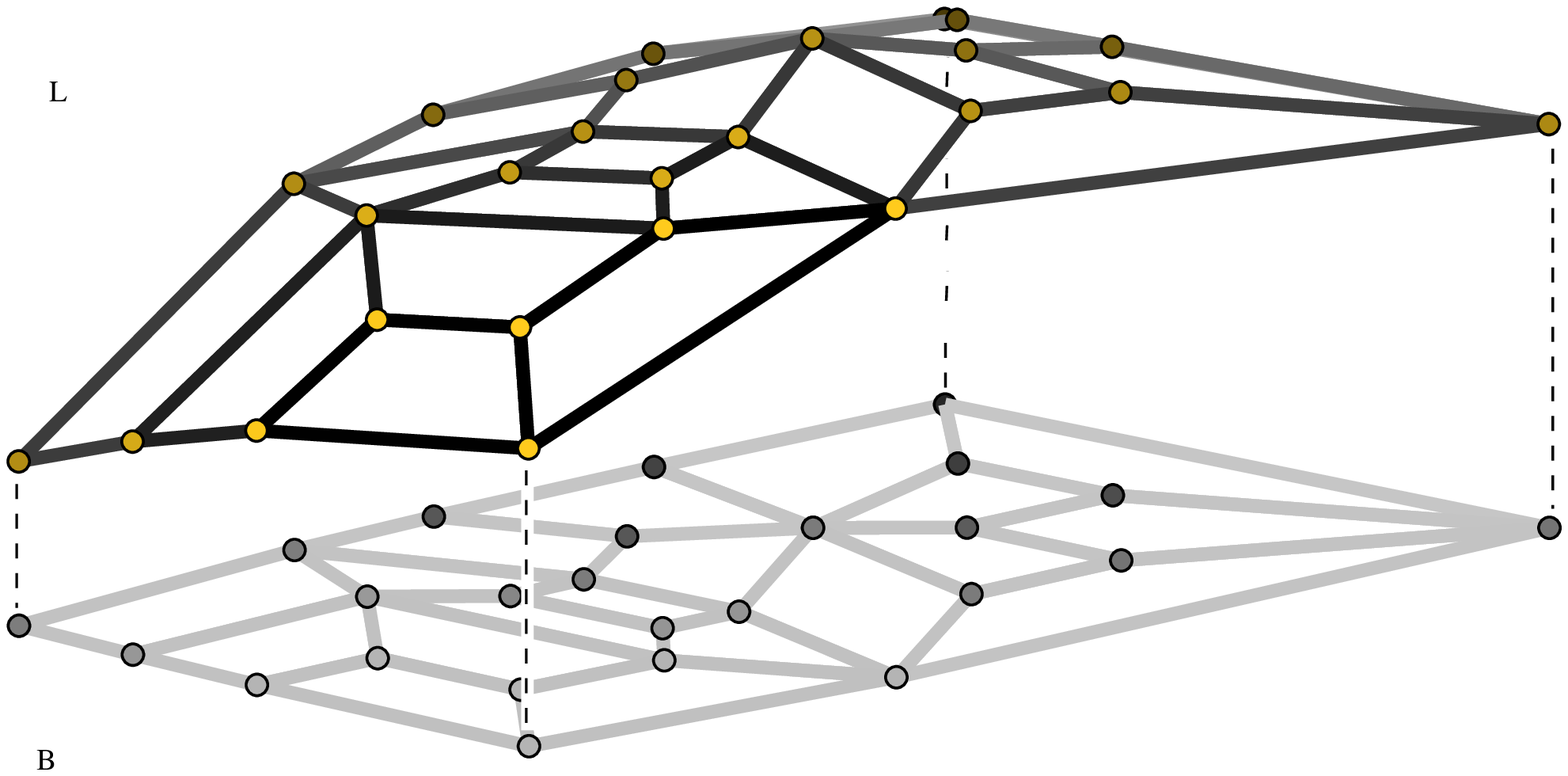}\qquad\qquad
    \psfrag{L}{\rput(2mm,-3mm){$Q$}}
    \includegraphics[width=.4\textwidth]{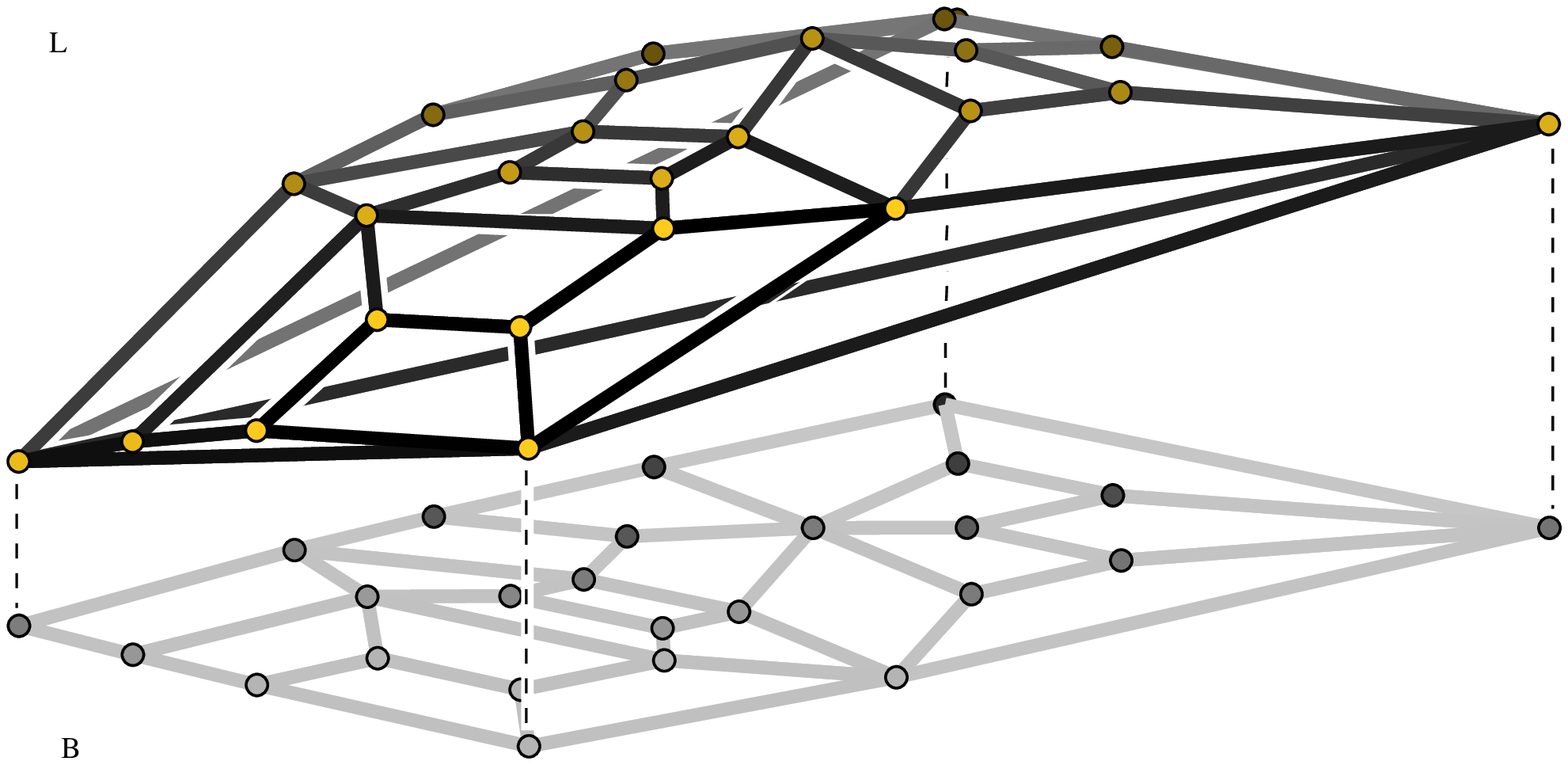}
        \caption{A lifted cubical ball $(\mathcal{B},h)$ and its
          lifted copy $\lifted{\mathcal{B}}{h}$. The figure on the
          right shows the convex hull~$Q=\convOP(\lifted{\mathcal{B}}{h})$. }
    \label{fig:lifted_ball}       
\end{Figure}

\begin{notation}  
  We identify $\R^d$ with $\R^d\times\{0\}\subset\R^{d+1}$, and
  decompose a point $\boldsymbol{x}\in\R^{d+1}$ as
  $\boldsymbol{x}=(\pi(\boldsymbol{x}),\gamma(\boldsymbol{x}))$, where
  $\gamma(\boldsymbol{x})$ is the last coordinate of~$\boldsymbol{x}$ and
  $\pi:\R^{d+1}\rightarrow\R^d$ is the projection that eliminates the
  last coordinate.
\end{notation}
 
Often a lifted ball $(\mathcal{B},\psi)$ is constructed as follows:
Let $P$ be a $d$-polytope (in $\R^d$) and $Q\subset\R^{d+1}$ a
$(d+1)$-polytope such that $\pi(Q)=P$.  Then the complex $\mathcal{B}'$ given
as the set of upper faces of~$Q$ determines a lifted polytopal
subdivision $(\mathcal{B},\psi)$ of $P$ (where
$\mathcal{B}:=\pi(\mathcal{B}')$ and $\psi$ is determined the vertex
heights $\gamma(\boldsymbol{v})$,
$\boldsymbol{v}\in\vertices{\mathcal{B}'}$). Hence
$\lifted{\mathcal{B}}{\psi}$ equals $\mathcal{B}'$.  Compare again
Figure~\ref{fig:lifted_ball}.

The \emph{lifted boundary subdivision} of a $d$-polytope $P$ is a pair
$(\mathcal{S}^{d-1},\psi)$ consisting of a polytopal subdivision
$\mathcal{S}^{d-1}$ of the boundary of~$P$ and a piece-wise linear
function $\psi:\support{\boundaryOP{P}}\rightarrow\R$ such that for
each facet $F$ of $P$ the restriction of $\psi$ to $F$ is a lifting
function of the induced subdivision $\mathcal{S}^{d-1}\cap F$ of $F$.

\subsection{The patching lemma}\label{sec:patching_lemma}

Often regular cubical balls are constructed from other regular balls.
The following ``patching lemma'', which appears
frequently in the construction of regular subdivisions (see
\cite[Cor.~1.12]{KKMS} or
\cite[Lemma~3.2.2]{BrunsGubeladzeTrung1997}) is a basic
tool for this.

\begin{notation}  For a $d$-polytope $P\subset\R^{d'}$, a
  polytopal subdivision $\mathcal{T}$ of $P$ and a hyperplane~$H$ in
  $\R^{d'}$, we denote by $\mathcal{T}\cap H$ the \emph{restriction} of
  $\mathcal T$ to~$H$, which is given by
  \[
  \mathcal{T}\cap H:=\setdef{F\cap H}{F\in\mathcal{T}}.
  \] 
  For two $d$-polytopes $P,Q$ with $Q\subset P$ and a
  polytopal subdivision $\mathcal{T}$ of $P$ we denote by
  $\mathcal{T}\cap Q$ the \emph{restriction} of $\mathcal T$ to~$Q$,
  which is given by
  \[
  \mathcal{T}\cap Q:=\setdef{F\cap Q}{F\in\mathcal{T}}.
  \] 
By $\facets{\mathcal{S}}$ we denote 
the set of facets of a complex $\mathcal{S}$.
\end{notation}

\begin{lemma}[``Patching lemma'']
  \label{lemma:patching_lemma}%
  Let~$Q$ be a $d$-polytope.  Assume we are given the following data:
   \begin{dense_itemize}
      \item A regular polytopal subdivision $\mathcal{S}$ of~$Q$
                  (the ``raw subdivision'').
       \item For each facet $F$ of~$\mathcal{S}$,
                a regular polytopal subdivision $\mathcal{T}_F$ of~$F$,\\
             such that $\mathcal{T}_F \cap F' = \mathcal{T}_{F'}\cap
                 F$ for all facets $F,F'$ of~$\mathcal{S}$.
                 
               \item For each facet $F$ of~$\mathcal{S}$,
                 a concave lifting function $h_F$ of $\mathcal{T}_F$,\\
                 such that $h_F(\boldsymbol x)=h_{F'}(\boldsymbol x)$
                 for all $\boldsymbol x \in F\cap F'$,
                 where $F,F'$ are facets of~$\mathcal{S}$.          
   \end{dense_itemize}
   Then this uniquely determines a regular polytopal
   subdivision~$\mathcal{U}=\bigcup_F \mathcal{T}_F$ of~$Q$ (the
   ``fine subdivision''). Furthermore, for every lifting
     function $g$ of~$\mathcal{S}$ there exists a small $\eps_0>0$
     such that for all $\eps$ in the range
     $\eps_0 > \eps > 0$ the function $g+\eps h $ is
     a lifting function of~$\mathcal{U}$, where $h$ is the piecewise
     linear function $h: \support{Q}\rightarrow\R$ which on each
     $F\in\mathcal{S}$ is given by $h_F$.
\end{lemma}

\begin{proof} Let $g$ be a lifting function of~$\mathcal{S}$.
  For a parameter $\eps>0$ we define a piece-wise linear function
  $\phi_\eps:\support{P}\rightarrow\R$ that on ${\boldsymbol x} \in F
  \in \facets{\mathcal{S}}$ takes the value $\phi_\eps(\boldsymbol{x})
  = g(\boldsymbol{x}) + \eps h_F(\boldsymbol{x})$. (It is well-defined
  since the $h_F$ coincide on the ridges of $\mathcal{S}$.) The
  domains of linearity of~$\phi_\eps$ are given by the facets of the
  ``fine'' subdivision~$\mathcal{U}$.  If $\eps$ tends to zero then
  $\phi_\eps$ tends to the concave function $g$.  This implies that
  there exists a small $\eps_0>0$ such that $\phi_\eps$ is concave and
  thus a lifting function of~$\mathcal{U}$, for $\eps_0 > \eps > 0$.
 \end{proof}

\subsection{Products and prisms}

\begin{lemma}[``Product lemma'']
  Let $(\mathcal{B}_1,h_1)$ be a lifted cubical $d_1$-ball in
  $\R^{d_1'}$ and $(\mathcal{B}_2,h_2)$ be a lifted cubical $d_2$-ball
  in $\R^{d_2'}$.\\
    Then the product $\mathcal{B}_1\times\mathcal{B}_2$ of
    $\mathcal{B}_1$ and $\mathcal{B}_2$ is a regular cubical
    $(d_1+d_2)$-ball in  $\R^{d_1'+d_2'}$.
\end{lemma}

\begin{proof}
  Each cell of $\mathcal{B}_1\times\mathcal{B}_2$ is a product of two
  cubes. Hence $\mathcal{B}_1\times\mathcal{B}_2$ is a cubical
  complex.  A lifting function $h$ of
  $\mathcal{B}_1\times\mathcal{B}_2$ is given by the sum of $h_1$ and
  $h_2$, that is, by $h((\boldsymbol{x},\boldsymbol{y})) :=
  h_1(\boldsymbol{x})+h_2(\boldsymbol{y})$, for
  $\boldsymbol{x}\in\support{\mathcal{B}_1},
  \boldsymbol{y}\in\support{\mathcal{B}_2}.$
\end{proof}

As a consequence, the \emph{prism} $\operatorname{prism}(\mathcal{C})$
over a cubical $d$-complex $\mathcal{C}$ yields a cubical
$(d+1)$-dimensional complex. Furthermore, the prism over a regular
cubical ball $\mathcal{B}$ yields a regular cubical $(d+1)$-ball.

\subsection{Piles of cubes}\label{subsec:piles_of_cubes}

For integers $\ell_1,\ldots,\ell_d\geq 1$, the \emph{pile of cubes}
$\PileOfCubes{d}{\ell_1,\ldots,\ell_d}$ is the cubical $d$-ball formed
by all unit cubes with integer vertices in the $d$-polytope
$P := [0,\ell_1] \times \ldots \times [0,\ell_d]$,
that is, the cubical $d$-ball formed by the set of all $d$-cubes
\[
     C(k_1,\ldots,k_d):= [k_1, k_1+1] \times \dots \times [k_d, k_d+1]
\]
for integers $0 \leq k_i < \ell_i$ together with their faces
\cite[Sect.~5.1]{Z35}.  

The pile of cubes $\PileOfCubes{d}{\ell_1,\ldots,\ell_d}$ is a product
of $1$-dimensional subdivisions, which are regular. Hence the product
lemma implies that $\PileOfCubes{d}{\ell_1,\ldots,\ell_d}$ is a
regular cubical subdivision of the $d$-polytope $P$.

\subsection{Connector polytope}\label{subsec:connector_cube}

The following construction yields a ``connector'' polytope that may be
used to attach cubical $4$-polytopes resp.\ regular cubical $4$-balls
without the requirement that the attaching facets are projectively
equivalent.

\begin{lemma}\label{lem:connector_cube}
  For any combinatorial $3$-cube~$F$ there is a combinatorial
  $4$-cube $C$ that has both (a projective copy of) $F$ and a 
  regular $3$-cube $F'$ as (adjacent) facets.
\end{lemma}

\begin{proof}
  After a suitable projective transformation we may assume that
  $F\subset\R^3$ has a unit square~$Q$ as a face.  Now the
  prism $F\times I$ over $F$ has $F$ and $Q\times I$ as
  adjacent facets, where the latter is a unit cube.
\end{proof}


\section{Basic construction techniques}\label{sec:constructions}

\subsection{Lifted prisms}\label{subsec:cubical_prisms}\label{subsec:lifted_prisms}

While there appears to be no simple construction that would produce a
cubical \mbox{($d+1$)-}\allowbreak polytope from a given cubical
$d$-polytope, we do have a simple prism construction 
that produces regular cubical ($d+1$)-balls from
regular cubical $d$-balls.

\begin{construction}{Lifted prism}
\begin{tabular}{@{}ll}
  \textbf{Input:}  & A lifted cubical $d$-ball $(\mathcal{B},h)$.\\
  \textbf{Output:} & A lifted cubical ($d+1$)-ball
  $\LiftedPrism{\mathcal{B}}{h}$\\
& 
   which is combinatorially isomorphic to the prism over $\mathcal{B}$.
\end{tabular}

We may assume that the convex lifting function $h$ defined
on $P:=\support{\mathcal{B}}$ is strictly posi\-tive.
 Then the lifted facets
of $\LiftedPrism{\mathcal{B}}{h}$ may be taken to be the sets
\[
  \widetilde{F}\ \ :=\ \ \{(\boldsymbol x,t,h(\boldsymbol x)): \boldsymbol x\in F,\ 
                               -h(\boldsymbol x) \le t \le +h(\boldsymbol x)\}, \qquad F \in \facets{\mathcal{B}}.
\vspace{-7mm}
\]
\end{construction}
  
If $\mathcal{B}$ does not subdivide the boundary of $P$, then
$\LiftedPrism{\mathcal{B}}{h}$ does not subdivide the boundary of
$\support{\LiftedPrism{\mathcal{B}}{h}}$.  In this case $\widehat
P:=\support{\LiftedPrism{\mathcal{B}}{h}}$ is a cubical
$(d+1)$-polytope whose boundary 
complex is combinatorially isomorphic to the boundary of the prism over~$\mathcal{B}$.
The $f$-vector of~$\widehat P$ is then given by
  \[ 
f_k(\widehat P)\ \ =\ \ 
\begin{cases} 
2 f_0(\mathcal{B})               & \textrm{ for }k=0,\\
2 f_k(\mathcal{B}) + f_{k-1}(\boundaryOP\mathcal{B})  & \textrm{ for }0<k\le d.
\end{cases}
\]

Figure~\ref{fig:lifted_prism} shows the lifted prism over a lifted
cubical $2$-ball.
\begin{Figure}
   \begin{psfrags}
    \psfrag{proj}{}
    \psfrag{base}{\rput(-26mm,-15mm){$\mathcal{B}$}}
    \psfrag{lifted}{\rput(-25mm,-8mm){$\lifted{\mathcal{B}}{h}$}}

   \includegraphics[width=.4\textwidth]{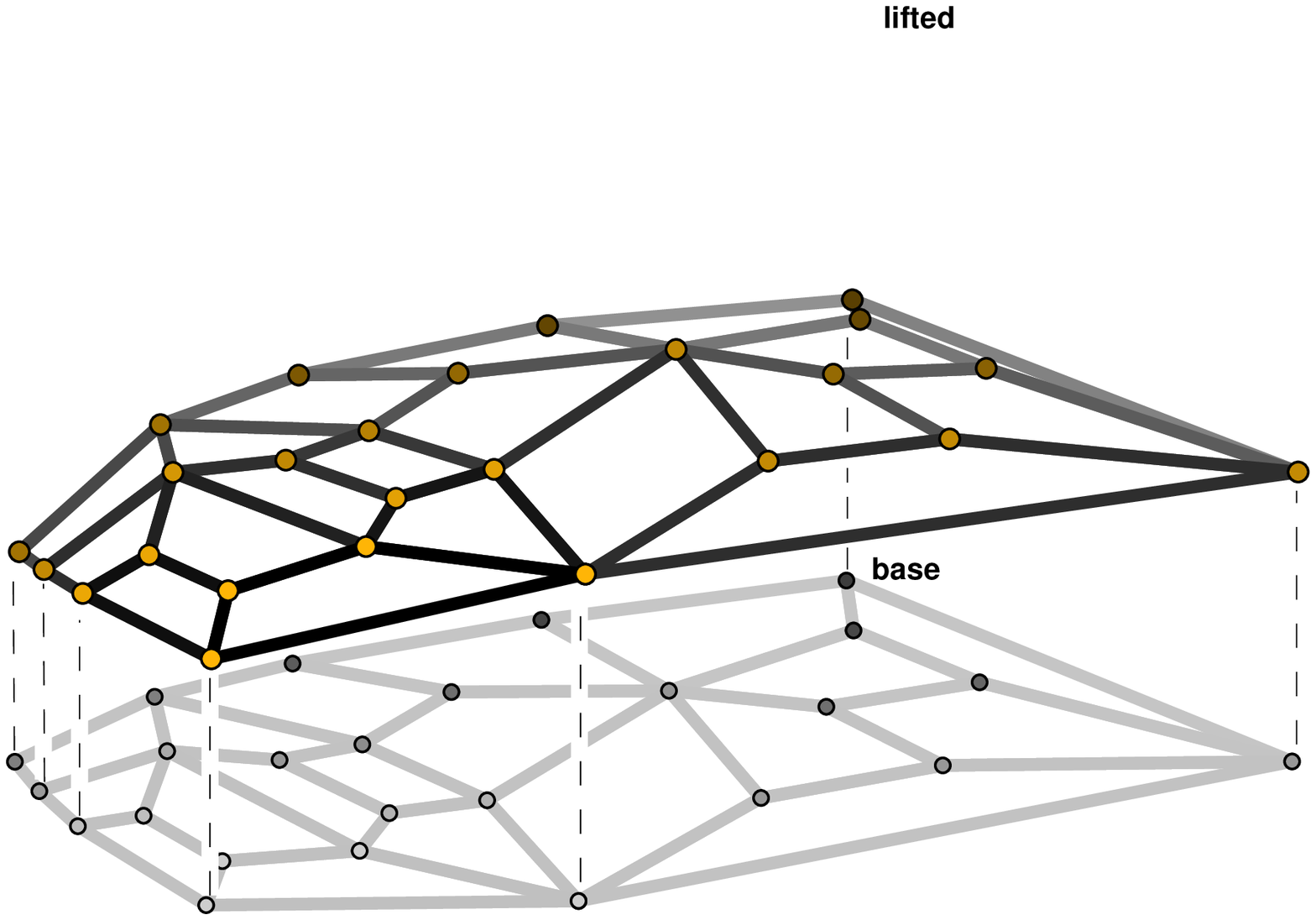}\quad
          \transformsToArrow{19mm}{10mm}\qquad
   \includegraphics[width=.4\textwidth]{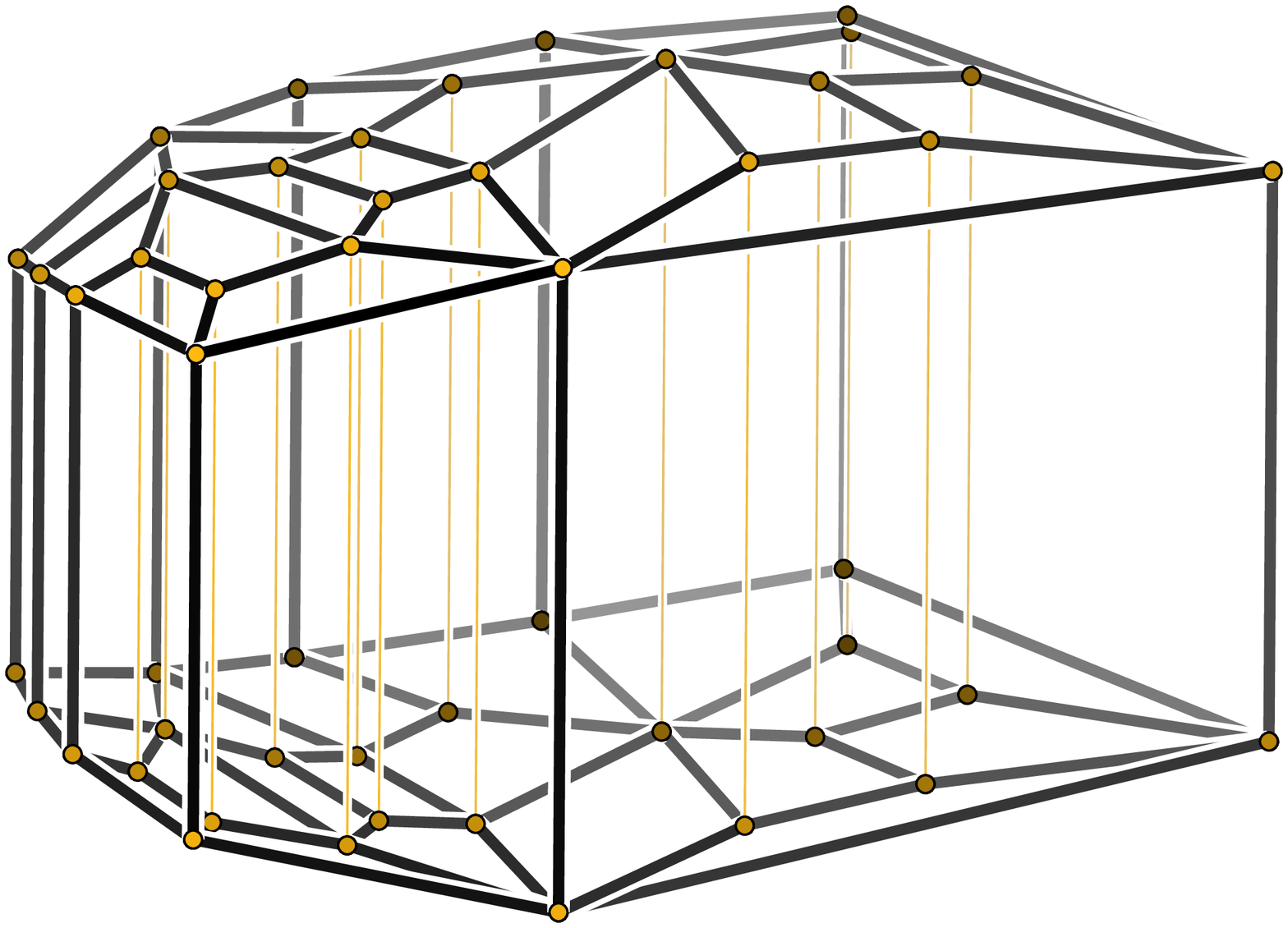}
   \rput(-24mm,-3mm){$\LiftedPrism{\mathcal{B}}{h}$}
        \end{psfrags}
   \vspace{4mm}
   \caption[The lifted prism over a lifted cubical ball]{%
     The lifted prism of a lifted cubical $d$-ball~$(\mathcal{B},h)$,
     displayed for $d=2$.  The result is a (regular) cubical $(d+1)$-ball
     that is combinatorially isomorphic to the prism over $\mathcal{B}$. }
   \label{fig:lifted_prism}
\end{Figure}

\begin{proposition}[Dual manifolds]\label{remark:dmfs_of_lifted_prism}
  Up to PL-homeomorphism, the cubical ball
  $\LiftedPrism{\mathcal{B}}{h}$ has the following dual manifolds:
  \begin{compactitem}
  \item $\mathcal{N}\times I$ for each dual manifold $\mathcal N$ 
    of~$\mathcal{B}$,
  \item one $(d-1)$-sphere combinatorially isomorphic 
    to~$\boundaryOP\mathcal{B}$. \qed
  \end{compactitem}
\end{proposition}

\subsection{Lifted prisms over two balls}\label{subsec:cubical_prisms_two_balls}

Another modification of this construction is to take two different lifted
cubical balls $(\mathcal{B}_1,h_1)$ and $(\mathcal{B}_2,h_2)$ with the same
lifted boundary complex (that is,
$\boundaryOP\mathcal{B}_1=\boundaryOP\mathcal{B}_2$ with
$h_1(\boldsymbol{x})=h_2(\boldsymbol{x})$ for all
$\boldsymbol{x}\in\boundaryOP\mathcal{B}_1=\boundaryOP\mathcal{B}_2$)
as input.  In this case the outcome is a cubical ($d+1$)-polytope which
may not even have a cubification.

\begin{construction}{Lifted prism over two balls}\label{constr:Lifted_prism_over_two_balls}
\begin{inputoutput}
  \item[Input:] Two lifted cubical $d$-balls 
  $(\mathcal{B}_1,h_1)$ and $(\mathcal{B}_2,h_2)$\\ with the same lifted
  boundary.
\item[Output:] A cubical ($d+1$)-polytope~$\LiftedPrismOverTwoBalls{\mathcal{B}_1}{h_1}{\mathcal{B}_2}{h_2}$\\
  with lifted copies of $\mathcal{B}_1$ and $\mathcal{B}_1$ in its boundary.
\end{inputoutput}

If both balls do not subdivide their boundaries, we set
$\mathcal{B}'_k:=\mathcal{B}_k$ and $h'_k:=h_k$ for $k\in\{1,2\}$.
Otherwise we apply the construction of the proof of
Lemma~\ref{lem:convexification} simultaneously to both lifted cubical
balls $(\mathcal{B}_1,h_1)$ and $(\mathcal{B}_2,h_2)$ to obtain two
lifted cubical $d$-balls $(\mathcal{B}'_1,h'_1)$ and
$(\mathcal{B}'_2,h'_2)$ with the same support
$Q=\support{\mathcal{B}_1}=\support{\mathcal{B}_2}$ which do not
subdivide the boundary of $Q$.

We can assume that $h'_1,h'_2$ are strictly positive. 
Then~$\widehat Q:=\LiftedPrismOverTwoBalls{\mathcal{B}_1}{h_1}{\mathcal{B}_2}{h_2}$
is defined as the convex hull of the points in
\[
    \{(\boldsymbol x,+h'_1(\boldsymbol x)): \boldsymbol x\in \support{\mathcal{B}'_1}\}\ \cup\   
    \{(\boldsymbol x,-h'_2(\boldsymbol x)): \boldsymbol x\in \support{\mathcal{B}'_2}\},
\]%

\label{p:lifted_prism_over_two_balls}%
Since $\mathcal{B}'_1$ and $\mathcal{B}'_2$ both do not subdivide their
boundaries, each of their proper faces yields a face of~$\widehat Q$.
Furthermore, $\widehat Q$ is a cubical $(d+1)$-polytope whose
$f$-vector is given by
\[
f_k(\widehat Q)\ \ =\ \ 
\begin{cases} 
f_0(\mathcal{B}_1) + f_0(\mathcal{B}_2)               & \textrm{ for }k=0,\\
f_k(\mathcal{B}_1) + f_k(\mathcal{B}_2) 
   + f_{k-1}(\boundaryOP\mathcal{B}_1)  & \textrm{ for }0<k\le d.
\end{cases}
\]
See Figure~\ref{fig:lifted_prism_over_two_balls}.
\end{construction}

\begin{Figure}
   \begin{psfrags}
    \psfrag{proj}{}
    \psfrag{base}{\rput(-26mm,-14mm){$\mathcal{B}_1$}}
    \psfrag{lifted}{\rput(-25mm,0mm){$\lifted{\mathcal{B}_1}{h_1}$}}
    \includegraphics[width=.3\textwidth]{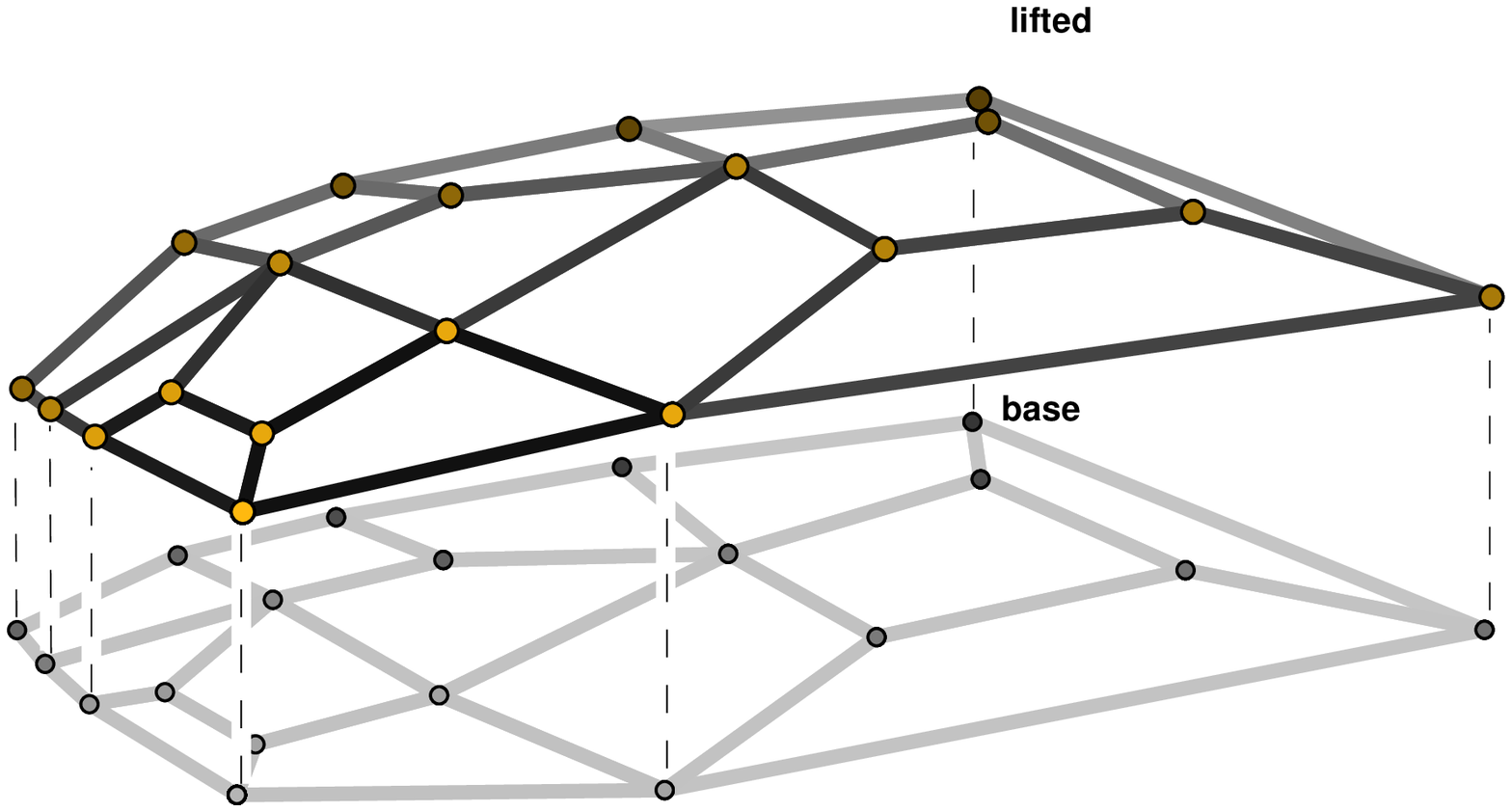}\quad
     \raisebox{15mm}{$+$}\quad
    \psfrag{base}{\rput(-26mm,-14mm){$\mathcal{B}_2$}}
    \psfrag{lifted}{\rput(-25mm,-7mm){$\lifted{\mathcal{B}_1}{h_2}$}}
   \includegraphics[width=.3\textwidth]{img/lifted_prism/convexified_lifted_with_base_corrected}\\
          \transformsToArrow{15mm}{10mm}\qquad
   \includegraphics[width=.3\textwidth]{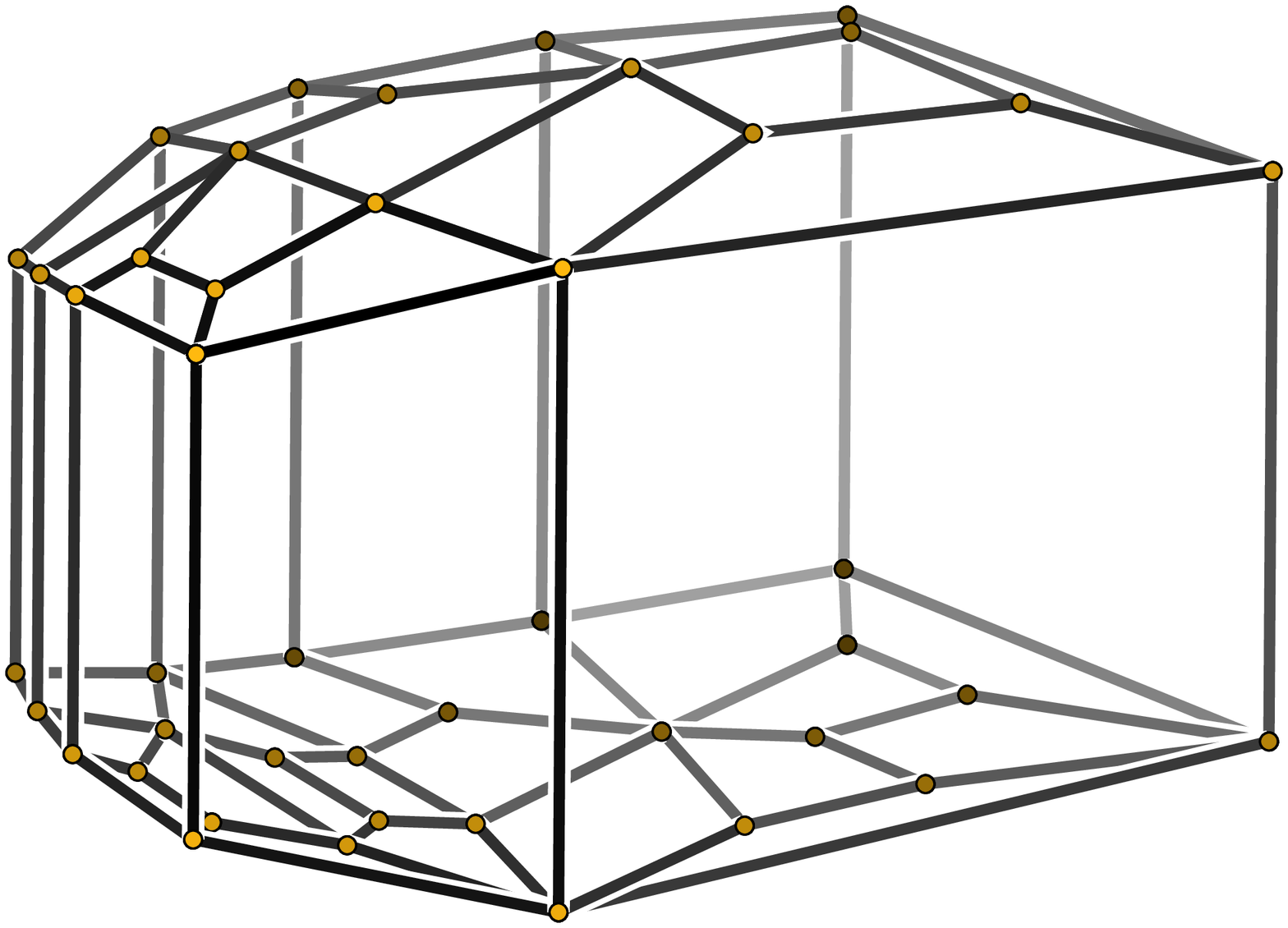}
   \rput(8mm,0mm){$\LiftedPrismOverTwoBalls{\mathcal{B}_1}{h_1}{\mathcal{B}_2}{h_2}$}
        \end{psfrags}
   \caption[The lifted prism over two lifted balls]{%
     The lifted prism over two lifted cubical $d$-balls 
     $(\mathcal{B}_1,h_1)$ and $(\mathcal{B}_2,h_2)$, displayed for $d=2$.
     The outcome is a cubical $(d+1)$-polytope.}
   \label{fig:lifted_prism_over_two_balls}
\end{Figure}

\subsection{Schlegel caps}\label{subsec:scap}

The following is a projective variant of the prism construction,
applied to a $d$-polytope $P$.

\begin{construction}{Schlegel cap}
\begin{inputoutput}
\item[Input:] An almost cubical $d$-polytope~$(P,F_0)$ 
\item[Output:] A regular cubical $d$-ball $\SchlegelCap(P,F_0)$, with
$P\subset\support{\SchlegelCap(P,F_0)}$
which is combinatorially isomorphic to the prism over $\Schlegel (P,F)$.
\end{inputoutput}

The construction of the Schlegel cap depends on two further pieces of
input data, namely on a point $\boldsymbol{x}_0\in{\mathbb R}^d$ 
beyond~$F_0$ (and beneath all other facets of $P$; 
cf.~\cite[Sect.~5.2]{Gr1-2}) and on 
a hyperplane~$H$ that separates $\boldsymbol{x}_0$ from~$P$.
In terms of projective transformations it is obtained as follows:
\begin{compactsteps}
\item Apply a projective transformation that moves $\boldsymbol{x}_0$
  to infinity while fixing~$H$ pointwise. This transformation moves
  the Schlegel complex $\ComplexOf{\boundaryOP P}{\setminus}\{F_0\}$
  to a new cubical complex $\mathcal{E}$.
\item Reflect the image $\mathcal{E}$ of the Schlegel complex 
  in $H$, and call its reflected copy $\mathcal{E}'$.
\item
      Build the polytope bounded by $\mathcal{E}$ and $\mathcal{E}'$.
\item Reverse the projective transformation of (1).
\end{compactsteps}
\begin{Figure}
    \psfrag{5}{\rput(-1.5mm,-3mm){$\boldsymbol{x}_0$}}
    \psfrag{0}{\rput(3mm,8mm){$F_0$}}
    \psfrag{1}{\rput(3.5mm,0){$P$}}
    \psfrag{2}{}
    \psfrag{3}{}
    \psfrag{4}{}
    \psfrag{H}{\rput(0mm,-3.2cm){$H$}}
   \includegraphics[height=3.2cm]{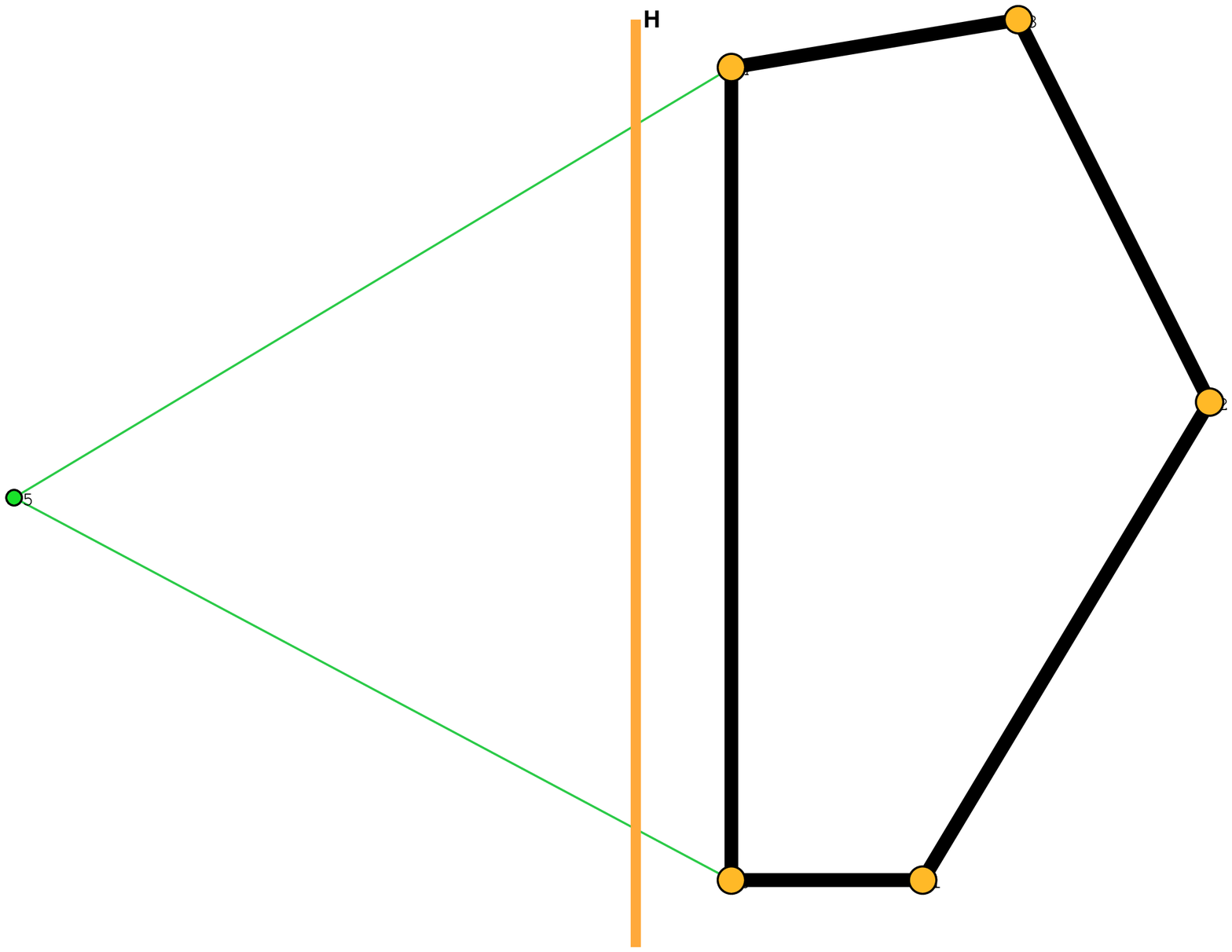}\ \transformsToArrow{13mm}{5mm}\ 
    \psfrag{H}{}
    \includegraphics[height=3.2cm]{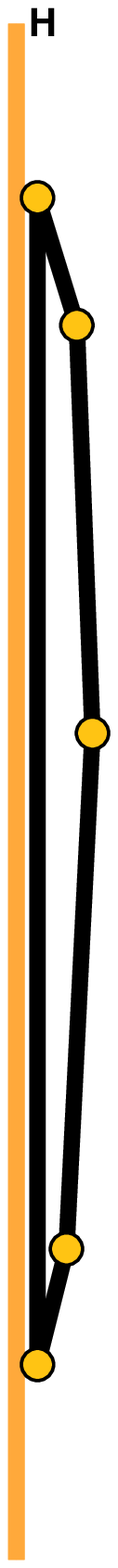}\ \transformsToArrow{13mm}{5mm}\quad
    \includegraphics[height=3.2cm]{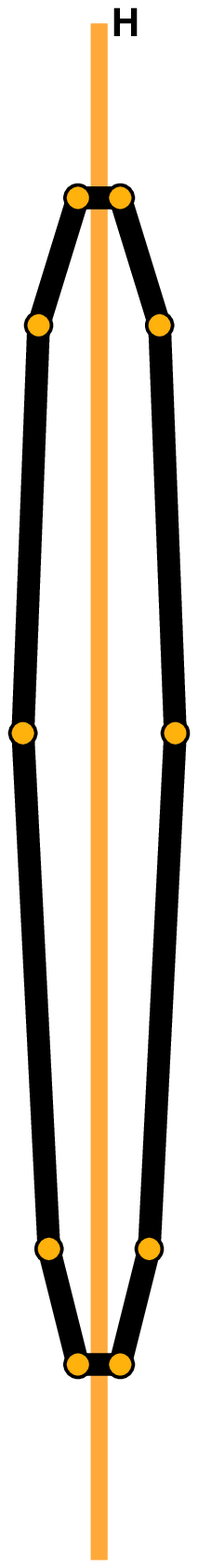}\transformsToArrow{13mm}{5mm}\qquad\ 
    \psfrag{H}{\rput(0mm,-3.4cm){$H$}}
    \psfrag{0}{}
    \psfrag{1}{}
    \includegraphics[height=3.5cm]{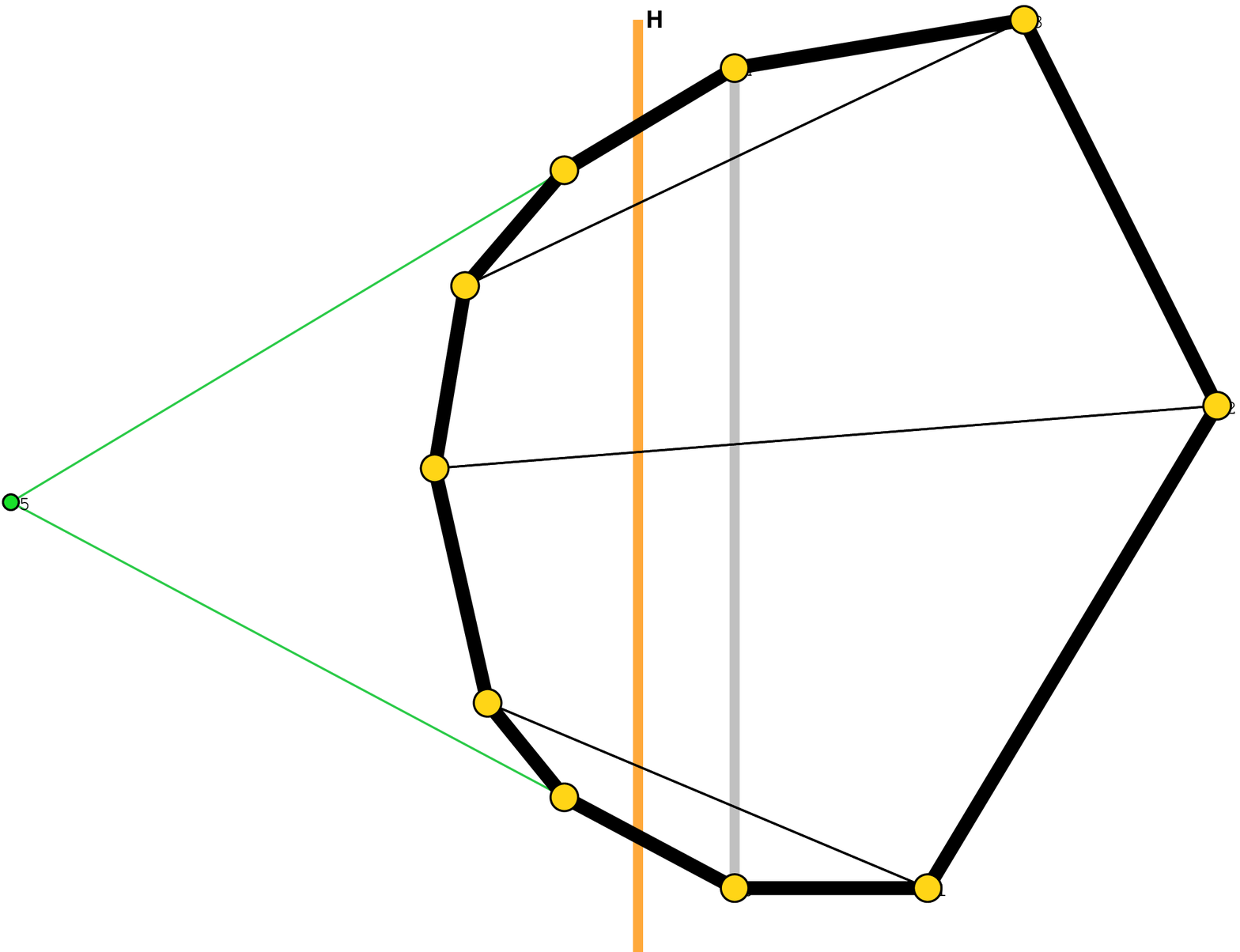}
   \vspace*{6mm}  
      \caption[The Schlegel cap construction]{
        Construction steps of the Schlegel cap over an almost cubical polytope.}
\end{Figure}
\vspace*{-2ex}
An alternative description, avoiding projective transformations, is
as follows:
\begin{compactsteps}
\item For each point $\boldsymbol{x}$ in the Schlegel complex 
  $\ComplexOf{\boundaryOP P}{\setminus}\{F_0\}$ let
  $\bar{\boldsymbol{x}}$ be the intersection point of~$H$ and the
  segment  $[\boldsymbol{x}_0,\boldsymbol{x}]$,
  and let $\boldsymbol x'$ be the point on the segment 
  $[\boldsymbol{x}_0,\boldsymbol{x}]$ such that
  $[\boldsymbol{x}_0,\bar{\boldsymbol{x}}; \boldsymbol{x}',\boldsymbol{x}]$ 
  form a harmonic quadruple (cross ratio~$-1$). 

  That is, if $\vec{\boldsymbol{v}}$ is a direction vector such
  that $\boldsymbol{x} = \boldsymbol{x}_0 + t \vec{\boldsymbol{v}}$ for 
  some $t>1$ denotes the difference $\boldsymbol{x}-\boldsymbol{x}_0$,
  while $\bar{\boldsymbol{x}} = \boldsymbol{x}_0 + \vec{\boldsymbol{v}}$ 
  lies on $H$, then 
  $\boldsymbol{x}' = \boldsymbol{x}_0 + \frac t{2t-1}\vec{\boldsymbol{v}}$.
\item 
  For each face $G$ of the Schlegel complex, 
  $G':=\{\boldsymbol{x}':\boldsymbol{x}\in G\}$ is the 
  ``projectively reflected'' copy of
  $G$ on the other side of $H$.
\item 
  The Schlegel cap $\SchlegelCap(P,F_0)$
  is the regular polytopal ball with faces $G$, $G'$ and $\conv(G\cup G')$ for
  faces $G$ in the Schlegel complex.
\end{compactsteps}
\vskip-5mm
\begin{Figure}
   \psfrag{x_0}{\rput(0mm,0){$\boldsymbol{x}'$}}
   \psfrag{x'}{\rput(0mm,0mm){$\boldsymbol{x}_0$}}
   \psfrag{bar_x}{\rput(.5mm,1.5mm){$\bar{\boldsymbol{x}}$}}
   \psfrag{x}{$\boldsymbol{x}$}
   \psfrag{F}{}
   \psfrag{P}{$P$}
   \psfrag{H}{$H$}
   \includegraphics[width=.37\textwidth]{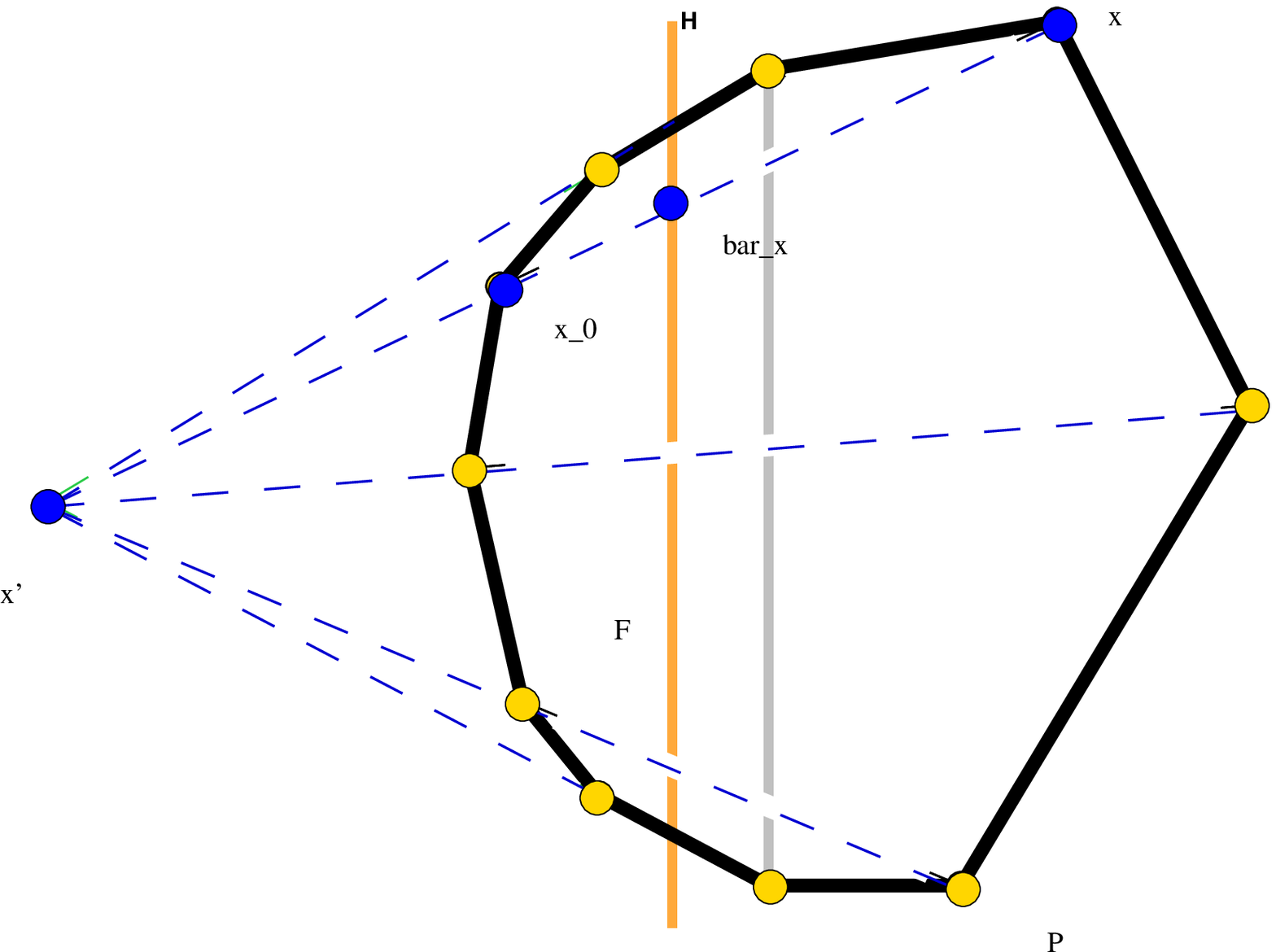}
        \caption{Constructing the Schlegel cap via cross ratios.}
\end{Figure}\vskip-4mm
\end{construction}

\section{A small cubical 4-polytope with a dual Klein bottle}\label{sec:klein}

In this section we present the first instance of a cubical
$4$-polytope with a non-orientable dual manifold. By
Proposition~\ref{prop:NEO} this instance is not edge-orientable.
Hence, its existence also confirms the conjecture of Hetyei
\cite[Conj.~2, p.~325]{hetyei95:_stanl}.  Apparently this is the first
example of a cubical polytope with a non-orientable dual manifold.

\begin{theorem}\label{thm:C4P_with_dual_klein}
  There is a cubical $4$-polytope $P_{72}$ with $f$-vector
  \[
     f(P_{72})\ \ =\ \ (72,\,196,\,186,\,62),
  \]
  one of whose dual manifolds is an immersed Klein bottle of
  $f$-vector $(80,\,160,\,80)$.
\end{theorem}


\paragraph{Step 1.}
We start with a cubical octahedron $O_8$, the smallest cubical $3$-polytope
    that is not a cube, with $f$-vector
\[
            f(O_8)\ \ =\ \ (10,\,16,\,8).
\]
\begin{figure}[H]\centering
    \includegraphics[width=60mm,height=55mm]{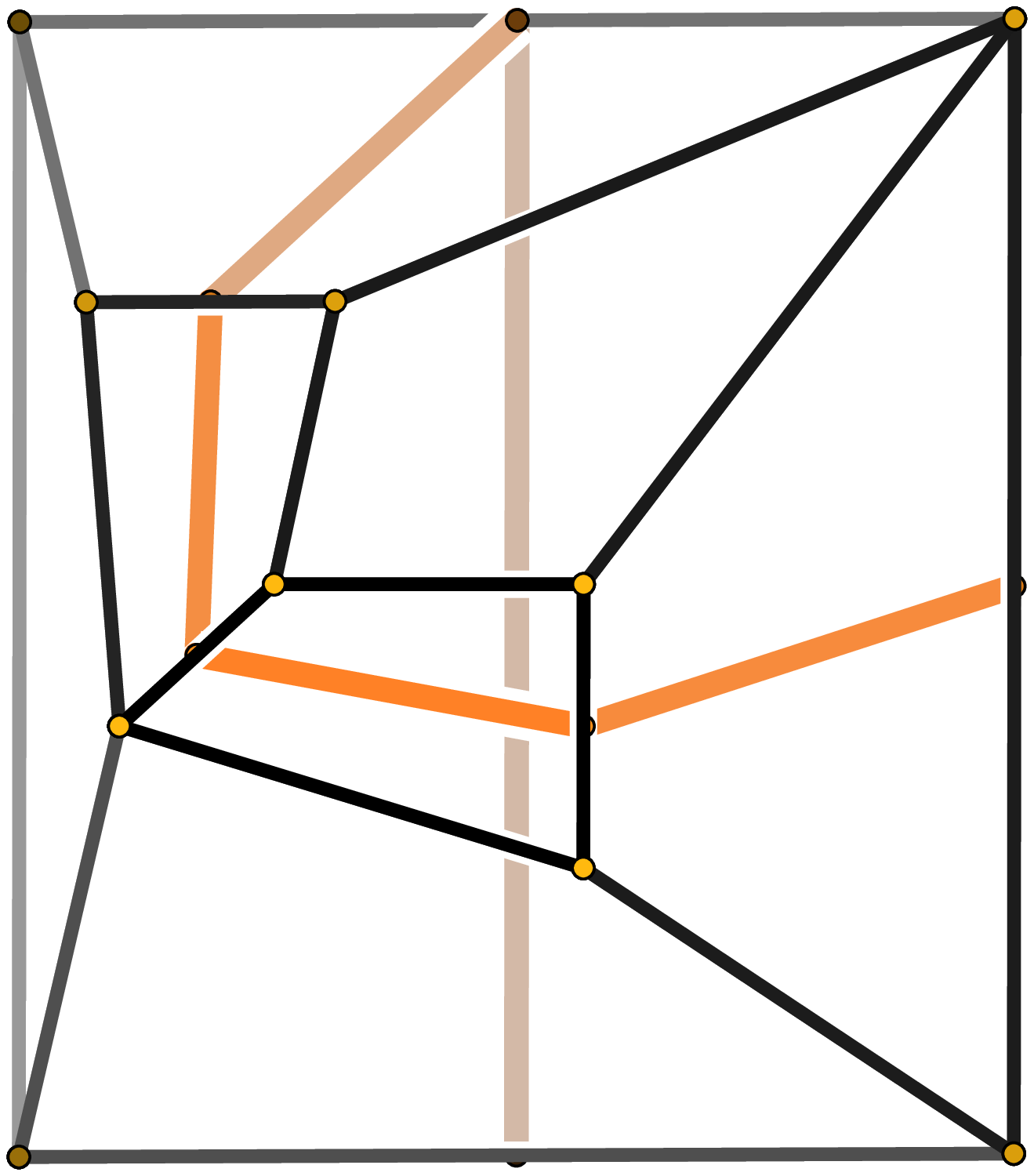}
   \caption[The cubical octahedron]{The cubical octahedron $O_8$ positioned in $\R^3$ 
     with a regular square base facet $Q$ and acute dihedral angles at
     this square base. A part of one dual manifold
     is highlighted.}
   \label{fig:O_8_with_dmf}
\end{figure}
We may assume that $O_8$ is already positioned in $\R^3$ 
with a regular square base facet $Q$ and acute dihedral angles
at this square base; compare the figure below.
The $f$-vector of any Schlegel diagram of $O_8$ is 
\[
           f(\Schlegel(O_8,Q))\ \ =\ \ (10,\,16,\,7).
\]
Let $O_8'$ be a congruent copy of $O_8$, obtained by reflection of
$O_8$ in its square base followed by a $90^\circ$ rotation around the
axis orthogonal to the base; compare the 
figure  below. This results in a regular
$3$-ball with cubical $2$-skeleton. Its $f$-vector is
\[
        f(\mathcal{B}_{2})\ \ =\ \ (16,\,28,\,15,\,2).
\]
The special feature of this complex is that it contains a cubical
M\"obius strip with parallel inner edges of length~$9$ in its
$2$-skeleton, as is illustrated in the figure.

\begin{figure}[H]\centering
  \psfrag{F}{$F$}
  \includegraphics[width=60mm,height=55mm]{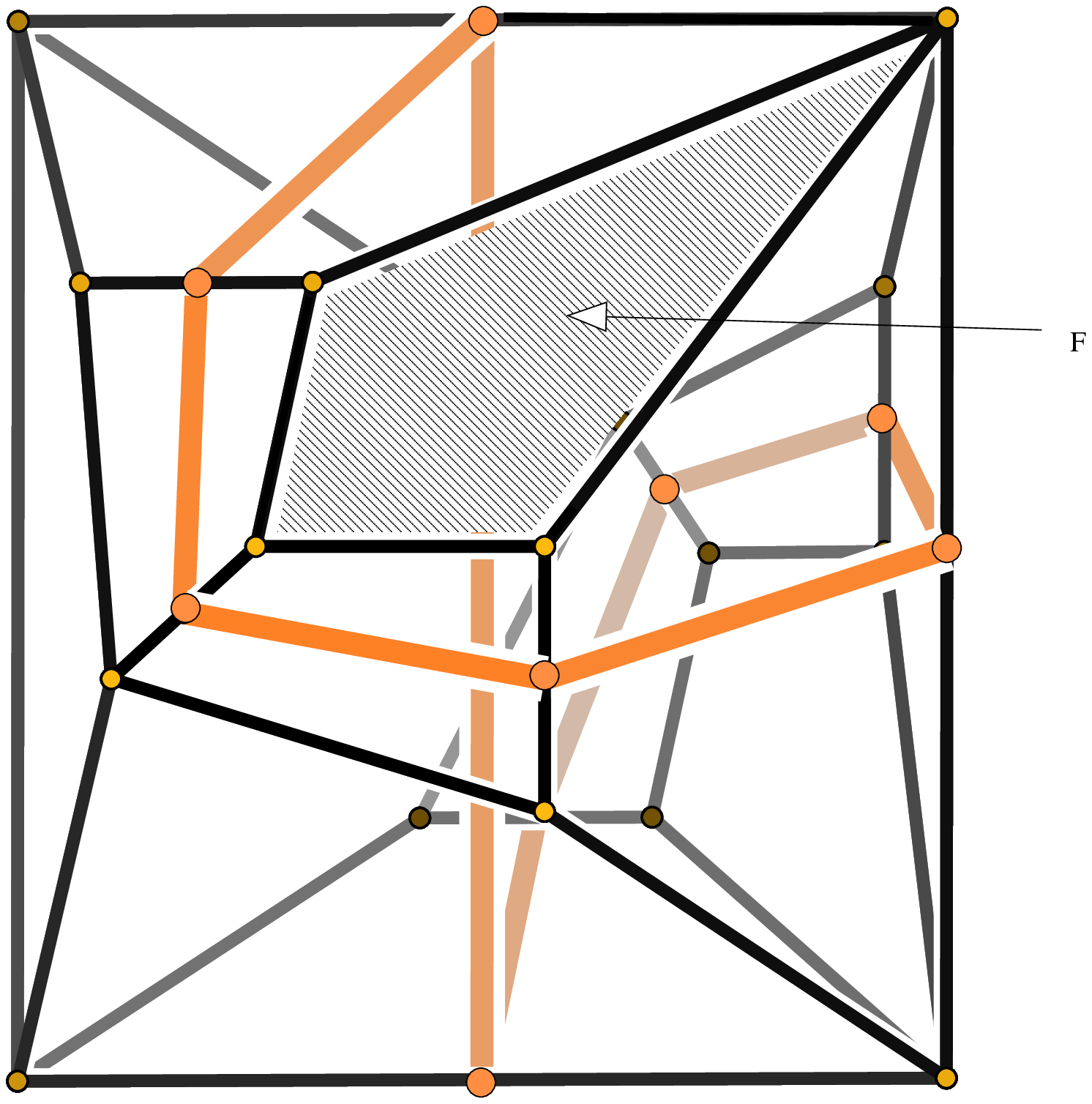}
   \caption{The outcome of
     step 1 of the construction: The $2$-cubical convex $3$-ball
     $\mathcal{B}_2$ which contains a M\"obius strip with parallel
     inner edges in the $2$-skeleton.}
\end{figure}

\paragraph{Step 2.}
Now we perform a Schlegel cap construction on $O_8$, based on the
(unique) facet~$F$ of~$O_8$ that is not contained in the M\"obius
strip mentioned above, and that is not adjacent to the square glueing
facet~$Q$.  This Schlegel cap has the $f$-vector
\[
     f(\mathcal{S}_7)\ \ =\ \ (20,\,42,\,30,\,7),
\]
while its boundary has the $f$-vector
\[
     f(\boundaryOP \mathcal{S}_7)\ \ =\ \ (20,\,36,\,18).
\]

\begin{figure}[ht]\centering
      \includegraphics[width=60mm,height=55mm]{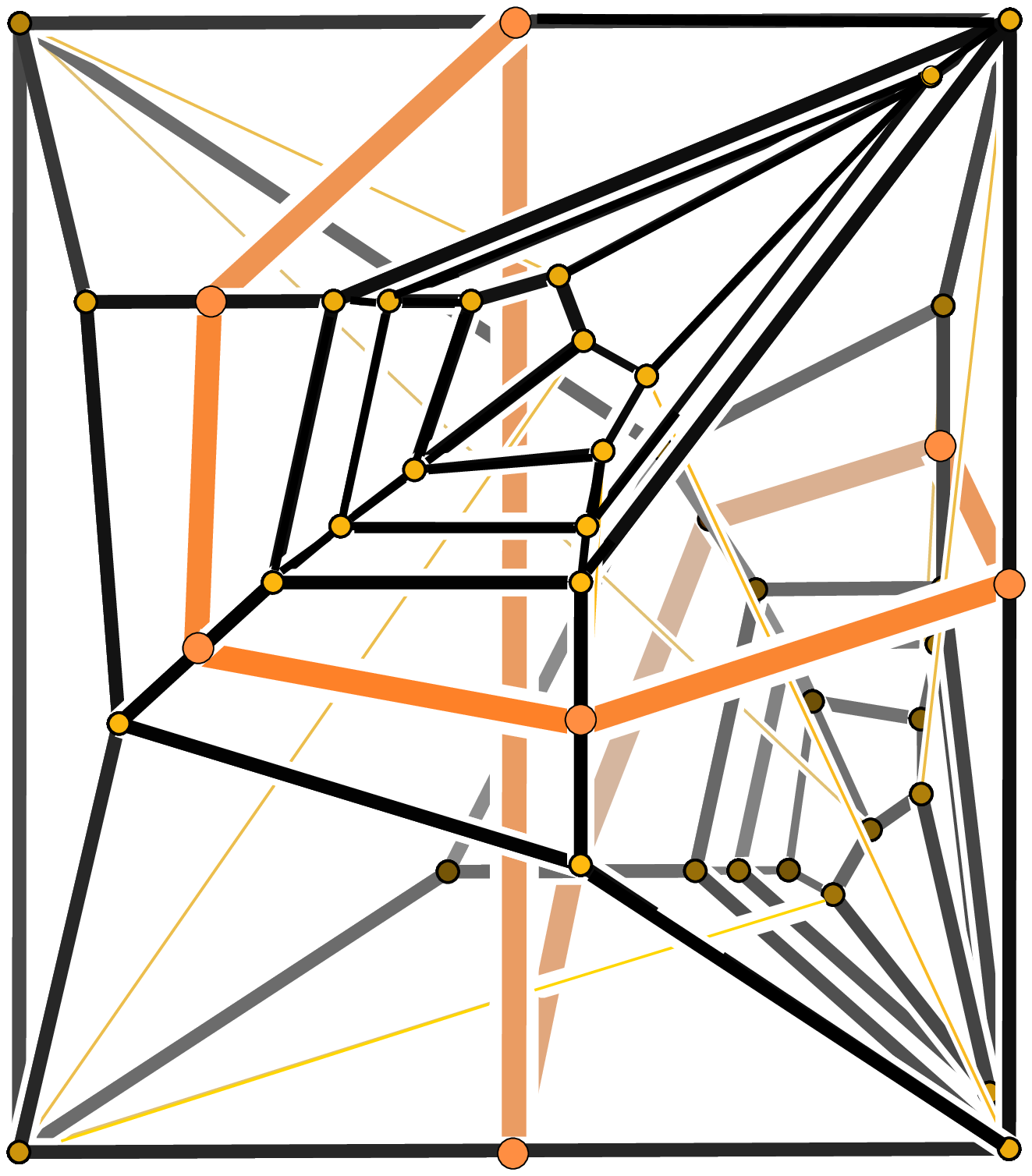}
   \caption[A not edge-orientable cubical $2$-ball]{
     The outcome of step 2 of the construction: The cubical convex
     $3$-ball~$\mathcal{B}_{14}$ which contains a M\"obius strip with
     parallel inner edges in the $2$-skeleton.}
  \label{fig:neo_3_cubical_ball}
\end{figure}

\paragraph{Step 3.}
The same Schlegel cap operation may be performed on the second copy
$O_8'$.  Joining the two copies of the Schlegel cap results in a
regular cubical $3$-ball $\mathcal{B}_{14}$ with $f$-vector
\[
    f(\mathcal{B}_{14})\ \ =\ \ (36,\,80,\,59,\,14)
\]
whose boundary has the $f$-vector
\[
f(\boundaryOP \mathcal{B}_{14})\ \ =\ \ (36,\,68,\,34).
\]
The ball $\mathcal{B}_{14}$ again contains the cubical M\"obius strip with
parallel inner edges of length $9$ as an embedded subcomplex in its
$2$-skeleton. Compare Figure~\ref{fig:neo_3_cubical_ball}.

\paragraph{Step 4.}
Now we build the prism over this regular cubical ball,
resulting in a regular cubical $4$-ball $\mathcal{B}$ whose
$f$-vector is
\[
        f(\mathcal{B})\ \ =\ \ (72,\,196,\,198,\,87,\,14)
\]
and whose support is a cubical $4$-polytope $P_{72}:=\support{\mathcal{B}}$
with two copies of the cubical M\"obius strip in its
$2$-skeleton. Its $f$-vector is
\[
         f(P_{72})\ \ =\ \ (72,\,196,\,186,\,62).
\]
A further (computer-supported)
analysis of the dual manifolds shows that there are six dual
manifolds in total: one Klein bottle of $f$-vector $(80,\,160,\,80)$,
and five $2$-spheres (four with $f$-vector $(20,\,36,\,18)$, one with
$f$-vector $(36,\,68,\,34)$).  All the spheres are embedded, while the
Klein bottle is immersed with five double-intersection curves
(embedded $1$-spheres), but with no triple points. \hfill{\qed}



\section{Constructing cubifications}\label{sec:cubifications}                  

A lot of construction techniques for cubifications (see
Section~\ref{subsec:cubifications}) are available in the CW category.
In particular, every cubical CW $(d-1)$-sphere $\mathcal{S}^{d-1}$
with an even number of facets admits a CW cubification, that is, a
cubical CW $d$-ball with boundary $\mathcal{S}^{d-1}$,
according to Thurston \cite{Thurston1993}, Mitchell \cite{Mitchell1996}, and
Eppstein~\cite{Eppstein1999}.

\subsection{The Hexhoop template}

Yamakawa \& Shimada~\cite{YamakawaShimada2001} have introduced 
an interesting polytopal construction in dimension $3$
called the \emph{Hexhoop template}\index{Hexhoop template};
see Figure~\ref{fig:Hexhoop_3d}.

\begin{Figure} 
  \psfrag{H}{\rput(-3mm,1mm){$H$}}
  \includegraphics[width=.45\textwidth]{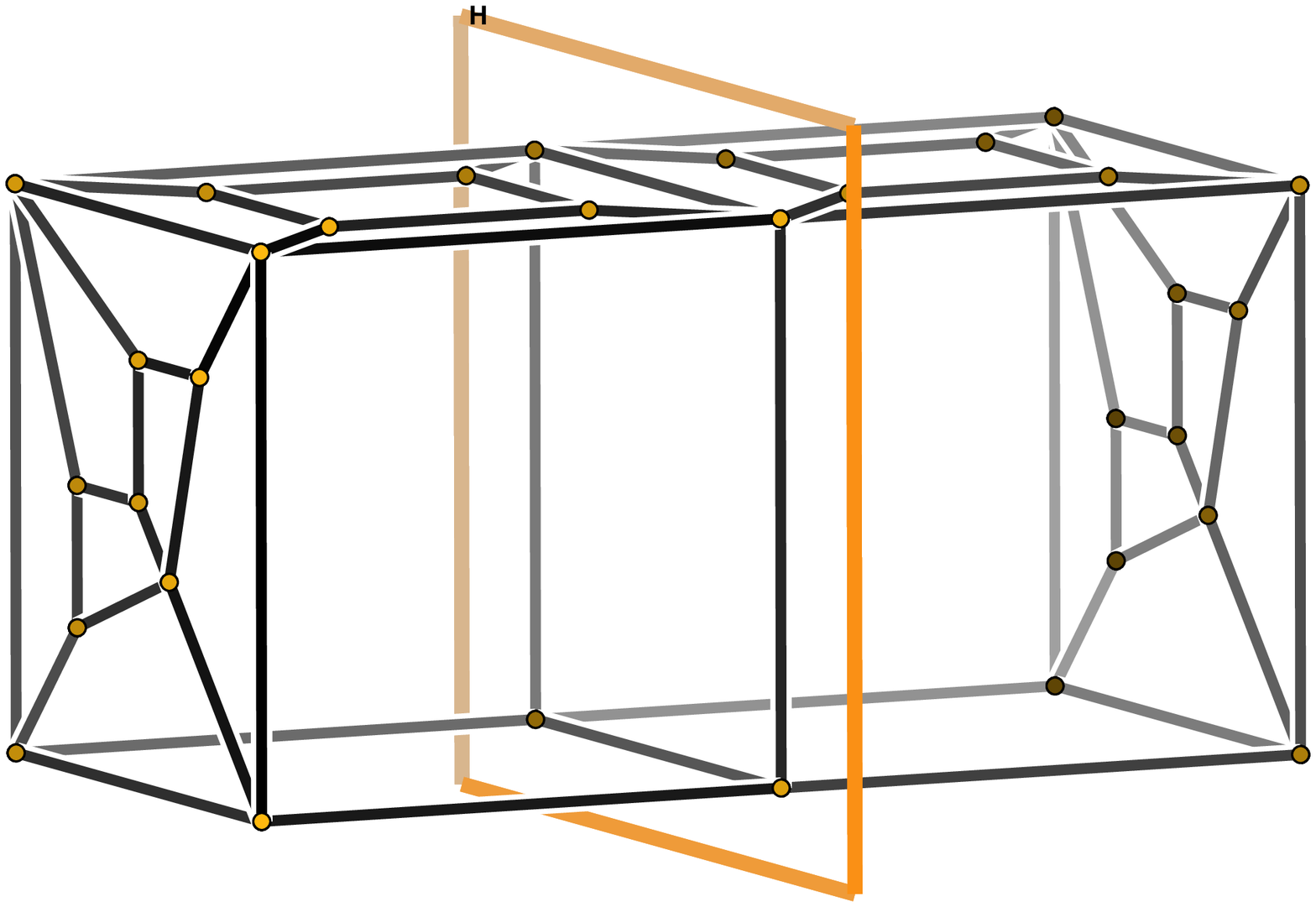}
  \rput(-2mm,4mm){$\mathcal{S}$}
  \transformsToArrow{20mm}{10mm}
  \includegraphics[width=.45\textwidth]{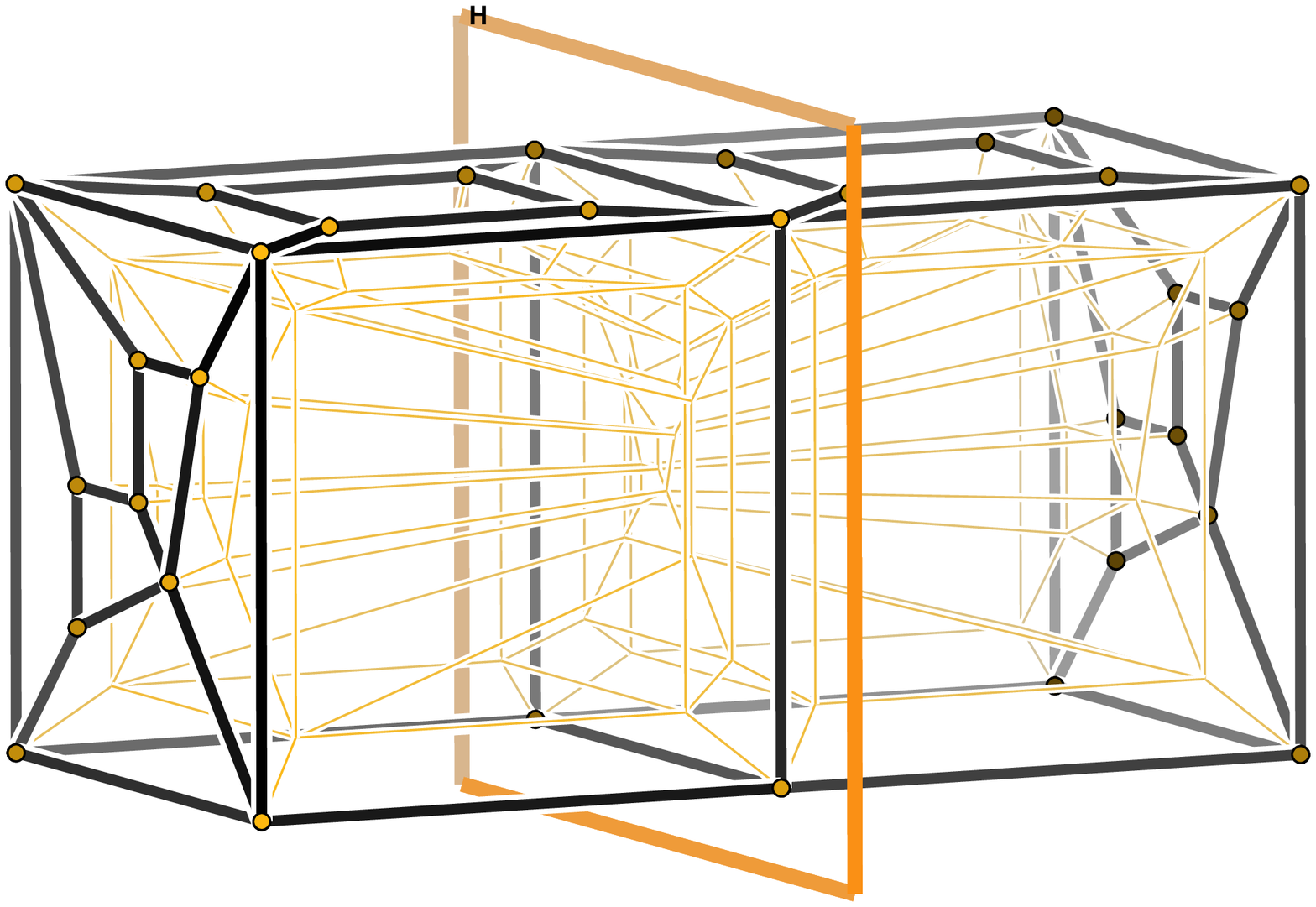}
  \rput(-2mm,4mm){$\mathcal{B}$}  \\
  \caption[The Hexhoop template]{
    The Hexhoop template of Yamakawa \&
    Shimada~\cite{YamakawaShimada2001}.}
  \label{fig:Hexhoop_3d}
\end{Figure}
Their construction takes as input a $3$-polytope $P$ that is affinely
isomorphic to a regular $3$-cube, a hyperplane $H$ and a cubical
subdivision $\mathcal{S}$ of the boundary complex of~$P$ such that
$\mathcal{S}$ is symmetric with respect to $H$ and $H$ intersects no
facet of~$\mathcal{S}$ in its relative interior. For such a
cubical PL $2$-sphere~$\mathcal{S}$ the Hexhoop template produces 
a cubification.
A $2$-dimensional version is shown in Figure~\ref{fig:Hexhoop_2d}.

\begin{FigureHere}
  \psfrag{H}{\rput(-3mm,-36mm){$H$}}
  \includegraphics[width=64mm]{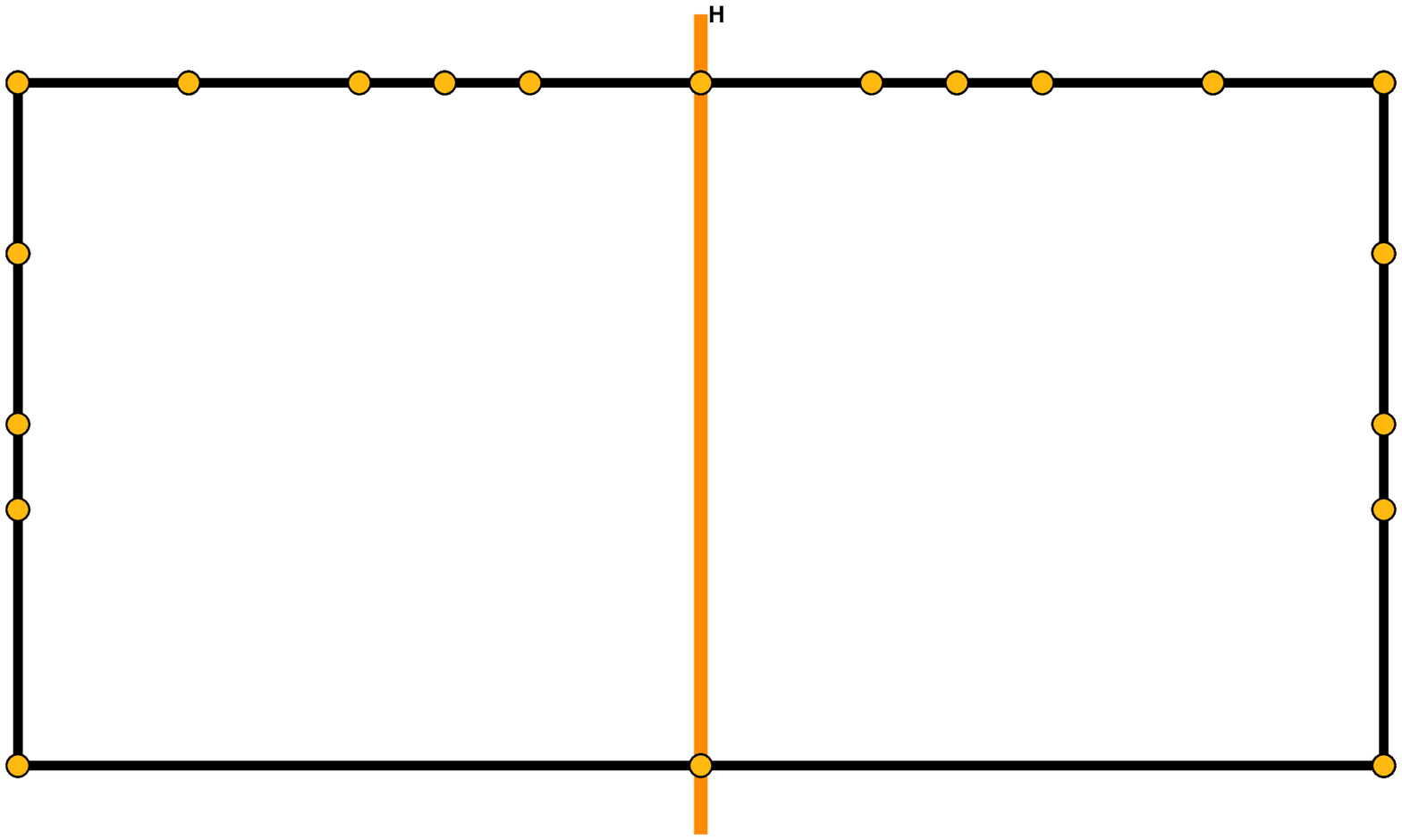}
  \rput(0mm,2mm){$\mathcal{S}$}
  \quad\transformsToArrow{15mm}{10mm}\qquad 
  \includegraphics[width=64mm]{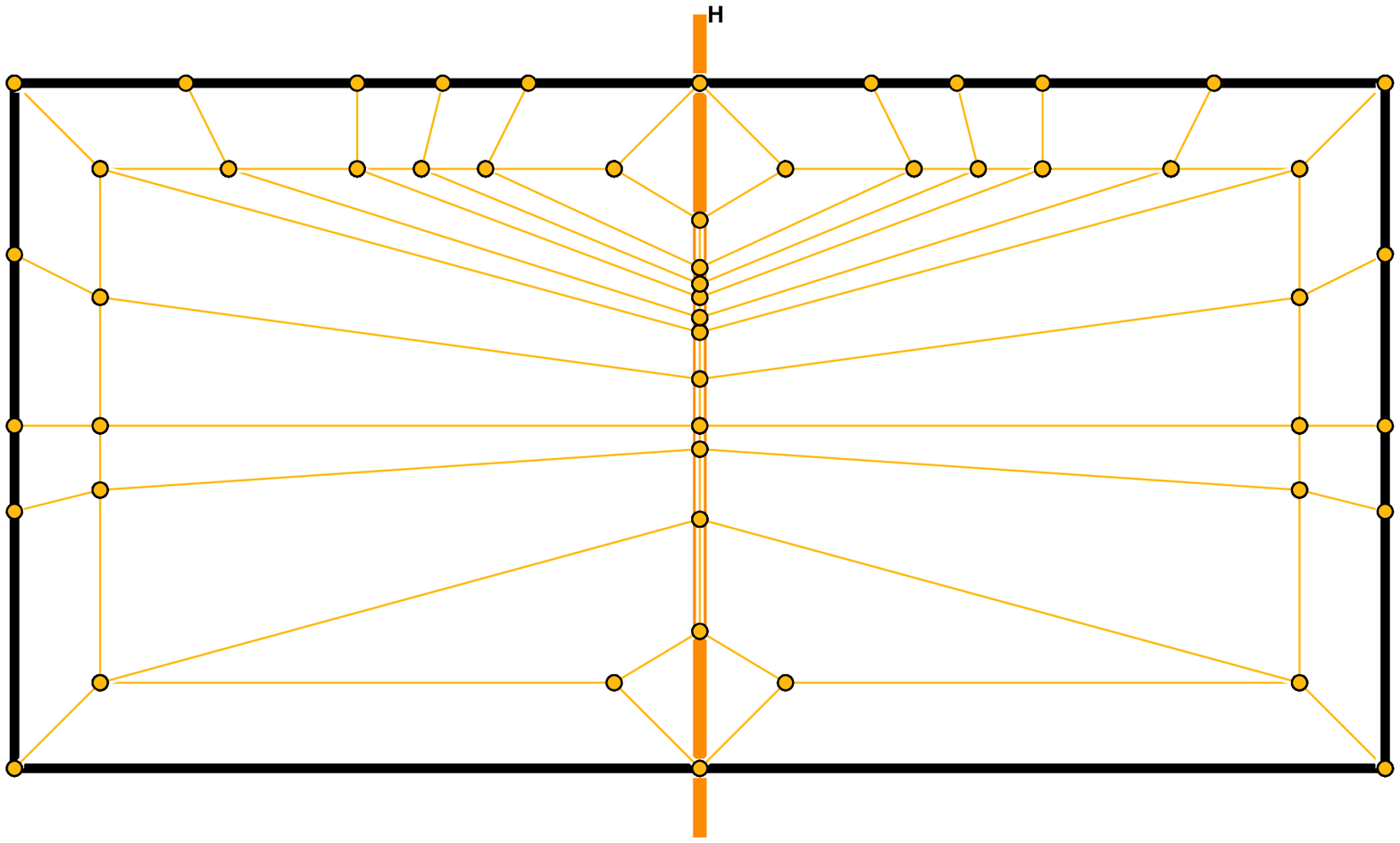}
  \rput(0mm,2mm){$\mathcal{B}$}\\
 \vspace*{2mm}
  \caption[$2$-dimensional Hexhoop]{
       A two-dimensional version of the  Hexhoop template.}
  \label{fig:Hexhoop_2d}
\end{FigureHere}

\subsection{The generalized regular Hexhoop --- overview}

In the following we present a 
\emph{generalized regular Hexhoop} construction.  It is a
generalization of the Hexhoop template in several directions: Our
approach admits arbitrary geometries, works in any dimension, and
yields regular cubifications with ``prescribed heights on the
boundary'' (with a symmetry requirement and with the requirement that
the intersection of the symmetry hyperplane and
the boundary subdivision is a subcomplex of the boundary subdivision).\\
Figure~\ref{fig:gen_reg_Hexhoop_outcome} displays a $2$-dimensional
cubification (of a boundary subdivision $\mathcal{S}$ of a
$2$-polytope such that $\mathcal{S}$ is symmetric with respect to a
hyperplane $H$) obtained by our construction.
\begin{Figure}
  \psfrag{H}{\rput(2mm,2mm){$H$}}
  \psfrag{S}{\rput(-3mm,0mm){$\mathcal{S}$}}
  \includegraphics[width=.33\textwidth]{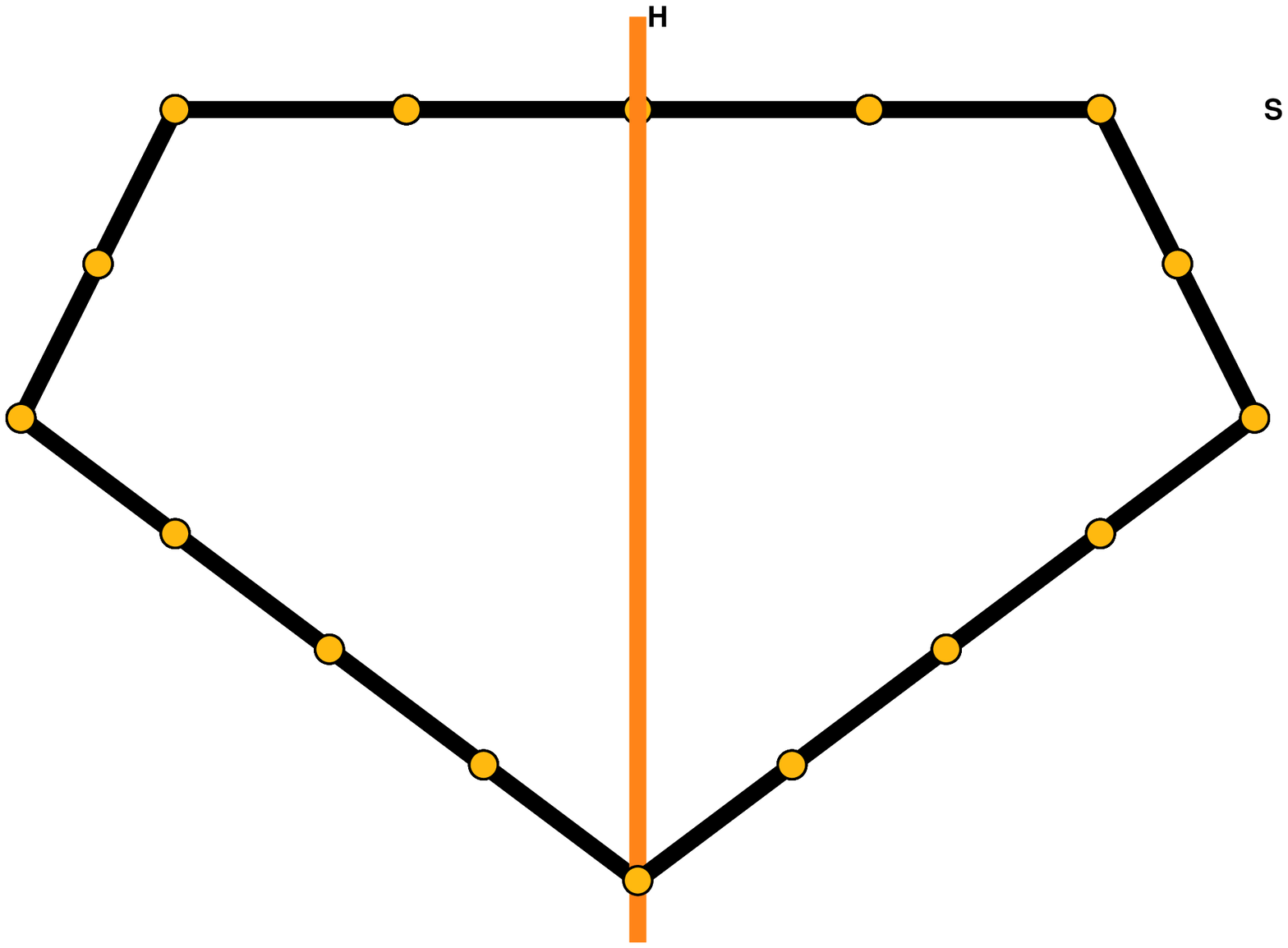}
 \quad \transformsToArrow{20mm}{10mm}\qquad
  \includegraphics[width=.33\textwidth]{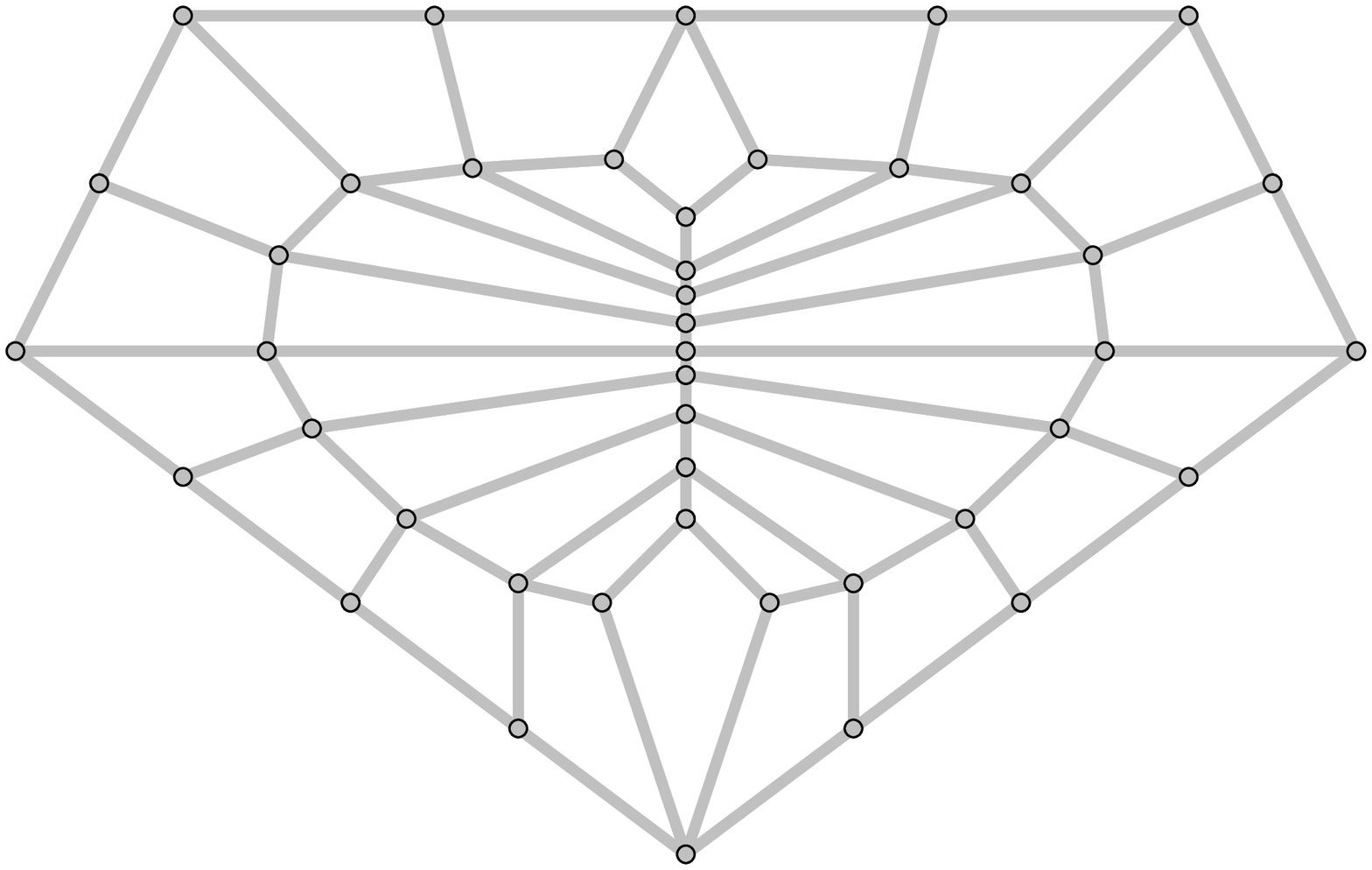}
  \caption[Outcome of the generalized regular Hexhoop]{
    A cubification of a boundary subdivision of a pentagon, produced
    by our \emph{generalized regular Hexhoop} construction.}
  \label{fig:gen_reg_Hexhoop_outcome}
\end{Figure}

Not only do we get a cubification, but we may also derive a symmetric
lifting function for the cubification that may be quite arbitrarily
prescribed on the boundary.
The input of our construction is a lifted cubical boundary subdivision
$(\mathcal{S}^{d-1},\psi)$ of a $d$-polytope~$P$, such that both $P$
and $(\mathcal{S}^{d-1},\psi)$ are symmetric with respect to a
hyperplane~$H$.
\begin{Figure} 
  \psfrag{P}{\rput(0mm,0mm){$P$}}
  \psfrag{S}{\rput(0mm,0mm){$\mathcal{S}$}}
  \psfrag{H}{\rput(1mm,2mm){$H$}}
  \psfrag{lift}{\rput(-3mm,0mm){$\lifted{\mathcal{S}}{\psi}$}}
  \includegraphics[width=.5\textwidth]{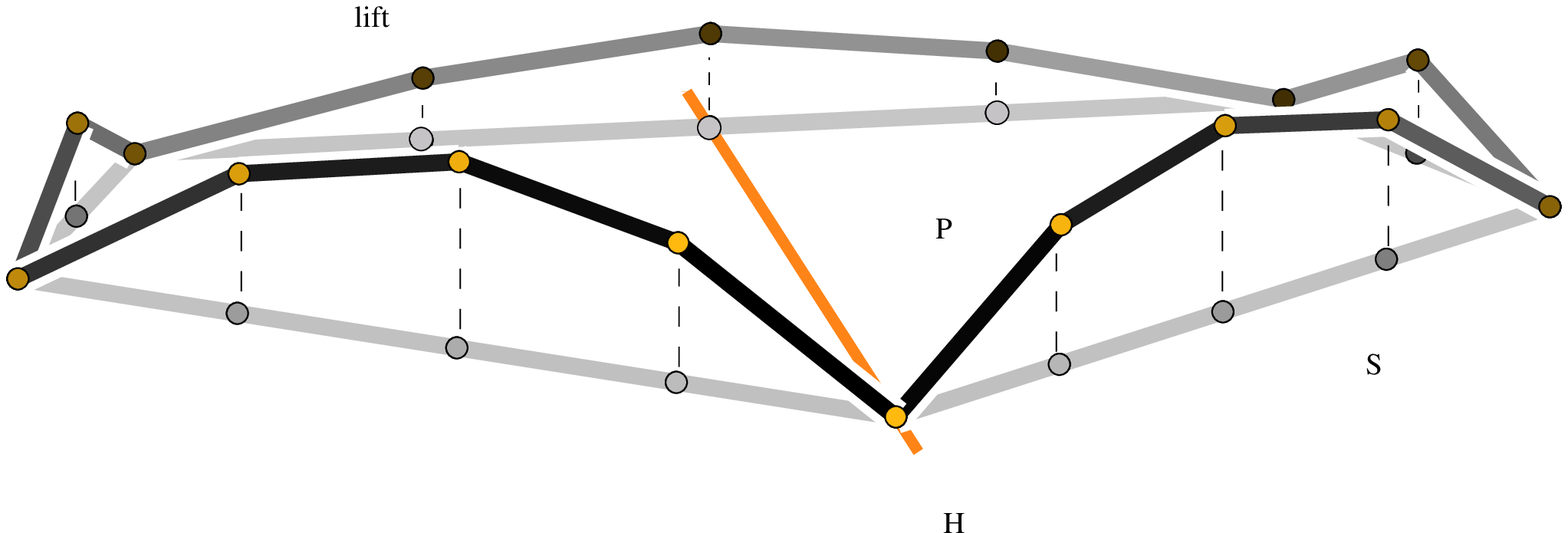}
  \caption[Input of generalized regular Hexhoop]{An
    input for the generalized regular Hexhoop construction.}
  \label{fig:input_of_gen_Hexhoop}
\end{Figure}

Our approach goes roughly as follows.
\begin{compactsteps}
\item 
  We first produce a $(d+1)$-polytope $T$ that is a \emph{symmetric tent}
  (defined in Section~\ref{subsec:tent_over_lifted_bd_subdiv}) over
  the given lifted boundary subdivision $(\mathcal{S},\psi)$ of the
  input $d$-polytope $P$. Such a tent is the convex hull of all
  `lifted vertices'
  $(\boldsymbol{v},\psi(\boldsymbol{v}))\in\R^{d+1}$,
  $\boldsymbol{v}\in\vertices{\mathcal{S}}$, and of two 
    \emph{apex points} $\boldsymbol{p}_L,\boldsymbol{p}_R$; compare
  Figure~\ref{fig:first_tent}.

   \begin{Figure} \bigskip
      \psfrag{H}{\rput(13.5mm,-18mm){$H$}}
      \psfrag{T}{\rput(-3mm,-20mm){$T$}}
      \psfrag{1}{\rput(-2mm,4mm){$\boldsymbol{p}_L$}}
      \psfrag{0}{\rput(1mm,4mm){$\boldsymbol{p}_R$}}
     \includegraphics[width=77mm]{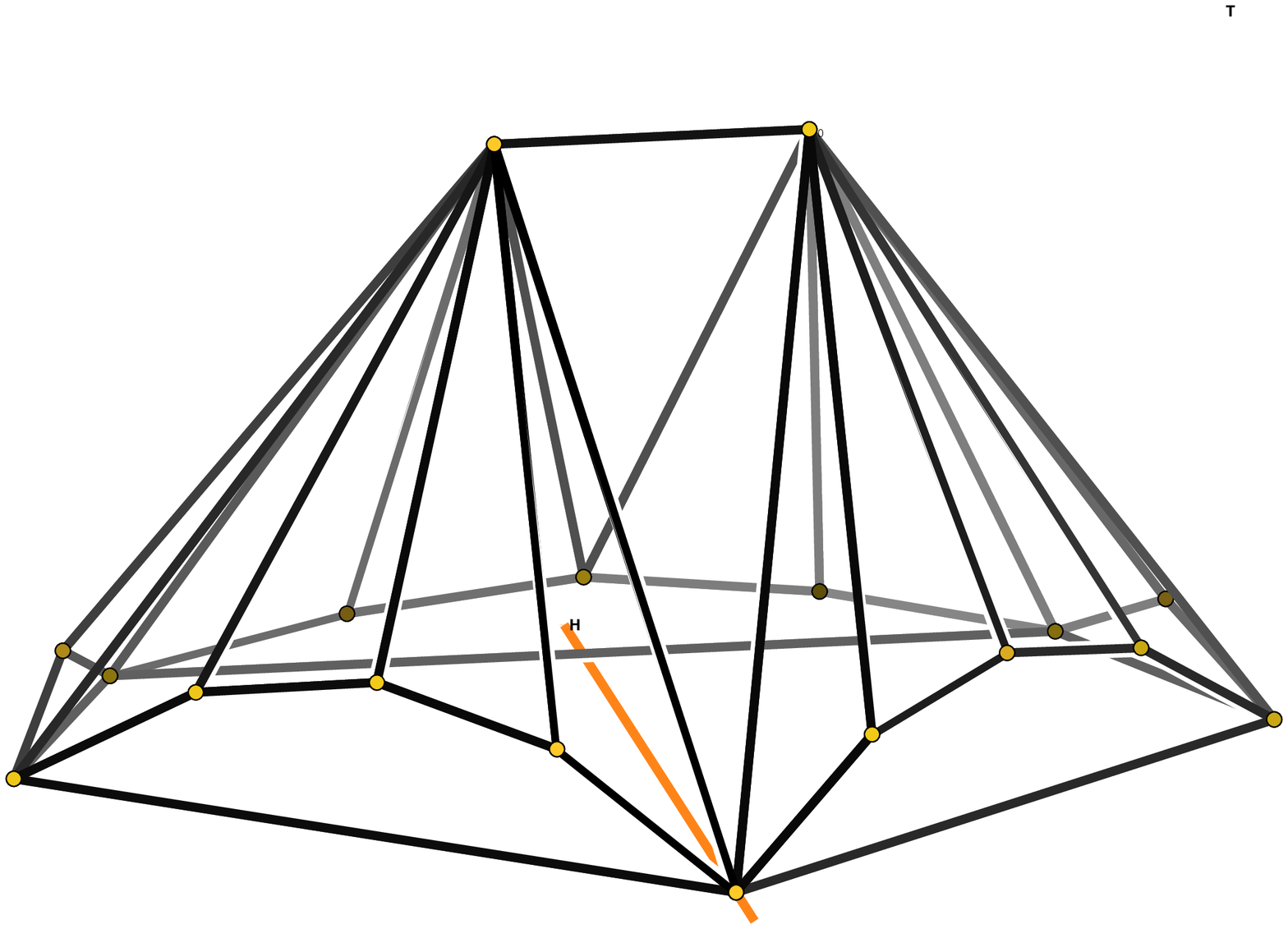}%
      \caption[A symmetric tent over a lifted boundary subdivision]{%
        A symmetric tent over the lifted boundary subdivision
        $(\mathcal{S},\psi)$ of the input $d$-polytope $P$.}
      \label{fig:first_tent}
    \end{Figure}
    
  \item Truncate $T$ by a hyperplane $H'$ parallel to
    $\affOP(P)=\R^d\subset\R^d\times\{0\}$ that separates the lifted
    points from the apex points, and remove the upper part.
    
  \item Add the polytope
    $R:=\coneOP(\boldsymbol{p}_L,Q)\cap\coneOP(\boldsymbol{p}_R,Q)\cap
    H'_{+}$, where $H'_{+}$ is the halfspace with respect to $H$ that
    contains $\boldsymbol{p}_L$ and $\boldsymbol{p}_R$. Compare
    Figure~\ref{fig:sketch_of_gen_reg_hexhoop}.
  \item Project the upper boundary complex of the resulting polytope to~$\R^d$.

\end{compactsteps}

The figures in this section illustrate the generalized regular
Hexhoop construction for the
$2$-dimensional input polytope of Figure~\ref{fig:input_of_gen_Hexhoop};
the generalized Hexhoop construction for $d=2$
yields $2$-dimensional complexes in $\R^3$. 
The extension to higher dimensions is immediate, and the case $d=3$ is
crucial for us (see Section~\ref{sec:immersions}). 
It is, however, also harder to visualize:
A $3$-dimensional generalized regular Hexhoop cubification is shown 
in Figure~\ref{fig:dmf_three_dim_gen_Hexhoop}. 

\begin{Figure}
    \psfrag{H}{\rput(2mm,1mm){$H$}}
    \psfrag{H'}{\rput(4mm,0mm){$H'$}}
    \psfrag{H+R}{\rput(11mm,-.5mm){$H+\R\unitvector{d+1}$}}
    \psfrag{1}{\rput(-3.5mm,4.3mm){$\boldsymbol{p}_L$}}
    \psfrag{0}{\rput(3mm,3.5mm){$\boldsymbol{p}_R$}}

    \includegraphics[width=\textwidth]{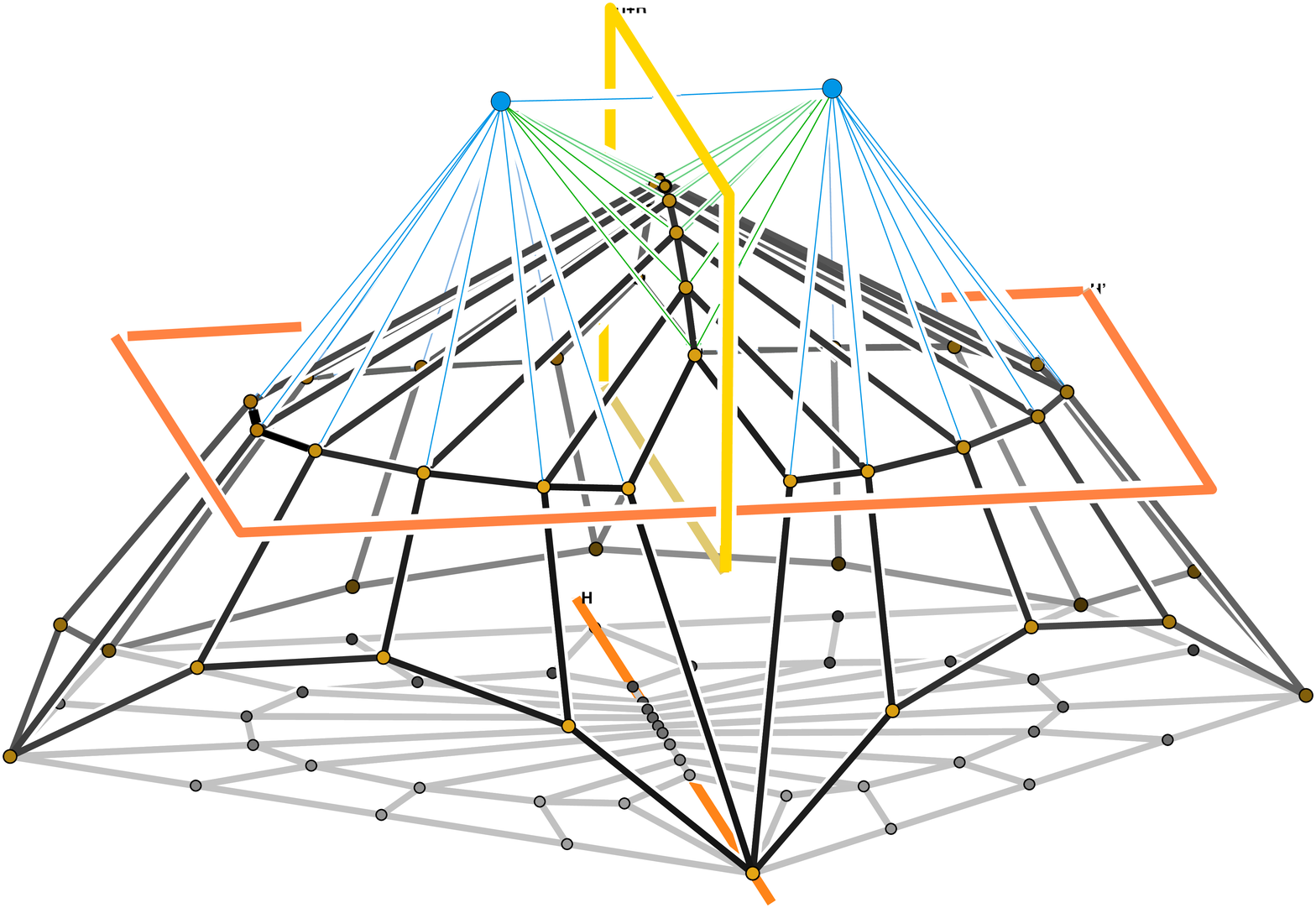}
    \caption[Sketch of the generalized regular Hexhoop construction]{
      Sketch of the generalized regular Hexhoop construction.}
              
    \label{fig:sketch_of_gen_reg_hexhoop}
\end{Figure}


\subsection{Symmetric tent over a lifted boundary subdivision}
\label{subsec:tent_over_lifted_bd_subdiv}

Let $P$ be a $d$-polytope that is symmetric with respect to a
hyperplane $H$ in $\R^{d}$. Choose a positive halfspace $H_{+}$ with
respect to $H$.  Let $(\mathcal{S},\psi)$ be a lifted boundary
subdivision of $P$ such that $\mathcal{S}\cap H$ is a subcomplex
of~$\mathcal{S}$. Define $\widetilde H:=H+\R\unitvector{d+1}$, which is
a symmetry hyperplane for~$P\subset\R^{d+1}$.
The positive  halfspace of $\widetilde{H}\subset\R^{d+1}$
will be denoted by $\widetilde{H}_+$.

The \emph{symmetric tent} over $(\mathcal{S},\psi)$ is the lifted
polytopal subdivision $(\mathcal{T},\phi)$ of $P$ given by the
upper faces of the polytope
\[   
    T\ :=\ \conv(P \cup \{\boldsymbol{p}_L,\boldsymbol{p}_R\})
\]
if $\boldsymbol{p}_L,\boldsymbol{p}_R\in\R^{d+1}$ are two apex points in
$\R^{d+1}$ that are symmetric with respect to the hyperplane
$\widetilde H$,
and the upper facets of $T$ are
\begin{compactitem}
  \item pyramids with apex point $\boldsymbol{p}_L$ over facets $F$ of
        $\lifted{\mathcal{S}}{\psi}$ such that $\pi(F)\subset H_{+}$,
  \item pyramids with apex point $\boldsymbol{p}_R$ over facets $F$ of
        $\lifted{\mathcal{S}}{\psi}$ such that $\pi(F)\subset H_{-}$, and
      \item $2$-fold pyramids with apex points
        $\boldsymbol{p}_L,\boldsymbol{p}_R$ over ridges $R$ of
        $\lifted{\mathcal{S}}{\psi}$ with $\pi(R)\subset H$.
\end{compactitem}

(This requires that $\boldsymbol{p}_L\not\in\affOP(P)$
 and $\pi(\boldsymbol{p}_L)\in\relint{P\cap H_{+}}$.)

\begin{lemma}\label{lemma:ex_of_tent}
    Assume we are given the following input.
    \begin{compactsymboldesc}[$\boldsymbol{q}_L,\boldsymbol{q}_R\in P$]

        \item[$P$]       a convex $d$-polytope in $\R^d$,

        \item[$(\mathcal{S},\psi)$] 
                         a lifted boundary subdivision of $P$,

        \item[$H$]       a hyperplane  in~$\R^d$  such that 
                  \begin{compactitem}
                     \item $P$ and  $(\mathcal{S},\psi)$ 
                               are both symmetric with respect to $H$, and
                     \item $\mathcal{S}\cap H$ is a 
                           subcomplex of $\mathcal{S}$, and
                  \end{compactitem}
        \item[$\boldsymbol{q}_L,\boldsymbol{q}_R$]
                           two points in $P\subset \R^{d}$ such that
                  \begin{compactitem}
                     \item $\boldsymbol{q}_L\in\relint{P\cap H_{+}}$, and
                     \item $\boldsymbol{q}_L,\boldsymbol{q}_R$ are 
                           symmetric with respect to $\widetilde H$.
                  \end{compactitem}
    \end{compactsymboldesc}
   
    Then for every sufficiently large height $h>0$ the $(d+1)$-polytope
    $T:=\conv\{\lifted{\mathcal{S}}{\psi},\boldsymbol{p}_L,\boldsymbol{p}_R\}$
    with $\boldsymbol{p}_L:=(\boldsymbol{q}_L,h)\in\widetilde {H}_+$ and
    $\boldsymbol{p}_R:=(\boldsymbol{q}_R,h)\notin\widetilde {H}_+$ 
    is a symmetric tent over $(\mathcal{S},\psi)$.
\end{lemma}

This can be shown for instance by using the Patching Lemma
(Lemma~\ref{lemma:patching_lemma}).

\subsection{The generalized regular Hexhoop in detail}%
\label{subsec:generalized_hexhoop}

In this section we specify our generalization of the 
Hexhoop template and prove
the following existence statement for cubifications.

\begin{theorem}\label{thm:ex_of_sym_cubification}
Assume we are given the following input.
    \begin{compactsymboldesc}[$(\mathcal{S}^{d-1},\psi)$]
    \item[$P$]       a convex $d$-polytope in $\R^d$,
    \item[$(\mathcal{S}^{d-1},\psi)$] 
                     a lifted cubical boundary subdivision of $P$, and
    \item[$H$]       a hyperplane  in~$\R^d$ such that 
      \begin{compactitem}
          \item $P$ and  $(\mathcal{S}^{d-1},\psi)$ are symmetric 
                with respect to $H$, and
          \item $\mathcal{S}^{d-1}\cap H$ is a 
                subcomplex of~$\mathcal{S}^{d-1}$.
      \end{compactitem}
   \end{compactsymboldesc}
Then there is a lifted cubification $(\mathcal{B}^d,\phi)$ 
of~$(\mathcal{S}^{d-1},\psi)$.
\end{theorem}

The proof relies on the following construction.

\begin{construction}{Generalized regular Hexhoop}
  \label{constr:gen_hexhoop}
  \begin{compactdesc}
  \item[Input:]\
    \begin{compactsymboldesc}[$(\mathcal{S}^{d-1},\psi)$]
    \item[$P$]       a convex $d$-polytope $P$ in $\R^d$.
    \item[$(\mathcal{S}^{d-1},\psi)$] 
                     a lifted cubical boundary subdivision of $P$.
    \item[$H$]       a hyperplane  in~$\R^d$ such that 
      \begin{compactitem}
          \item $P$ and  $(\mathcal{S}^{d-1},\psi)$ 
                are symmetric with respect to $H$, and
          \item $\mathcal{S}^{d-1}\cap H$ is a subcomplex 
                of~$\mathcal{S}^{d-1}$.
      \end{compactitem}
      
    \end{compactsymboldesc}

  \item[Output:] \ 
    \begin{symboldesc}[$(\mathcal{S}^{d-1},\psi)$]
        \item[$(\mathcal{B}^d,\phi)$] a symmetric lifted cubification 
           of $(\mathcal{S}^{d-1},\psi)$ given by                                      a cubical $d$-ball $\mathcal{C}'$ in $\R^{d+1}$.
    \end{symboldesc}

  \end{compactdesc}

  \begin{steps}
  
  \item Choose a positive halfspace $H_+$ with respect to $H$, and a
    point $\boldsymbol{q}_L \in \relint{P\cap H_{+}}$.  Define
    $\boldsymbol{q}_{R}:=\boldsymbol{p}_{L}^M$, where the upper index
    $^M$ denotes the mirrored copy with respect to
    $\widetilde H = H+\R\unitvector{d+1}$.

      By Lemma~\ref{lemma:ex_of_tent} there is a height $h>0$ such that
    \[T:=\conv\{\lifted{\mathcal{S}^{d-1}}{\psi},\boldsymbol{p}_L,\boldsymbol{p}_R\}\] with $\boldsymbol{p}_L:=(\boldsymbol{q}_L,h)$ and
    $\boldsymbol{p}_R:=(\boldsymbol{q}_R,h)$ 
    forms a symmetric tent over $(\mathcal{S}^{d-1},\psi)$.

  \item Choose a hyperplane $H'$ parallel to $\affOP(P)\subset\R^d$
    that separates $\{\boldsymbol{p}_L,\boldsymbol{p}_R\}$ and
    $\lifted{\mathcal{S}^{d-1}}{\psi}$.  Let $H'_{+}$ be the halfspace
    with respect to~$H'$ that contains $\boldsymbol{p}_L$
    and~$\boldsymbol{p}_R$.
\smallskip

    \begin{Figure}
      \psfrag{H}{\rput(11mm,-16mm){$H$}}
      \psfrag{H'}{\rput(4mm,0mm){$H'$}}
      \psfrag{T}{\rput(-12mm,-10mm){$T$}}
      \psfrag{1}{\rput(-3.5mm,2mm){$\boldsymbol{p}_L$}}
      \psfrag{0}{\rput(3mm,2mm){$\boldsymbol{p}_R$}}
      \includegraphics[width=.4\textwidth]{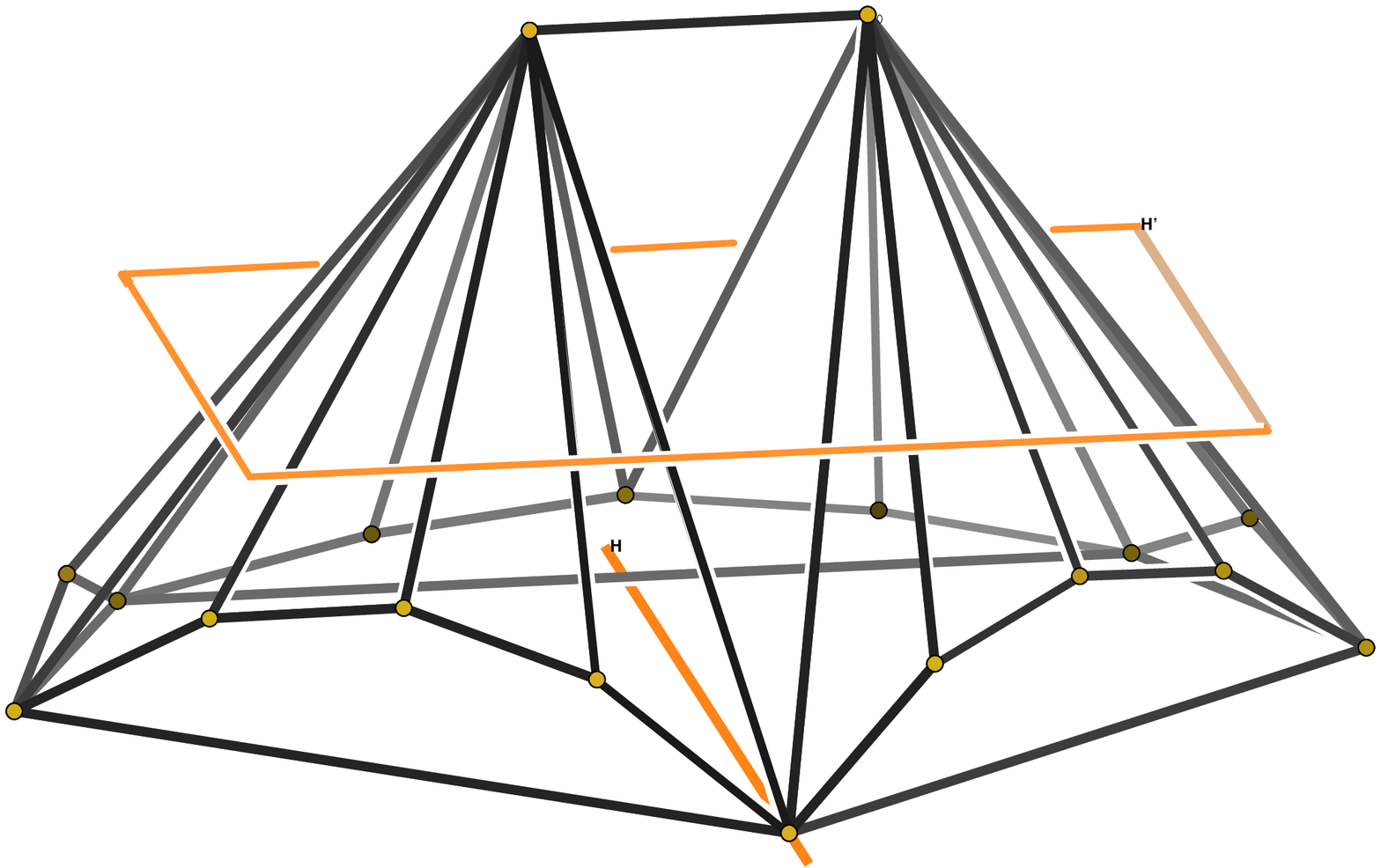}
      \caption[Step 2 of generalized regular Hexhoop construction]{%
        Step 2. The hyperplane $H'$ 
        separates $\{\boldsymbol{p}_L,\boldsymbol{p}_R\}$ from
        $\lifted{\mathcal{S}^{d-1}}{\psi}$.  }
    \end{Figure}

   \item Define the ``lower half'' of the tent $T$ as
\[
              T_-\ :=\ T \cap H'_{-},
\]
           whose ``top facet'' is the convex $d$-polytope $Q:=T\cap H'$.

   \begin{Figure} \vspace*{-5mm}
      \psfrag{H'}{\rput(6mm,0mm){$H'$}}
      \psfrag{H}{}
      \psfrag{T}{\rput(-15mm,-3mm){$T$}}
      \psfrag{T_}{\rput(-2mm,0mm){$T_{-}$}}
      \psfrag{Q}{\rput(1.5mm,.5mm){$Q$}}
      \includegraphics[width=.4\textwidth]{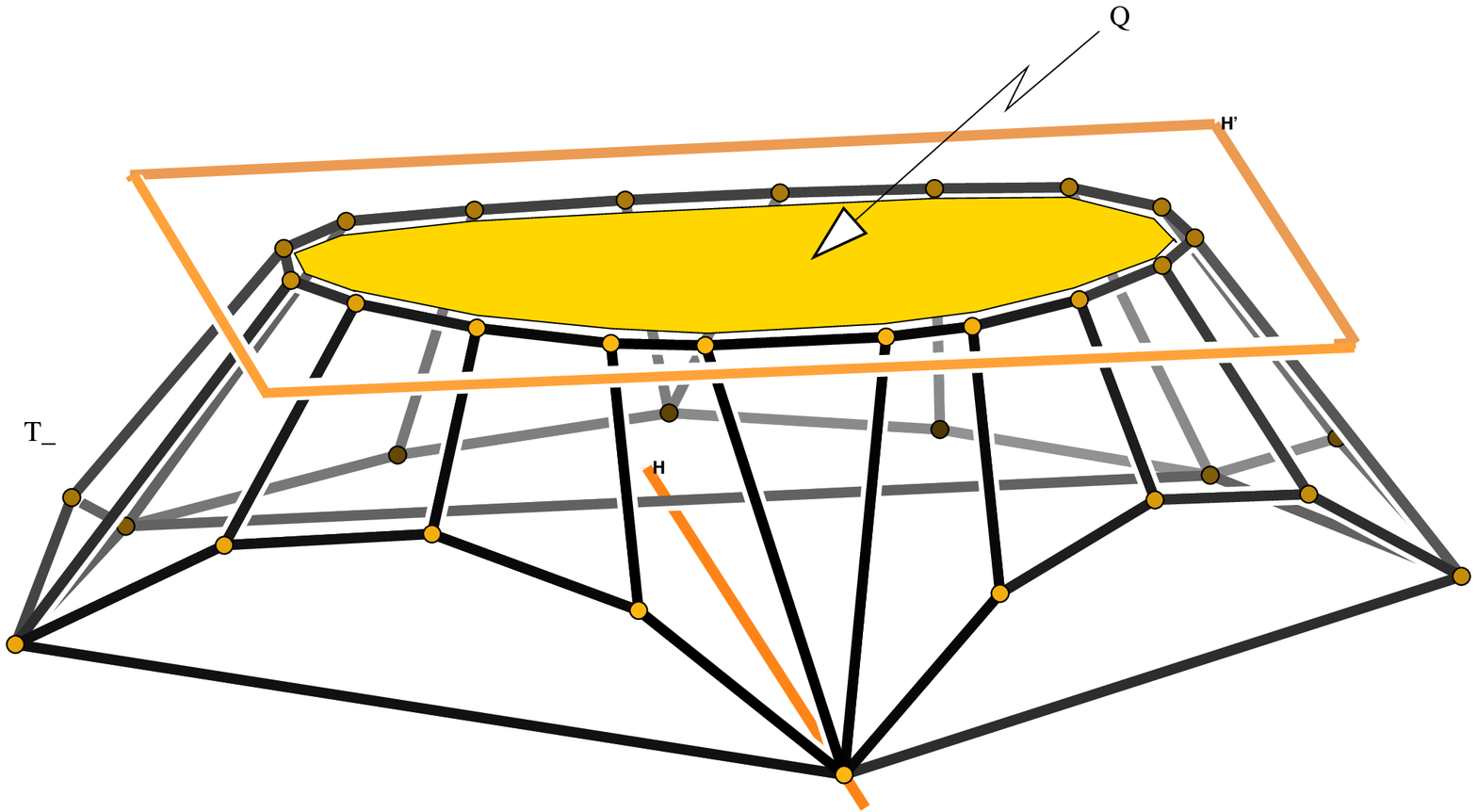}
      \caption[Step 3 of generalized regular Hexhoop construction]{%
        Step 3.
        The ``lower half'' $T_{-}$ of $T$.}
    \end{Figure}

   \item Define the two $d$-polytopes 
           \begin{eqnarray*}
                Q_L &:=& \conv{\setdef{\boldsymbol{v} \in\vertices{Q}}{\boldsymbol{v} \in H_{+}}},\\
                Q_R &:=& \conv{\setdef{\boldsymbol{v} \in\vertices{Q}}{\boldsymbol{v} \in H_{-}}}.
           \end{eqnarray*}
           Let $F_L:=H'\cap\convOP(\boldsymbol{p}_L,P\cap H)$,
           the unique facet of $Q_L$ that is not a facet of~$Q$.

   \begin{Figure} 
      \psfrag{H'}{\rput(6mm,0mm){$H'$}}
      \psfrag{H}{}
      \psfrag{T}{\rput(-15mm,-3mm){$T$}}
      \psfrag{T_}{\rput(-2mm,0mm){$T_{-}$}}
      \psfrag{Q_L}{\rput(1mm,1.5mm){$Q_L$}}
      \psfrag{Q_R}{\rput(3mm,2mm){$Q_R$}}
      \psfrag{F_L}{\rput(2mm,1.5mm){$F_L$}}
      \includegraphics[width=.4\textwidth]{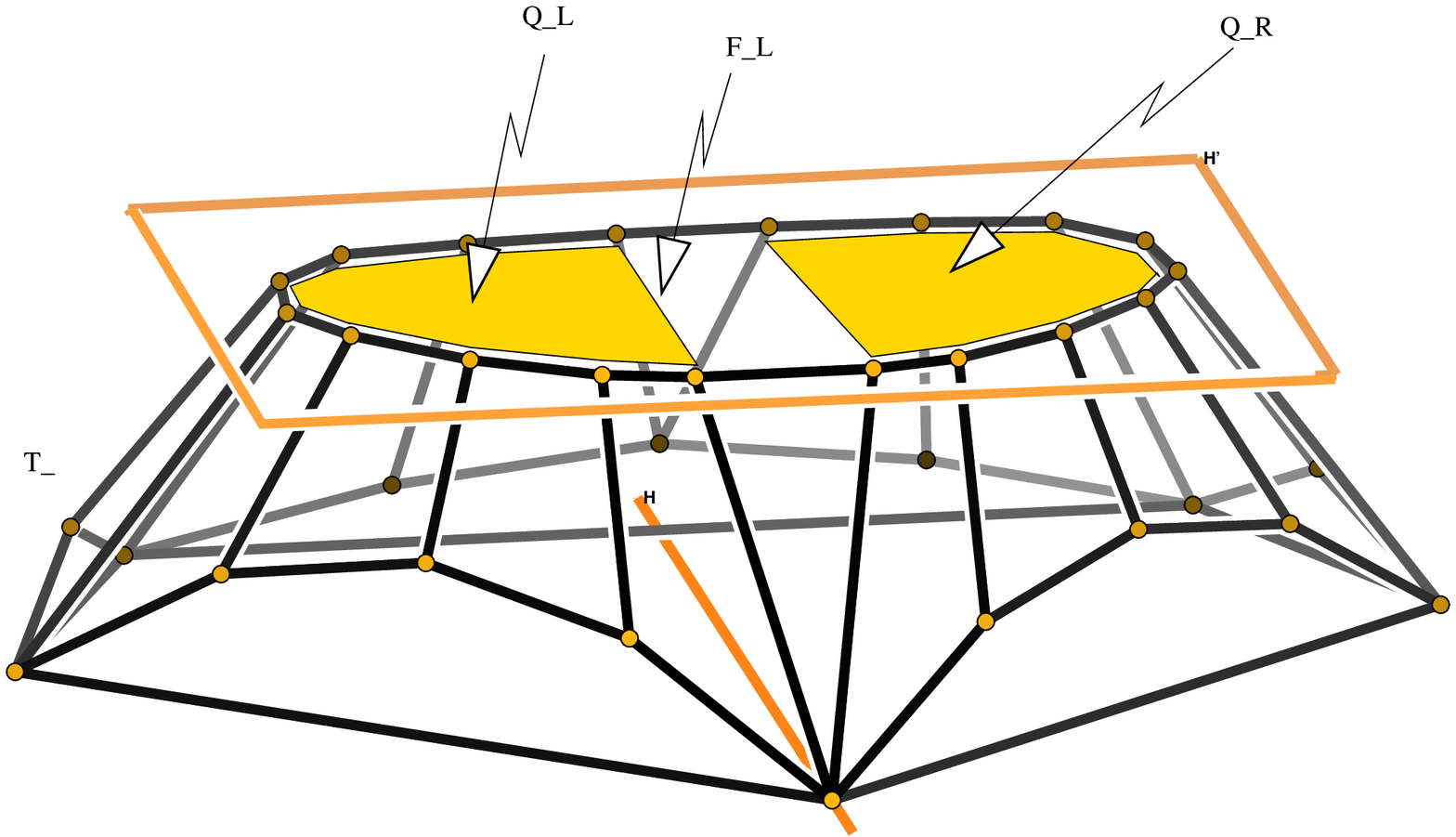}
      \caption[Step 4 of generalized regular Hexhoop construction]{%
        Step 4.
        Define $Q_L$ and $Q_R$.}
    \end{Figure}
           
  \item Construct the polytope
\[
R\ \ :=\ \ 
\coneOP(\boldsymbol{p}_L,Q)\ \cap\ \coneOP(\boldsymbol{p}_R,Q)\ \cap\ H'_{+}.
\]

   \begin{Figure} \vspace*{2mm}
       \psfrag{p_L}{\rput(2mm,0mm){$\boldsymbol{p}_L$}}
       \psfrag{p_R}{\rput(2mm,0mm){$\boldsymbol{p}_R$}}
      \psfrag{H'}{\rput(4mm,0mm){$H'$}}
      \psfrag{H}{\rput(6mm,0mm){$H$}}
      \psfrag{H''}{\rput(2mm,3.5mm){$\widetilde H = H+\R\unitvector{d+1}$}}
      \psfrag{R}{\rput(0mm,0mm){$R$}}
      \includegraphics[width=.7\textwidth]{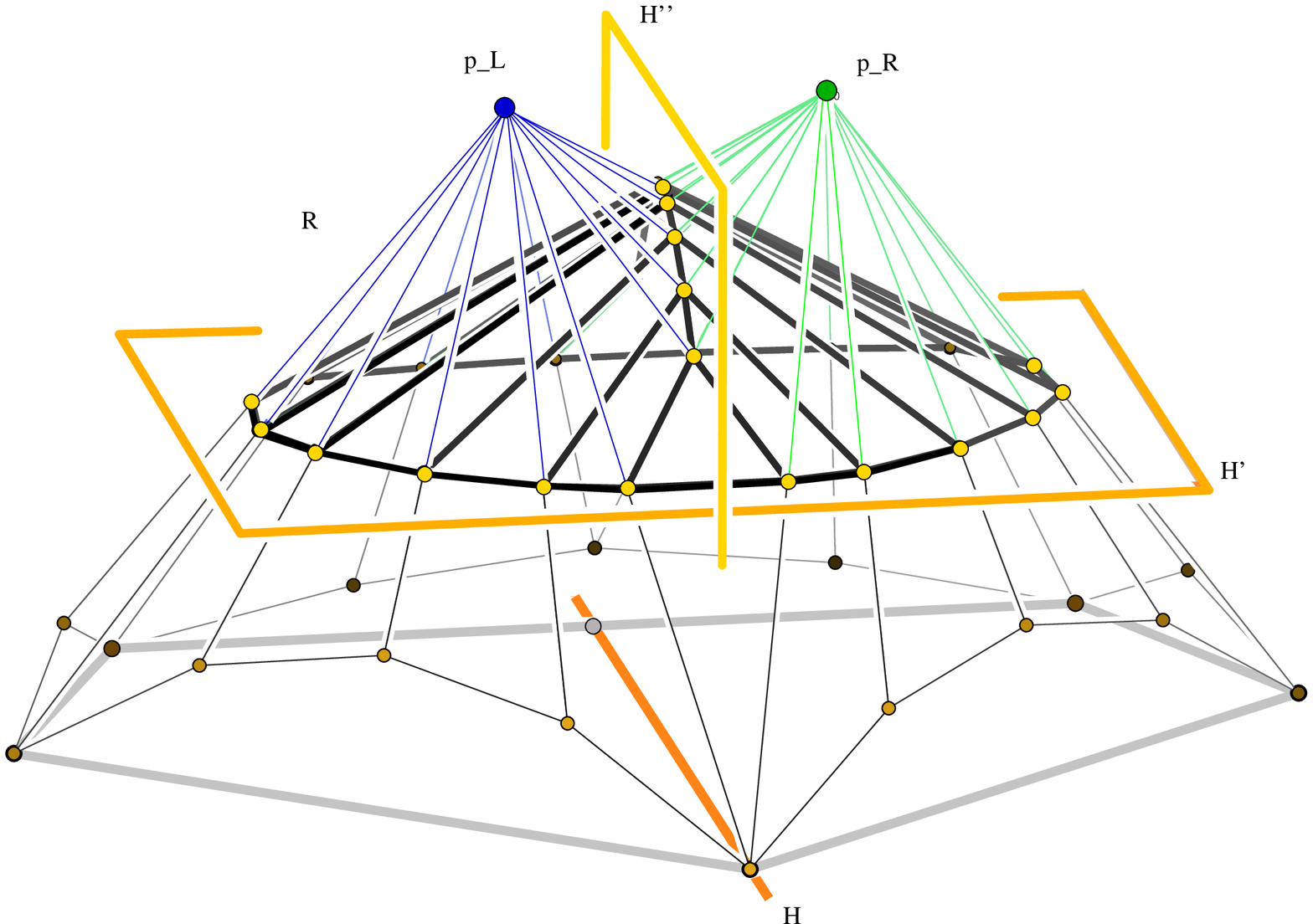}
      \caption[Step 5 of generalized regular Hexhoop construction]{%
        Step 5. 
        The polytope $R:=\coneOP(\boldsymbol{p}_L,Q)\cap\coneOP(\boldsymbol{p}_R,Q)\cap
    H'_{+}$.}
    \end{Figure}
         
  \end{steps}
The complex $\mathcal{C}'$ in question is given by the upper facets of
the $(d+1)$-polytope 
\[
U\ :=\ T_{-} \cup R.
\]
\vskip-3mm
\end{construction}

\begin{lemma}[Combinatorial structure of \boldmath$Q$]~\\
   The vertex set of $Q$ consists of
   \begin{compactitem}[~~$\bullet$~]
   \item the points $\convOP(\boldsymbol{p}_L,\boldsymbol{v})\cap H'$
     for vertices
     $\boldsymbol{v}\in\vertices{\lifted{\mathcal{S}}{\psi}}$ such
     that $\pi(\boldsymbol{v})\subset H_+$, and
      \item the points $\convOP(\boldsymbol{p}_R,\boldsymbol{v})\cap
        H'$ for vertices
        $\boldsymbol{v}\in\vertices{\lifted{\mathcal{S}}{\psi}}$ such
        that $\pi(\boldsymbol{v})\subset H_-$.
   \end{compactitem}
   The facets of $Q$ are
   \begin{compactenum}[\ \rm(a)]
       \item the combinatorial cubes 
             $\convOP(\boldsymbol{p}_L,F)\cap H'$ 
             for facets $F$ of $\lifted{\mathcal{S}}{\psi}$ such that $F\subset \widetilde H_{+}$,
       \item the combinatorial cubes 
             $\convOP(\boldsymbol{p}_R,F)\cap H'$ 
             for facets $F$ of $\lifted{\mathcal{S}}{\psi}$ such that $F\subset \widetilde H_{-}$,
             
           \item the combinatorial cubes
             $\convOP(\boldsymbol{p}_L,\boldsymbol{p}_R,F)\cap H'$ for
             $(d-2)$-faces $F$ of $\lifted{\mathcal{S}}{\psi}$ with
             $F\subset \widetilde H$.
             
  \end{compactenum}
\end{lemma}
\begin{proof}
  By the definition of a symmetric tent, upper facets of the
  symmetric tent $T$ are
  \begin{compactitem}[~~$\bullet$~]
    \item the pyramids with apex point $\boldsymbol{p}_L$ over facets $F$ of
        $\lifted{\mathcal{S}}{\psi}$ such that $F \subset \widetilde H_{+}$,
    \item the pyramids with apex point $\boldsymbol{p}_R$ over facets $F$ of
        $\lifted{\mathcal{S}}{\psi}$ such that $F\subset \widetilde H_{-}$, and
    \item the $2$-fold pyramids with apex points
        $\boldsymbol{p}_L,\boldsymbol{p}_R$ over ridges $R$ of
        $\lifted{\mathcal{S}}{\psi}$ with $R\subset \widetilde H$.
  \end{compactitem}
  Since $Q$ is the intersection of $T$ with $H$, the polytope $Q$ 
  has the vertices and facets listed above.  
  It remains to show that the facets of type (c) are combinatorial
  cubes.  Let $F$ be a $(d-2)$-face of $\lifted{\mathcal{S}}{\psi}$
  such that $F\subset \widetilde H$. 
  Every point on the facet lies in the convex hull of $F$ with a
  unique point on the segment $[\boldsymbol{p}_L,\boldsymbol{p}_R]$. 
  Thus the facet is combinatorially isomorphic to a prism over~$F$.
\end{proof}


Let a \emph{$d$-dimensional half-cube} be the product of a
combinatorial $(d-2)$-cube and a triangle. A \emph{combinatorial
  half-cube} is a polytope combinatorially isomorphic to a half-cube.

\begin{lemma}[Combinatorial structure of \boldmath$T_{-}$]~\\
  The vertices of $T_{-}$ are the vertices of $\lifted{\mathcal{S}}{\psi}$ and the
  vertices of~$Q$.
  Furthermore, the upper facets of $T_{-}$ are
  \begin{compactenum}[\ \rm(a)]
    \item the combinatorial cubes 
               $\coneOP(\boldsymbol{p}_L,F)\cap H'_{-}\cap(\R^d\times\R_{+})$ 
          for facets $F$ of $Q$ such that $F \subset \widetilde H_{+}$,
    \item the combinatorial cubes 
               $\coneOP(\boldsymbol{p}_R,F)\cap H'_{-}\cap(\R^d\times\R_{+})$ 
          for facets $F$ of $Q$ such that $F\subset \widetilde H_{-}$,
    \item the combinatorial half-cubes 
              $\coneOP(\boldsymbol{p}_L,F)\cap\coneOP(\boldsymbol{p}_R,F)\cap H'_{-}$
          for facets $R$ of $Q$ that intersect~$\widetilde H$, and 
    \item $Q$.
  \end{compactenum}
   The facet defining hyperplanes of the upper facets of $T_{-}$ are
  \begin{compactenum}[\ \rm(a)]
     \item $\affOP(\boldsymbol{p}_L,F)$ for facets $F$ of $Q$ such that $F\subset \widetilde H_{+}$,
    \item $\affOP(\boldsymbol{p}_R,F)$ for facets $F$ of $Q$
    such that $F\subset \widetilde H_{-}$,
    \item $\affOP(\boldsymbol{p}_L,\boldsymbol{p}_R,F)$
      for facets $F$ of $Q$ that intersect  $\widetilde H$, and 
    \item $\affOP(Q)$.
  \end{compactenum}
\end{lemma}
\begin{proof}
  Since $T_{-}$ is the intersection of $T$ with $H'_{-}$,
  the upper facets of $T_{-}$ are given by $Q$ plus the intersections
  of the upper facets of $T$ with $H'_{-}$, and the
  vertices of $T_{-}$ are the vertices of $T$ and the vertices of~$Q$.
\end{proof}

\begin{lemma}[Combinatorial structure of \boldmath$R$]~\\  
   The set of vertices of $R$ consists of the vertices of $Q$ and
   all points in $V'':=\vertices{R}\setminus\vertices{Q}$.
   Furthermore, the set of (all)
   facets of $R$ consists of
  \begin{compactenum}[\ \rm(a)]
  \item the combinatorial cubes $\convOP(\boldsymbol{p}_R,F)\cap \widetilde H_{+}$ for facets $F$ of $Q$
    such that $F\subset \widetilde H_{+}$,
    \item the combinatorial cubes $\convOP(\boldsymbol{p}_L,F)\cap \widetilde H_{+}$ for facets $F$ of $Q$
    such that $F\subset \widetilde H_{-}$,
    \item the combinatorial half-cubes $\convOP(\boldsymbol{p}_R,F)\cap\convOP(\boldsymbol{p}_L,F)$ 
      for facets $F$ of $Q$ that intersect  $\widetilde H$, and
    \item $Q$.
  \end{compactenum}
  The set of facet defining hyperplanes of the
   facets of $R$ consists of
  \begin{compactenum}[\ \rm(a)]
  \item $\affOP(\boldsymbol{p}_R,F)$ for facets $F$ of $Q$
    such that $F\subset \widetilde H_{+}$,
    \item $\affOP(\boldsymbol{p}_L,F)$ for facets $F$ of $Q$
    such that $F\subset \widetilde H_{-}$,
    \item $\affOP(\boldsymbol{p}_L,\boldsymbol{p}_R,F)$
      for facets $F$ of $Q$ such that $F$ intersects  $\widetilde H$, and
    \item $\affOP(Q)$.

  \end{compactenum}
\end{lemma}



\begin{proof}[Proof of Theorem~\ref{thm:ex_of_sym_cubification}.]
  We show that the complex $\mathcal{C}'$ given by the upper facets of
  the polytope~$U$ of
  Construction~\ref{constr:gen_hexhoop} determines a lifted
  cubification $(\mathcal{B}^d,\phi)$ of $(\mathcal{S}^{d-1},\psi)$.

  First observe that no vertex of $T_{-}$ is beyond a facet of $R$,
  and no vertex of $R$ is beyond a facet of $T_{-}$.  Hence
  the boundary of $U=\convOP(T_{-}\cup R)$ is the union of the two
  boundaries of the two polytopes, excluding the relative
  interior of $Q$.  \\
  Define the vertex sets $V:=\vertices{\lifted{\mathcal{S}t}{\psi}}$, $V':=\vertices{Q}$ and  
  $V'':=\vertices{R}\setminus V'$. Then
  \begin{compactitem}[$\bullet$~]
      \item each vertex of $V$ is 
       beneath each facet of $R$ that is of type (a) or (b), and
      \item each vertex of $V''$ is 
       beneath each facet of $T_{-}$ that is of type (a) or (b). 
  \end{compactitem}
  Hence these four types of facets are facets of $U$ that are
  combinatorial cubes, and the set of vertices of $U$ is given by the
  union of $V,V'$ and $V''$. It remains to show that each hyperplane
  $\affOP(\boldsymbol{p}_L,\boldsymbol{p}_R,F)$, where $F$ is a facet
  of $Q$ that intersects $\widetilde H$, is the affine hull of a
  cubical facet of~$U$.  To see this, observe that there are two
  facets $F_+$, $F_-$ of $R$, $T_-$ respectively, that are both
  contained in the affine hull of $F$. These two facets $F_+$, $F_-$
  are both half-cubes that intersect in a common $(d-1)$-cube, namely
  $F$. Furthermore, all vertices of $F_+$ and of $F_-$ that are not
  contained in $\affOP(F)$ are contained in $\widetilde{H}$. Hence the
  union of $F_+$ and $F_-$ is a combinatorial cube.
  
  Thus every upper facet of $U$ is a combinatorial cube.
  Furthermore, $\pi(R)=\pi(Q)$ and $\pi(T_-)=\support{P}$, so
  the upper facets of $U$ determine a lifted cubical
  subdivision of~$(\mathcal{S}^{d-1},\psi)$.
\end{proof}
                          
%
\begin{proposition}[Dual manifolds]\label{prop:dual_mf_of_hescap}
  Up to PL-homeomorphism, the generalized regular Hexhoop
  cubification $\mathcal{B}^d$ of $\mathcal{S}^{d-1}$
  has the following dual manifolds:
  \begin{compactitem}
  \item  $\mathcal{N}\times I$ for each dual manifold $\mathcal N$ 
   (with or without boundary) 
   of~$\mathcal{S}_L=\mathcal{S}^{d-1}\cap \widetilde H_+$,
  \item two $(d-1)$-spheres ``around'' $Q$, $Q^M$, respectively,
    where the upper index $^M$ denotes the mirrored copy.
  \end{compactitem}
\end{proposition}

\begin{proof}
  The ``main part'' of the complex $\mathcal{B}^d$ may be 
  viewed as a prism of height~$4$, whose dual manifolds are
  of the form $\mathcal{N}\times I$, as well as
  four $(d-1)$-balls.
  This prism is then modified by glueing a full torus
  (product of the $(d-2)$-sphere $\mathcal{S}^{d-1}\cap\mathcal{H}$
  with a square $I^2$) into its ``waist.''
  This extends the dual manifolds $\mathcal{N}\times I$ without
  changing the PL-homeomorphism type, while closing the four
  $(d-1)$-balls into two intersecting, embedded spheres.

  We refer to Figure~\ref{fig:dms_of_gen_hexhoop}
  (case $d=2$) and Figure~\ref{fig:dmf_three_dim_gen_Hexhoop} ($d=3$)
  for geometric intuition.

\begin{Figure} \vskip-5mm
  \psfrag{p_1}{}
  \psfrag{p_2}{}
  \psfrag{P}{}
  \psfrag{H2}{}
  \psfrag{not a cube}{}
  \subfigure[{Dual product manifolds $\mathcal{N}\times I$}]
   {\includegraphics[width=.35\textwidth]{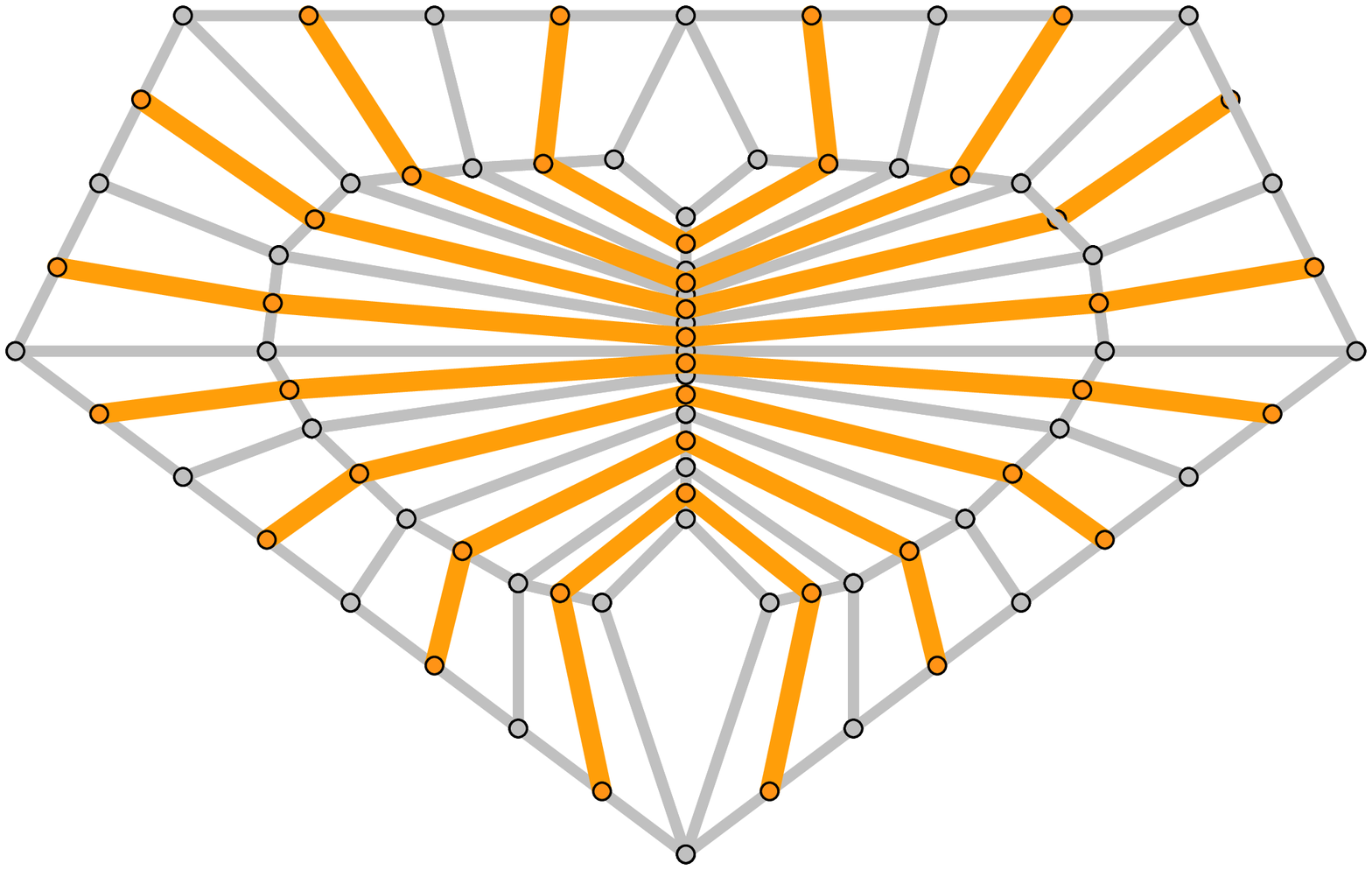}}\qquad\qquad
  \subfigure[Dual spheres]{\includegraphics[width=.35\textwidth]{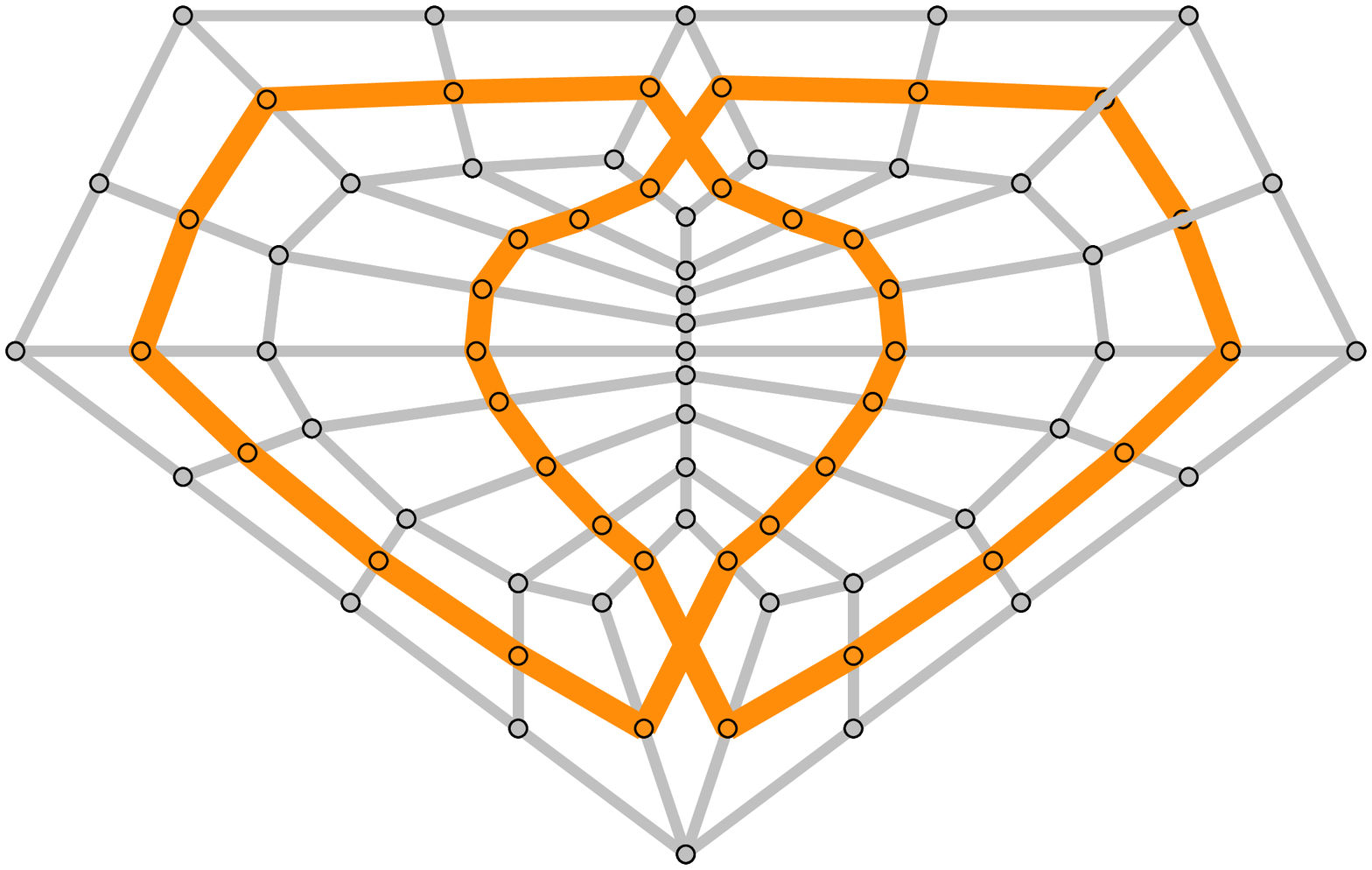}}
  \vskip-5mm
      \caption[Dual manifolds of a $2$-dim.~generalized regular
      Hexhoop]{%
         The dual manifolds of a $2$-dimensional 
        generalized regular Hexhoop. }
      \label{fig:dms_of_gen_hexhoop}
\end{Figure}

\vskip-5mm
\begin{Figure} 
  \includegraphics[width=.45\textwidth]{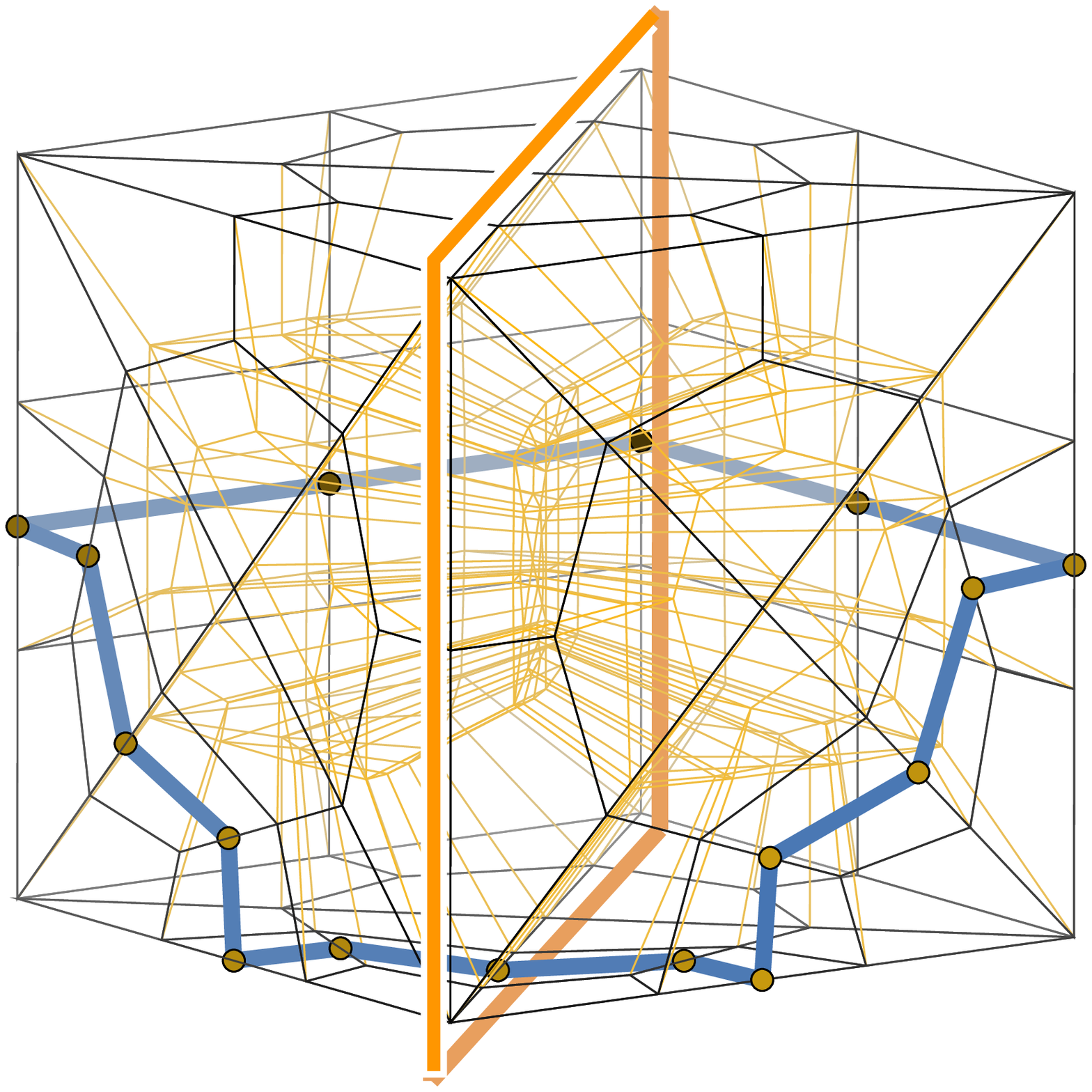}
  \rput(-27mm,60mm){$H$}
  \rput(-25mm,5mm){$\mathcal{N}$}
  \rput(-70mm,55mm){$\mathcal{B}^d$}
  \transformsToArrow{30mm}{10mm}
  \includegraphics[width=.45\textwidth]{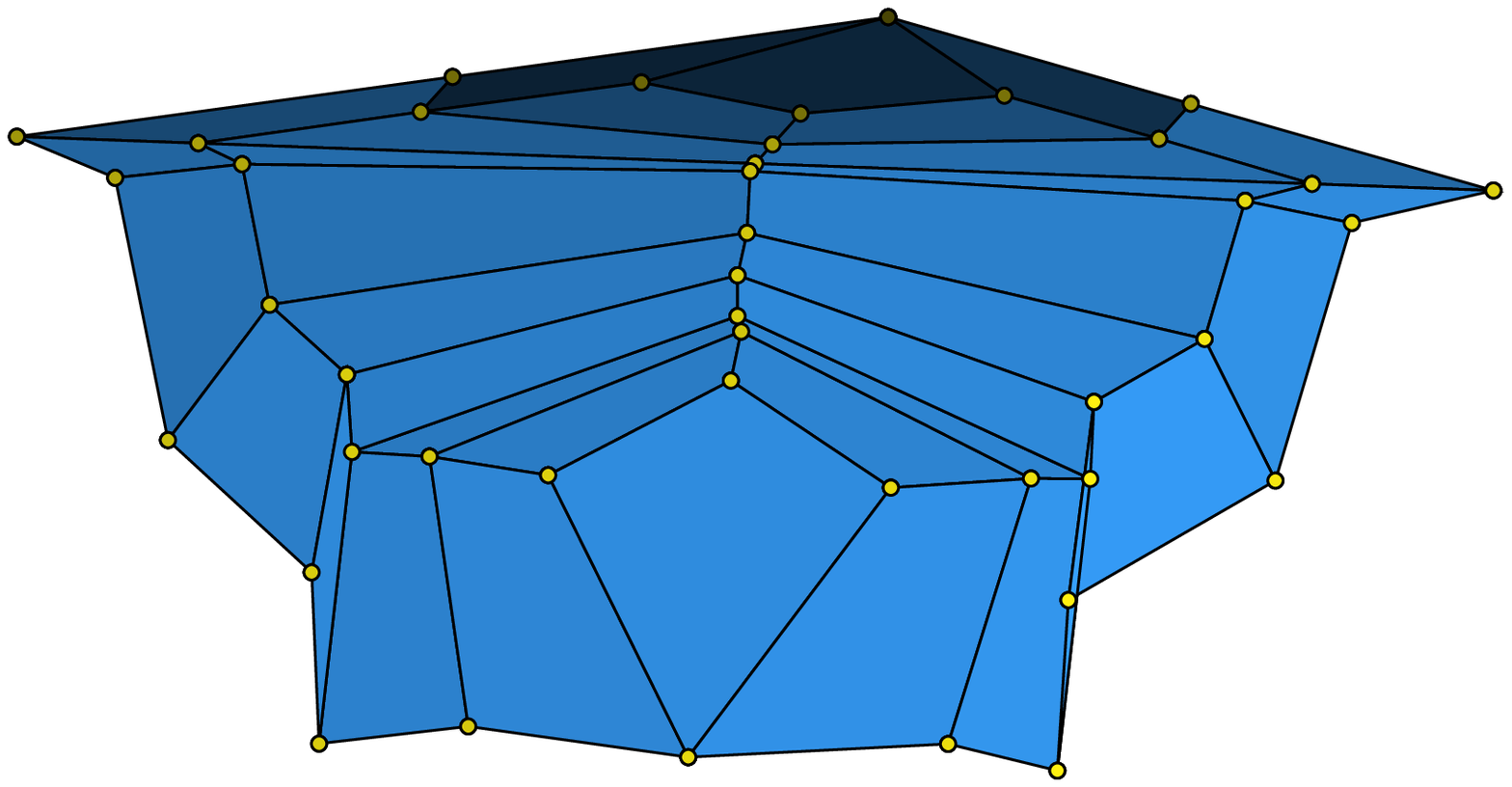}
  \rput(-27mm,60mm){}
  \rput(-25mm,5mm){$\mathcal{N}'$}
  \rput(-70mm,55mm){}
\vskip-5mm
      \caption[Dual manifolds of a $3$-dim.~generalized regular Hexhoop]{%
        A $3$-dimensional cubification produced by the generalized
        regular Hexhoop construction. For every embedded dual circle
        $\mathcal{N}$ which intersects $H_{+}\setminus H$ and
        $H_{-}\setminus H$, there is an embedded dual $2$-ball
        $\mathcal{N}'$ with boundary~$\mathcal{N}$ in the the
        cubification. (This is a cubification for the case ``single5'' 
        introduced in Section~\ref{sec:immersions}.)}  
      \label{fig:dmf_three_dim_gen_Hexhoop}
\end{Figure}
\end{proof}


\newcommand{\InvariantConsistency}[1]{$(\text{I}_{#1}1)$}
\newcommand{\InvariantPLEquivalence}[1]{$(\text{I}_{#1}2)$}
\newcommand{\InvariantSymmetry}[1]{$(\text{I}_{#1}3)$}
\newcommand{\InvariantDiagSubcomplex}[1]{$(\text{I}_{#1}4)$}

\section{Cubical 4-polytopes with prescribed dual manifold immersions}
\label{sec:immersions}
 
Now we use our arsenal of cubical construction techniques
for the construction of cubifications 
with prescribed dual $2$-manifold immersions, and 
thus approach our main theorem.

For this we ask for our input to be given by normal crossing PL-immersions
whose local geometric structure is rather special: 
We assume that
$\mathcal{M}^{d-1}$ is a ($d-1$)-dimensional cubical
PL-manifold, and $j:\mathcal{M}^{d-1}\immersed\R^{d}$ is a
\emph{grid immersion},\index{grid immersion}\index{immersion!grid immersion} 
a cubical normal crossing codimension one
immersion into $\R^{d}$ equipped with the standard unit cube
structure.

\subsection{From PL immersions to grid immersions}

In view of triangulation and approximation methods available in PL and
differential topology, the above assumptions are not so
restrictive. (See, however, Dolbilin et al.~\cite{Dolbilin1995}
for extra problems and obstructions that may arise without
the PL assumption, and if we do not admit subdivisions,
even for the high codimension embeddings/immersions.)

\begin{proposition}\label{prop:PL_to_grid_imm}
  Every locally flat normal crossing immersion of a compact
  $(d-1)$-manifold into~$\R^{d}$ is PL-equivalent to a 
  grid immersion of a cubification of the manifold into the
  standard cube subdivision of~$\R^{d}$.
\end{proposition}

\begin{proof}
  may replace any PL-immersion of~$\mathcal{M}^{d-1}$ by a simplicial
  immersion into a suitable triangulation of~$\R^{d}$. The vertices
  of~$j(\mathcal{M}^{d-1})$ may be perturbed into general position.
  
  Now we overlay the polyhedron $j(\mathcal{M}^{d-1})$ with a cube
  structure of~$\R^{d}$ of edge length $\varepsilon$ for suitably
  small $\varepsilon>0$, such that the vertices
  of~$j(\mathcal{M}^{d-1})$ are contained in the interiors of distinct
  $d$-cubes.
  
  Then working by induction on the skeleton, within each face of the
  cube structure, the restriction of~$j(\mathcal{M}^{d-1})$ to a
  $k$-face --- which by local flatness consists of one or several
  $(k-1)$-cells that intersect transversally --- is replaced by a
  standard cubical lattice version that is supposed to run through the
  interior of the respective cell, staying away distance
  $\varepsilon'$ from the boundary of the cell; here we take different
  values for $\varepsilon'$ in the situation where the immersion is
  not embedded at the vertex in question, that is, comes from several
  disjoint neighborhoods in~$\mathcal{M}^{d-1}$.
\begin{Figure}
   \includegraphics[height=60mm,width=.8\textwidth]{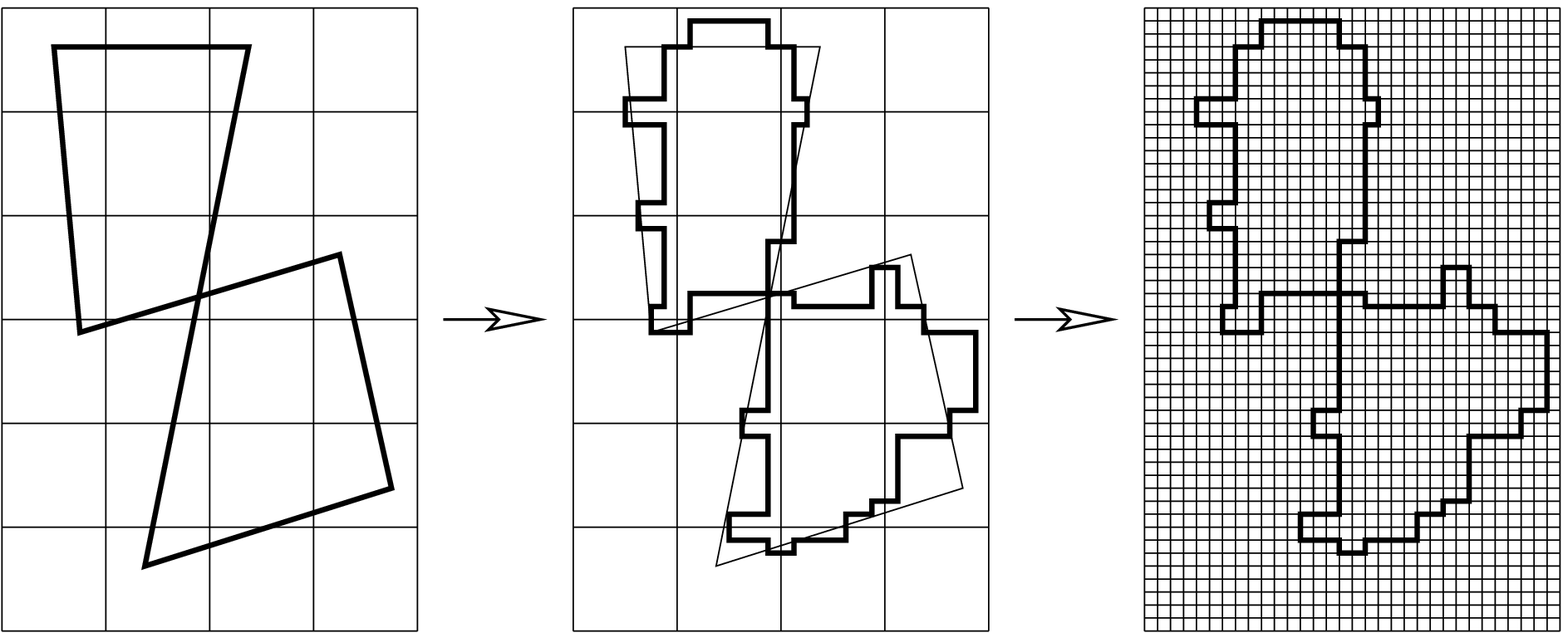}
   \caption{Illustration of the proof of Proposition~\ref{prop:PL_to_grid_imm}.}
   \label{fig:pl_to_grid}
\end{Figure}
The resulting modified immersion into $\R^{d}$ will be cellular with
respect to a standard cube subdivision of edge length
$\frac1N\varepsilon$ for a suitable large~$N$.
Figure~\ref{fig:pl_to_grid} illustrates this for~$d=2$.
\end{proof}

\subsection{Vertex stars of grid immersions of surfaces}

From now on, we restrict our attention to the case of $d=3$,
that is, $2$-manifolds and $4$-polytopes.

 There are nine types of vertex stars of grid immersions of surfaces,
 namely the following five vertex stars of a regular vertex,
     \begin{Figure}
        \subfig{single3}{\includegraphics[width=.15\textwidth]{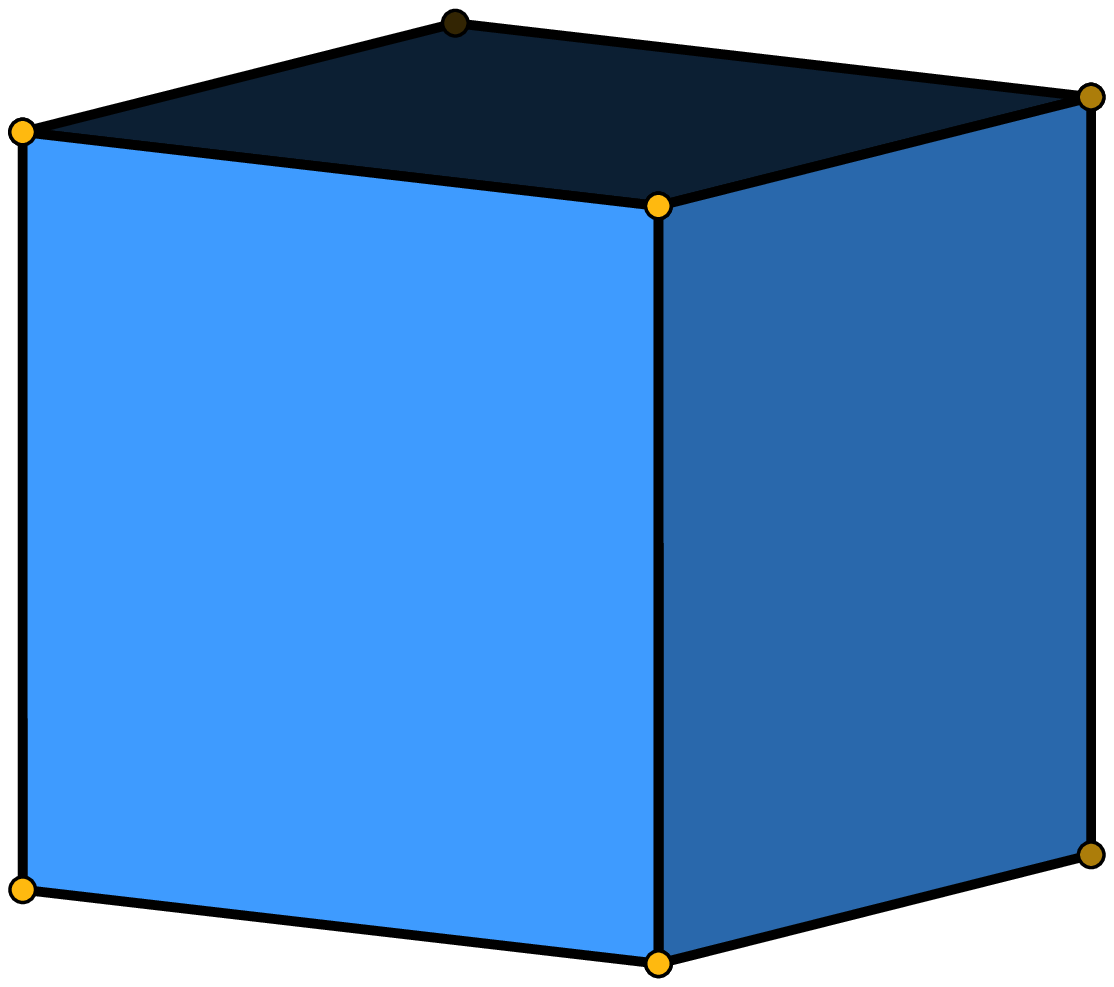}}
        \subfig{single4a}{\includegraphics[width=.15\textwidth]{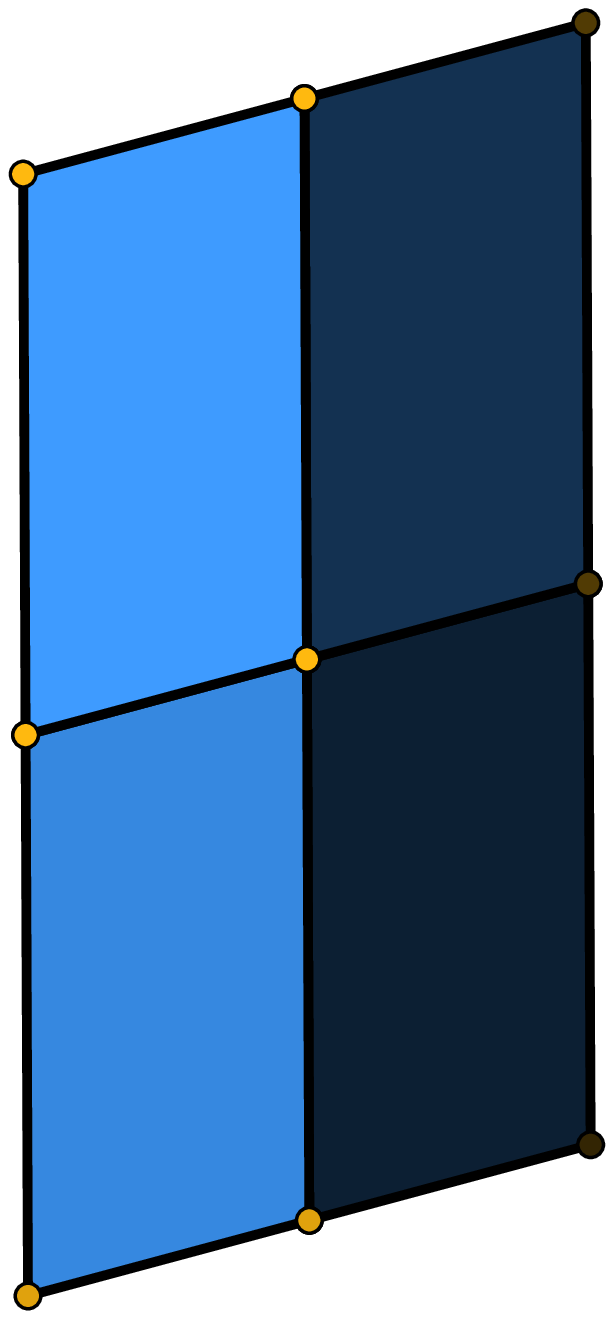}}
        \subfig{single4b}{\includegraphics[width=.15\textwidth]{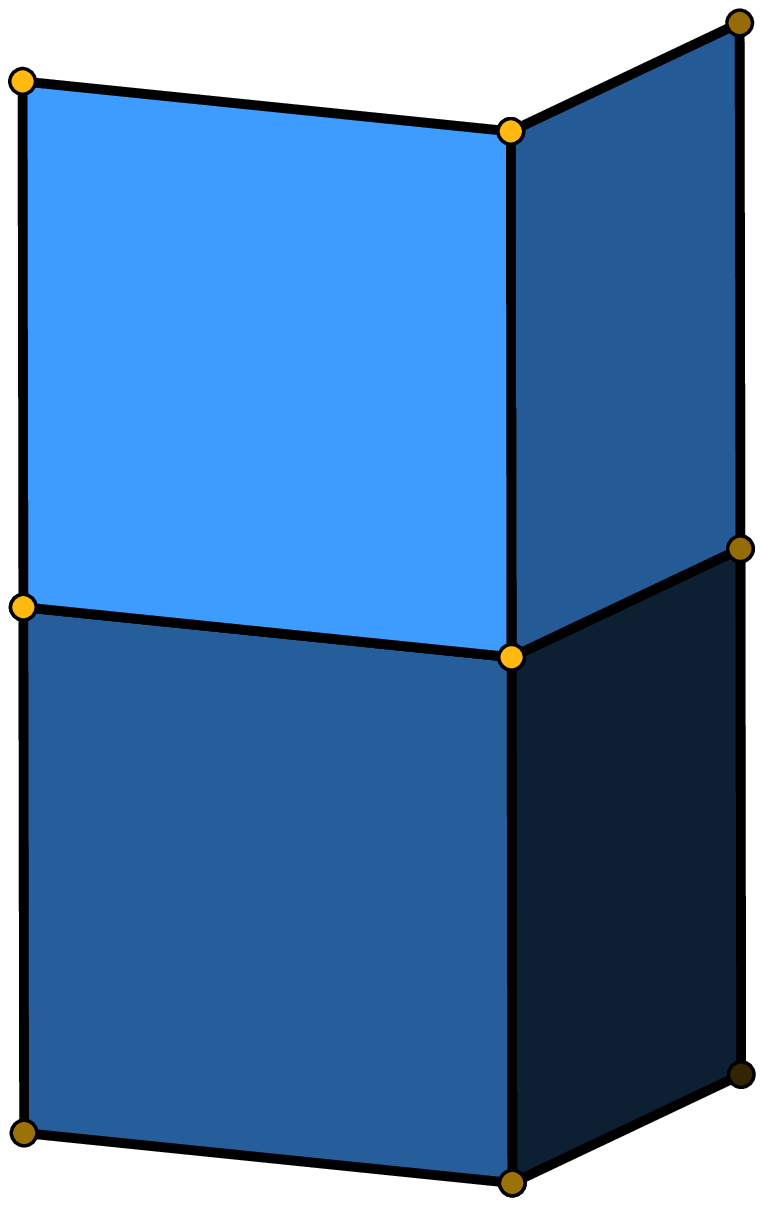}}
        \subfig{single5}{\includegraphics[width=.15\textwidth]{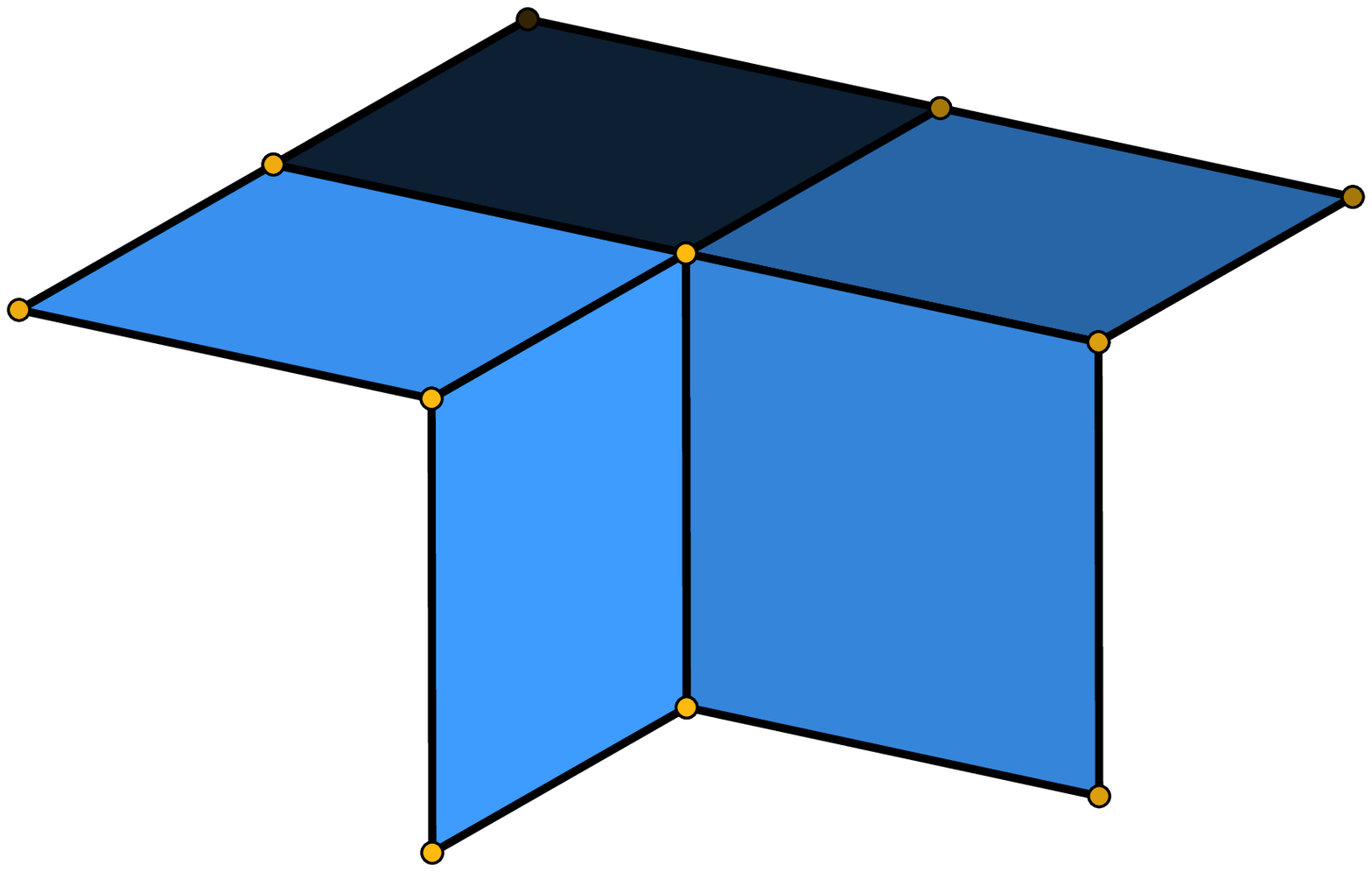}}
        \subfig{single6a}{\includegraphics[width=.15\textwidth]{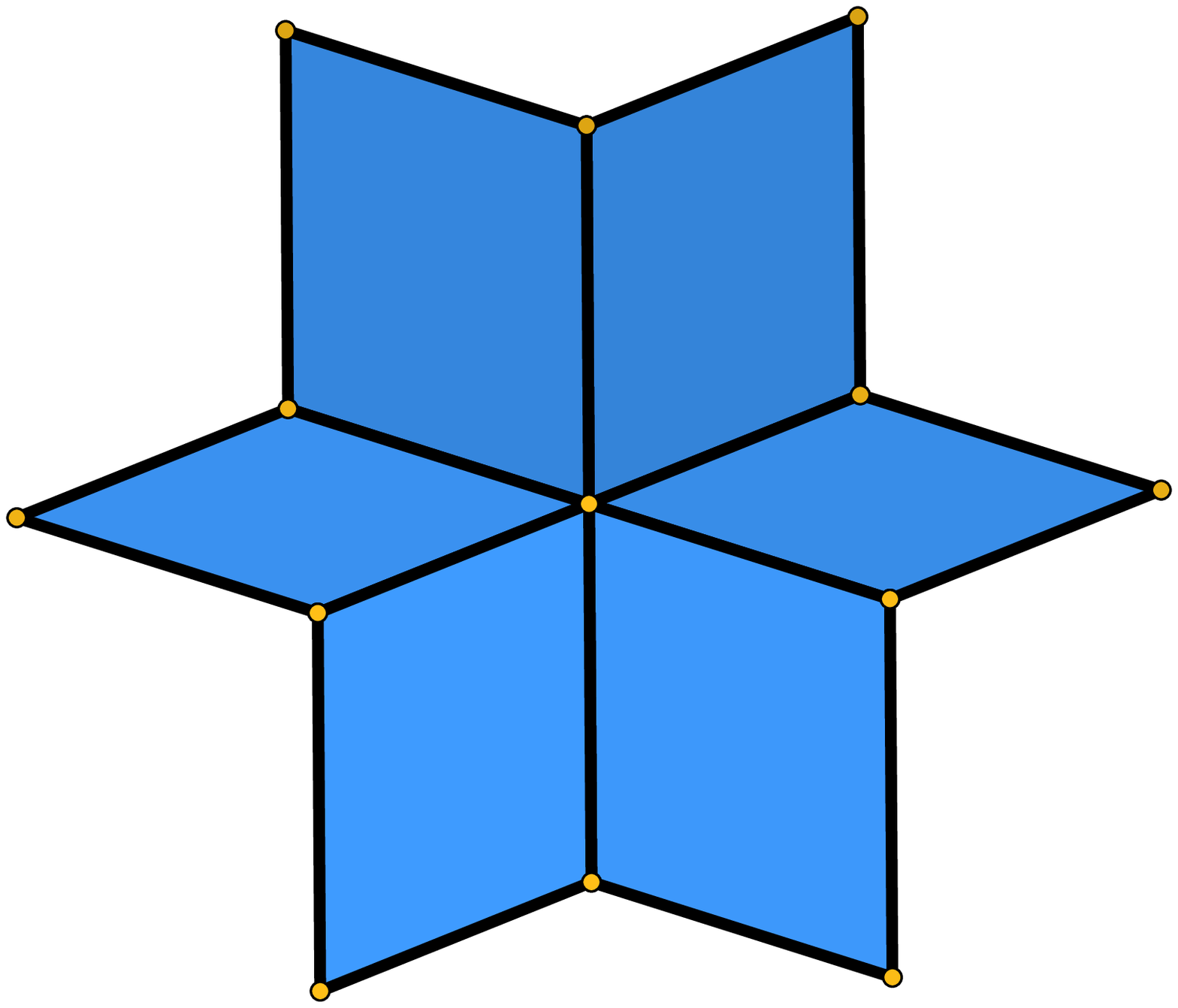}}
        \subfig{single6b}{\includegraphics[width=.15\textwidth]{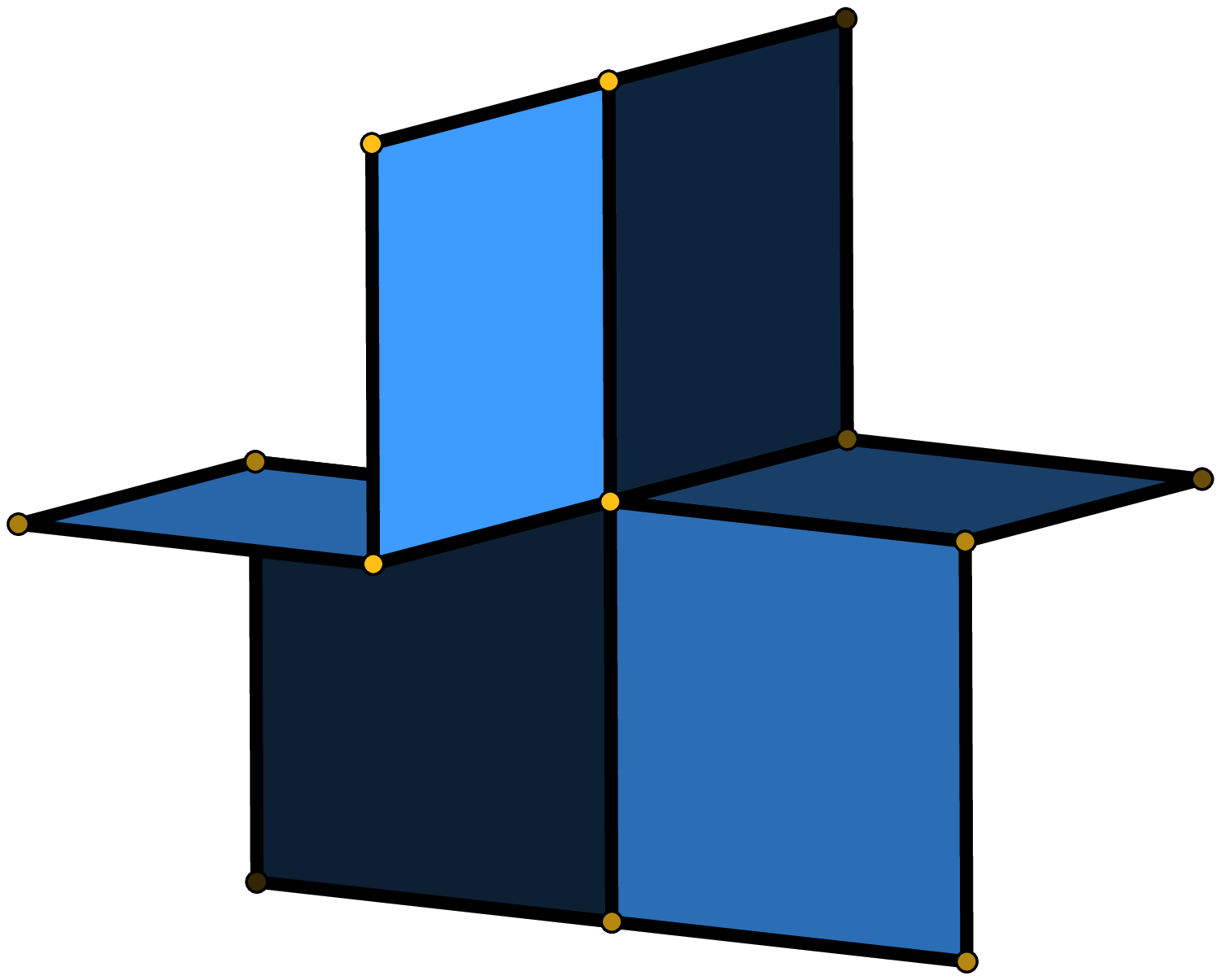}}
     \end{Figure}
  plus two vertex stars with double intersection and the vertex star of a triple
  intersection point:
     \begin{Figure}
        \subfig{double8a}{\includegraphics[width=.15\textwidth]{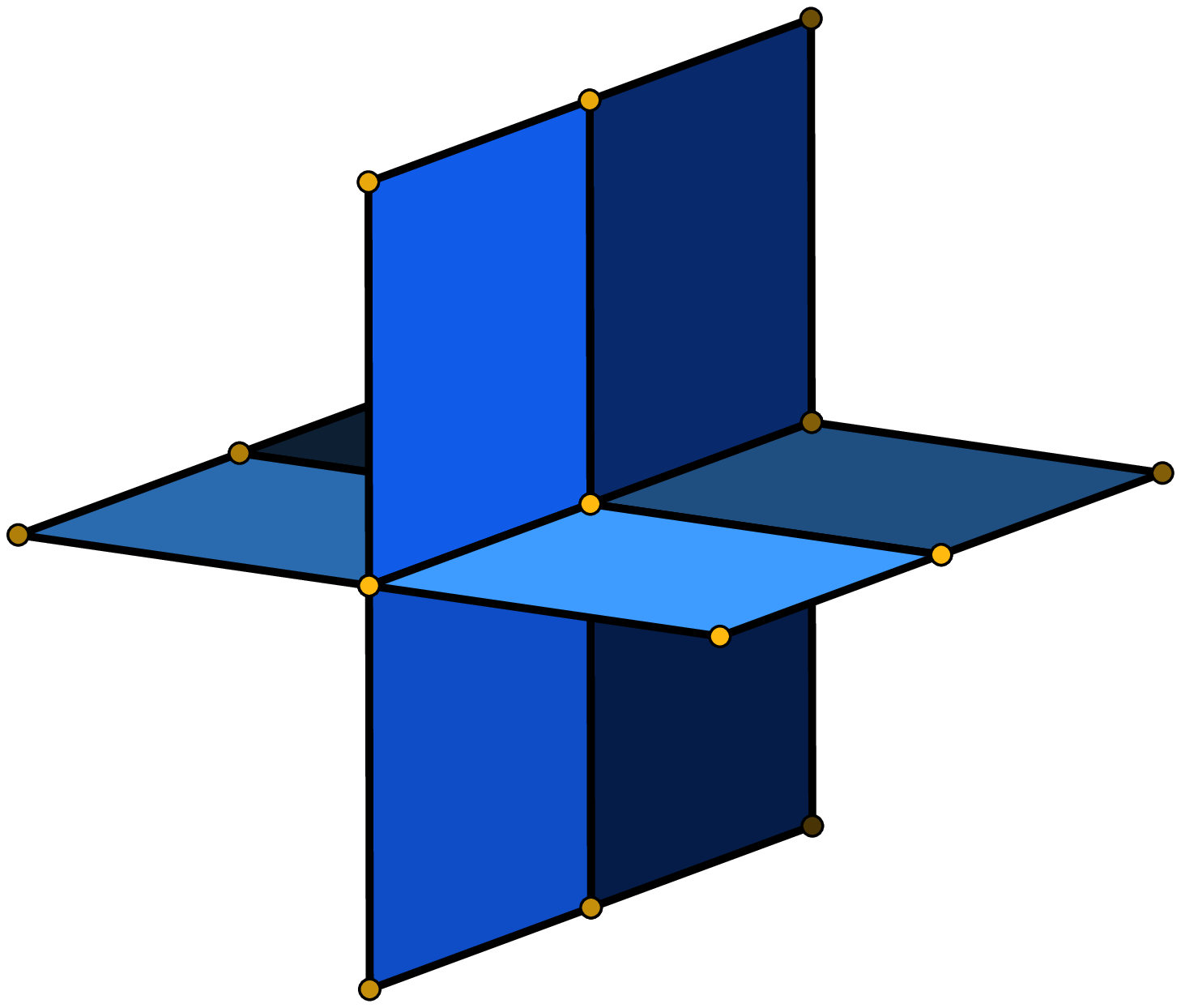}}
        \subfig{double8b}{\includegraphics[width=.15\textwidth]{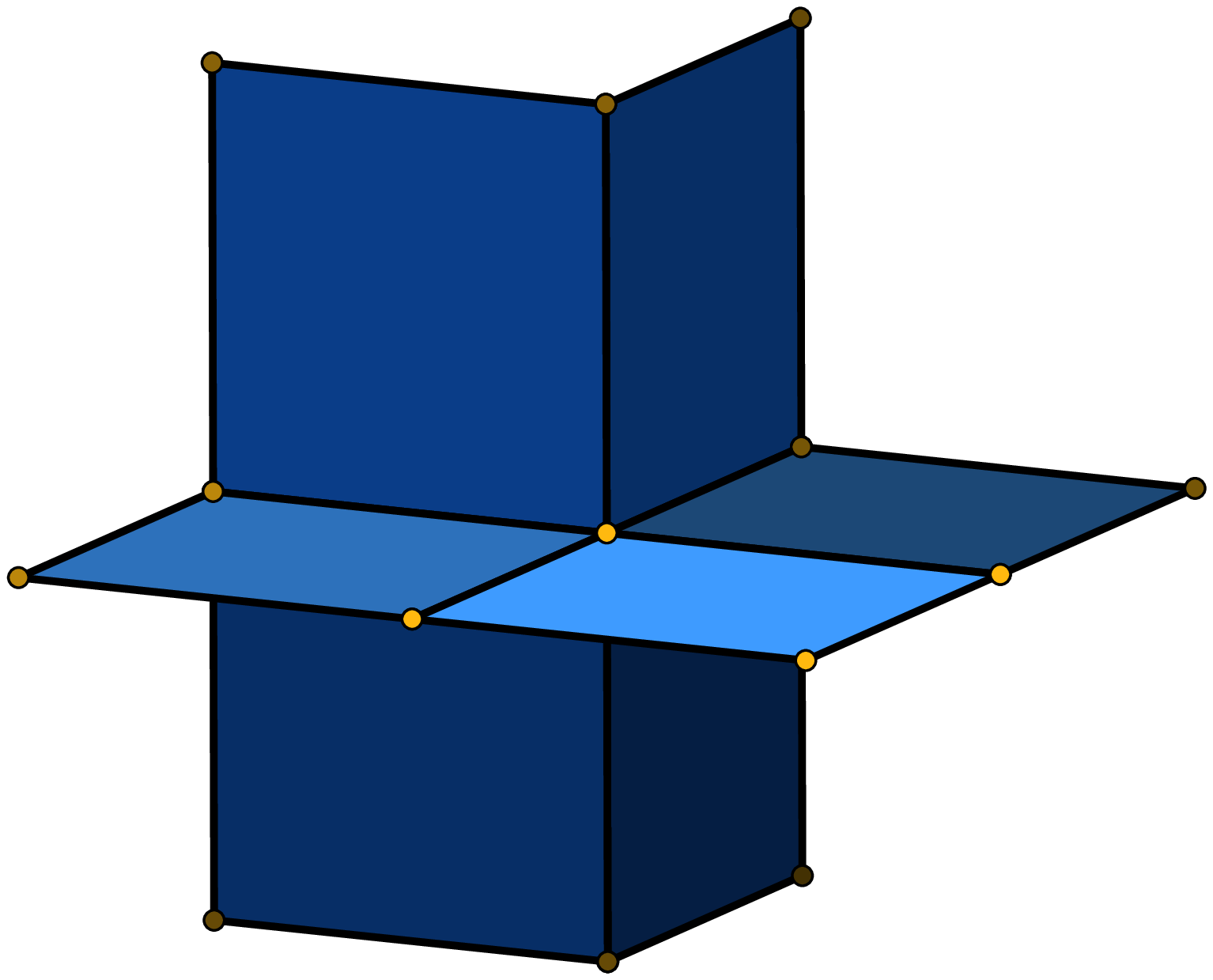}}
        \subfig{triple12}{\includegraphics[width=.15\textwidth]{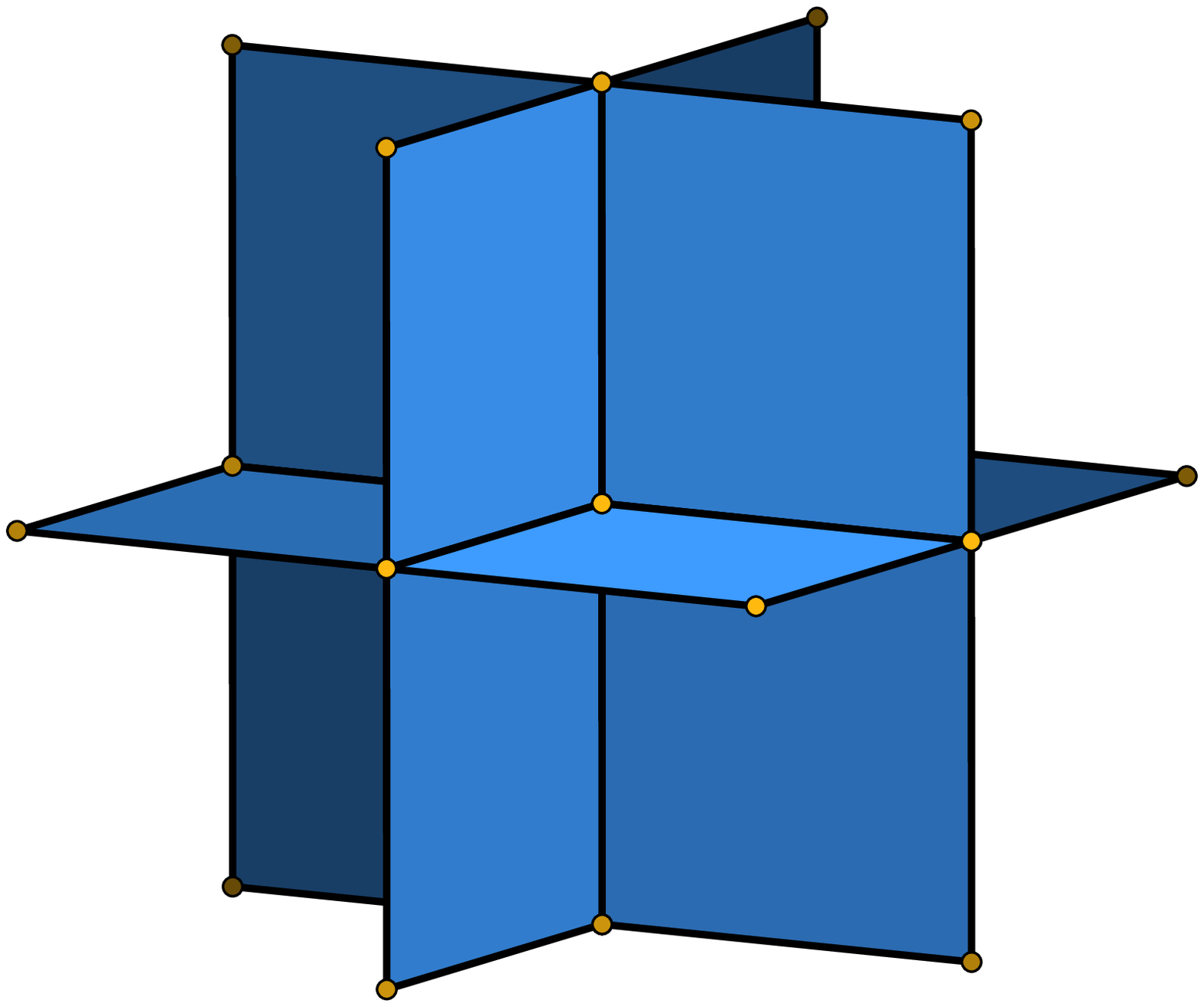}}
      \end{Figure}

For the constructions below we will require that the grid immersion 
$j:\mathcal{M}^2\immersed\R^3$ that we start out with is 
\emph{locally symmetric},\index{immersion!locally symmetric} that is, that at
     each vertex~$\boldsymbol{w}$ of $j(\mathcal{M}^2)$ there is
     plane $H$ through $\boldsymbol{w}$ such that for each
     vertex~$\boldsymbol{v}$ with $j(\boldsymbol{v})=\boldsymbol{w}$
     the image of the vertex star of $\boldsymbol{v}$ is symmetric
     with respect to~$H$.  Thus we require that $H$ is a symmetry
     plane separately for each of the (up to three) local sheets
     that intersect at~$\boldsymbol{w}$.
     Such a plane $H$ is necessarily of the form $x_i=k$, $x_i+x_j=k$,
     or $x_i- x_j=k$.  In the first case we say $H$ is a
     \emph{coordinate hyperplane}, and in other cases it is
     \emph{diagonal}.

\begin{proposition}\label{prop:local_symmetry}
  Any grid immersion of a compact cubical $2$-manifold into~$\R^3$ is
  equivalent to a locally symmetric immersion of the same type.
\end{proposition}

\begin{proof}
  All the vertex stars displayed above satisfy the 
  local symmetry condition, with a single exception, namely the star 
  ``single6b'' of a regular vertex with six adjacent quadrangles.
  As indicated in
  Figure~\ref{fig:get_rid_of_s6a}, a local modification of the surface
  solves the problem (with a suitable refinement of the standard cube
  subdivision).
\end{proof}
     \begin{Figure}
     \raisebox{2cm}{\subfig{single6b}{\includegraphics[height=.25\textwidth]{img/grid/cases/single_6b}}}
          \transformsToArrow{.075\textheight}{8mm}\qquad
    \includegraphics[height=.25\textwidth]{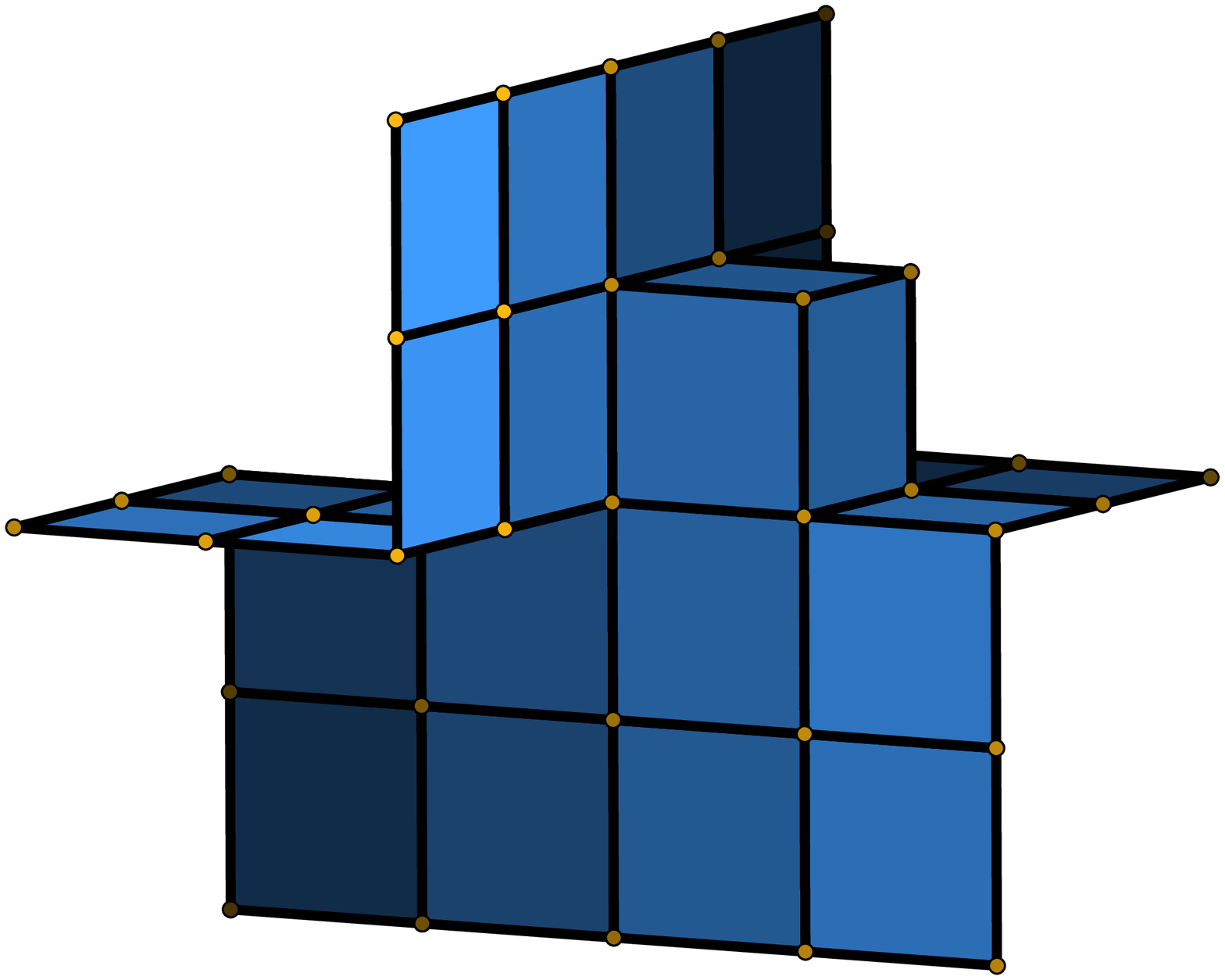} 
        \caption{Local modification used to ``repair'' the case ``single6a.''}
          \label{fig:get_rid_of_s6a}
     \end{Figure}

\subsection{Main theorem (2-manifolds into cubical 4-polytopes)}

\begin{theorem}\label{thm:C4P_with_prescribed_dmf} \label{thm:gridimm_to_C4P}
  Let $j\colon\mathcal{M}\immersed\mathbbm{R}^3$ be a locally flat
  normal crossing immersion of a compact $2$-manifold (without boundary)
  $\mathcal{M}$ into $\R^3$.
 
   Then there is a cubical $4$-polytope $P$ with a 
   dual manifold $\mathcal{M}'$ and associated immersion
   $y:\mathcal{M}'\immersed\support{\boundaryOP P}$ 
   such that the following conditions are satisfied:
   \begin{thm_enumerate_i}
   \item \label{main_theorem:PL_euqiv} $\mathcal{M}'$ is a cubical
     subdivision of $\mathcal{M}$, and the immersions $j$ (interpreted
     as a map to $\R^3\cup\{\infty\}\cong S^3$) and $y$ are
     PL-equivalent.
     \item  \label{main_theorem:parity_of_f_3}
       The number of facets of~$P$ is congruent modulo~$2$
       to the number~$t(j)$ of triple points of the immersion~$j$.%
       \label{eq:parity-of-d-dim-intersection-points}
       
     \item \label{main_theorem:odd}
       If the given surface $\mathcal{M}$ is non-orientable and of
       odd genus, then the cubical $4$-polytope~$P$ has an odd number
       of facets.
    \end{thm_enumerate_i}
\end{theorem}


The core of our proof is the following construction
of cubical $3$-balls with a prescribed dual manifold immersions.
\medskip

\begin{construction}{Regular cubical $3$-ball with a prescribed dual manifold}
\label{constr:cubical-3ball_with_prescribed_dmf}
\begin{inputoutput}
\item[Input:] A $2$-dimensional closed (that is, compact and without boundary) 
  cubical PL-surface $\mathcal{M}$, and a {locally symmetric}
  codimension one grid immersion
  $j\colon\mathcal{M}\immersed\support{\PileOfCubes{3}{\ell_1,\ell_2,\ell_3}}\subset\mathbbm{R}^3$.
\item[Output:] A regular convex $3$-ball $\mathcal{B}$ with a dual
  manifold $\mathcal{M}'$ and associated immersion
  $y:\mathcal{M}'\immersed\support{\mathcal{B}}$ such that the
  following conditions are satisfied:
   \begin{thm_enumerate_i}
     \item  $\mathcal{M}'$ is a cubical subdivision of $\mathcal{M}$,
       and the immersions $j$ and $y$ are PL-equivalent.
     \item The number of facets of~$\mathcal{B}$ is congruent modulo
       two to the number~$t(j)$ of triple points of the immersion
       $j$.
    \end{thm_enumerate_i}
\end{inputoutput}

\begin{enumerate}[\bf(1)]
\item \textbf{Raw complex.}
      Let $\mathcal{A}$ be a copy of the pile of
      cubes $\PileOfCubes3{\ell_1+1,\ell_2+1,\ell_3+1}$ with all
      vertex coordinates shifted by $-\tfrac{1}{2}\mathbbm{1}$.
      (Hence
$x_i\in\{-\tfrac{1}{2},\tfrac{1}{2},\tfrac{3}{2},\dots,\ell_i+\tfrac{1}{2}\}$
      for each vertex $\boldsymbol{x}\in\vertices{A}$.)

      Due to the local symmetry of the immersion, and the choice of the vertex
      coordinates of~$\mathcal{A}$, the following holds:
      \begin{dense_itemize}
         \item Each vertex of~$j(\mathcal{M})$ is the barycenter of a
           $3$-cube~$C$ of~$\mathcal{A}$.
         
         \item For each $3$-cube $C$ of~$\mathcal{A}$ the restriction
           $(C,j(\mathcal{M})\cap C)$ is locally symmetric.\\

      \end{dense_itemize}

\item \textbf{Local subdivisions.}
We construct the lifted cubical subdivision $\mathcal{B}$ of~$\mathcal{A}$
by induction over the skeleton:
For $k=1,2,3$, $\mathcal{B}^k$ will be a lifted cubical subdivision
of the $k$-skeleton $\faces{k}{\mathcal{A}}$,
with the final result $\mathcal{B}:=\mathcal{B}^3$.  
For each $k$-face $F\in\mathcal{A}$ we take for the restriction
$\mathcal{B}^k\cap F$ a congruent copy 
from a finite list of \emph{templates}. 

Consider the following invariants (for $k\in\{1,2,3\}$).
\begin{compactsymboldesc}[\InvariantConsistency{kk}]
  \item[\InvariantConsistency{k}] \emph{Consistency requirement}.\\ For every
    $k$-face $Q\in\faces{k}{\mathcal{A}}$ and every facet $F$ of $Q$,
    the induced subdivision $\mathcal{B}^k \cap F$ equals
    $\mathcal{B}^{k-1} \cap F$.
  \item[\InvariantPLEquivalence{k}] \emph{PL equivalence requirement}.\\
    For every $k$-face $Q\in\faces{k}{\mathcal{A}}$ and every dual
    manifold $\mathcal{N}$ of $Q$ (with boundary) the cubical
    subdivision $\mathcal{B}^k \cap Q$ has a dual manifold that is
    PL-equivalent to~$j(\mathcal{N})\cap Q$.
  \item[\InvariantSymmetry{k}] \emph{Symmetry requirement}.\\
    Every symmetry of $(Q,j(\mathcal{M})\cap Q)$ for a $k$-face
    $Q\in\faces{k}{\mathcal{A}}$ that is a symmetry of each sheet
    of~$j(\mathcal{M})\cap Q$ separately is a symmetry of $(Q,
    \mathcal{B}^k \cap Q)$.
  \item[\InvariantDiagSubcomplex{k}] \emph{Subcomplex  requirement}.\\
    For every diagonal symmetry hyperplane $H_Q$ of a facet $Q$ of $\mathcal{A}$ 
    and every facet $F$ of $Q$
    the (lifted) induced subdivision $\mathcal{B}^k \cap (F\cap H)$ 
    is a (lifted) subcomplex of~$\mathcal{B}^k$.

\end{compactsymboldesc}
These invariants are maintained while iteratively constructing
$\mathcal{B}^1$ and $\mathcal{B}^2$. The resulting lifted cubical
subdivision $\mathcal{B}^3$ of $\mathcal{A}$ will satisfy
\InvariantConsistency{3} and \InvariantPLEquivalence{3}, but not
in general the other two conditions.

\item \textbf{Subdivision of edges.}
      Let $e$ be an edge of~$\mathcal{A}$.
      \begin{compactitem}[$\bullet$~]
      \item If $e$ is not intersected by the immersed manifold, then
        we subdivide the edge by an affine copy $\mathcal{B}^1_e$ of
        the following lifted subdivision $\mathcal{U}_2:=(\mathcal{U}'_2,h)$
        of $\PileOfCubes{1}{2}$: 
      \begin{Figure}
         \psfrag{0}{$0$}
         \psfrag{1/4}{\rput(0mm,0mm){$\frac{1}{4}$}}
         \includegraphics[width=.6\textwidth]{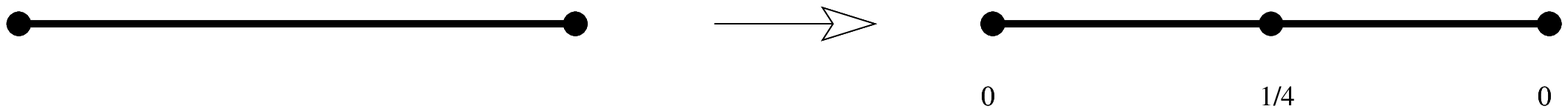}
     \end{Figure}    
      
   \item If $e$ is intersected by the immersed manifold, then we
     subdivide the edge by an affine copy $\mathcal{B}^1_e$ of the
     following lifted subdivision $\mathcal{U}_3:=(\mathcal{U}'_3,h)$
     of~$\PileOfCubes{1}{3}$:
      \begin{Figure}
         \psfrag{0}{$0$}
         \psfrag{1}{\rput(.2mm,0mm){$\frac{1}{6}$}}
         \includegraphics[width=.6\textwidth]{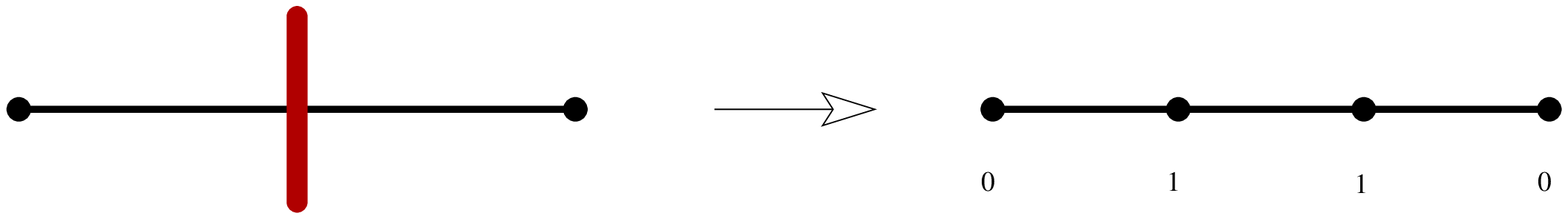}
      \end{Figure}    
    \end{compactitem}
    Observe that
     \InvariantConsistency{1}--\InvariantDiagSubcomplex{1} are satisfied.
 \bigskip

    \item \textbf{Subdivision of 2-faces.}
      Let $Q$ be a quadrangle of~$\mathcal{A}$,
      and $\boldsymbol{w}$ the unique vertex of~$j(\mathcal{M})$ that
      is contained in $Q$.  There are four possible types of restrictions
      of the grid immersion to~$Q$:

       \begin{Figure}
       \subfig{empty}{\includegraphics[width=.1\textwidth]{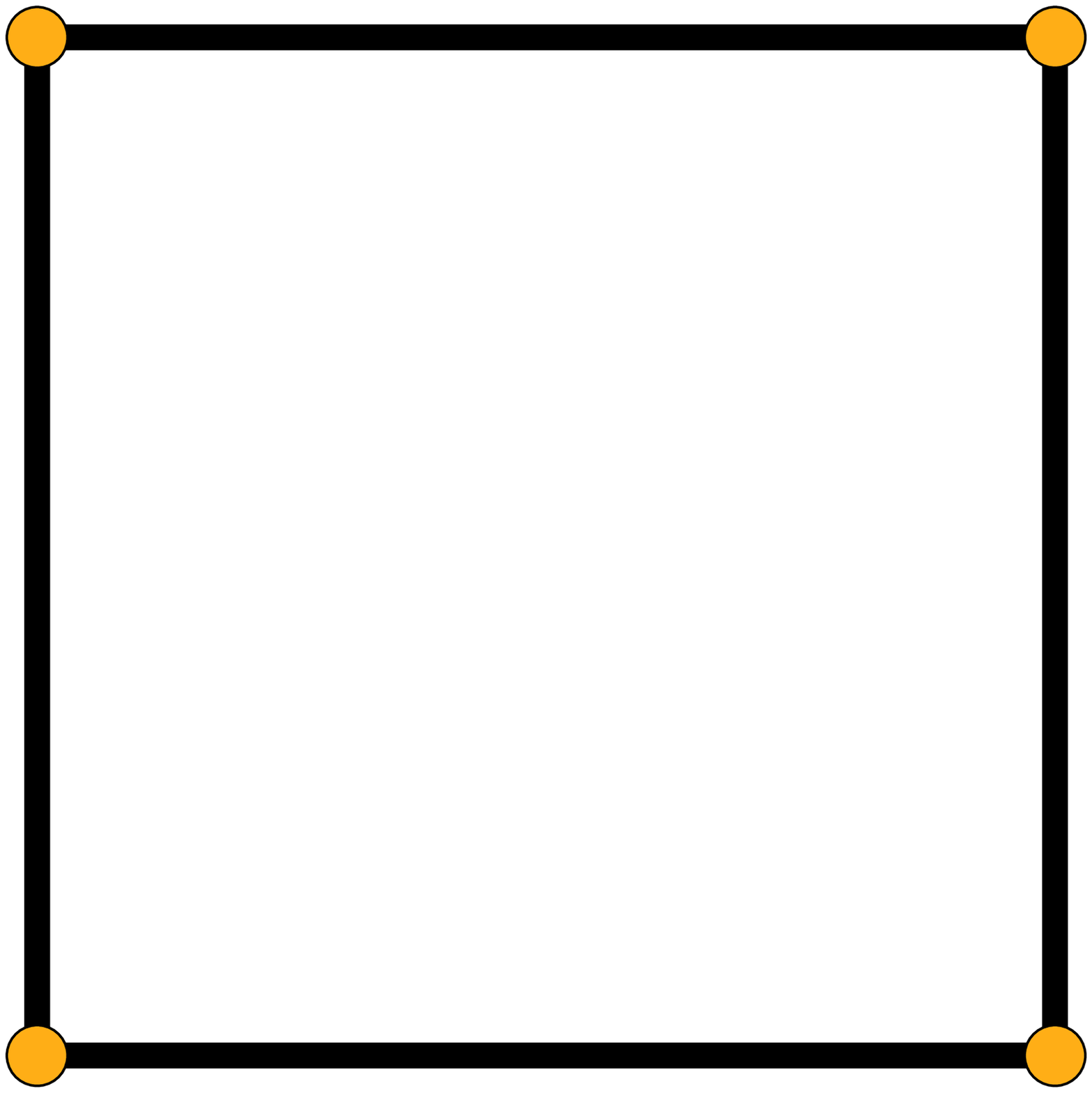}}\qquad
       \subfig{single2a}{\includegraphics[width=.1\textwidth]{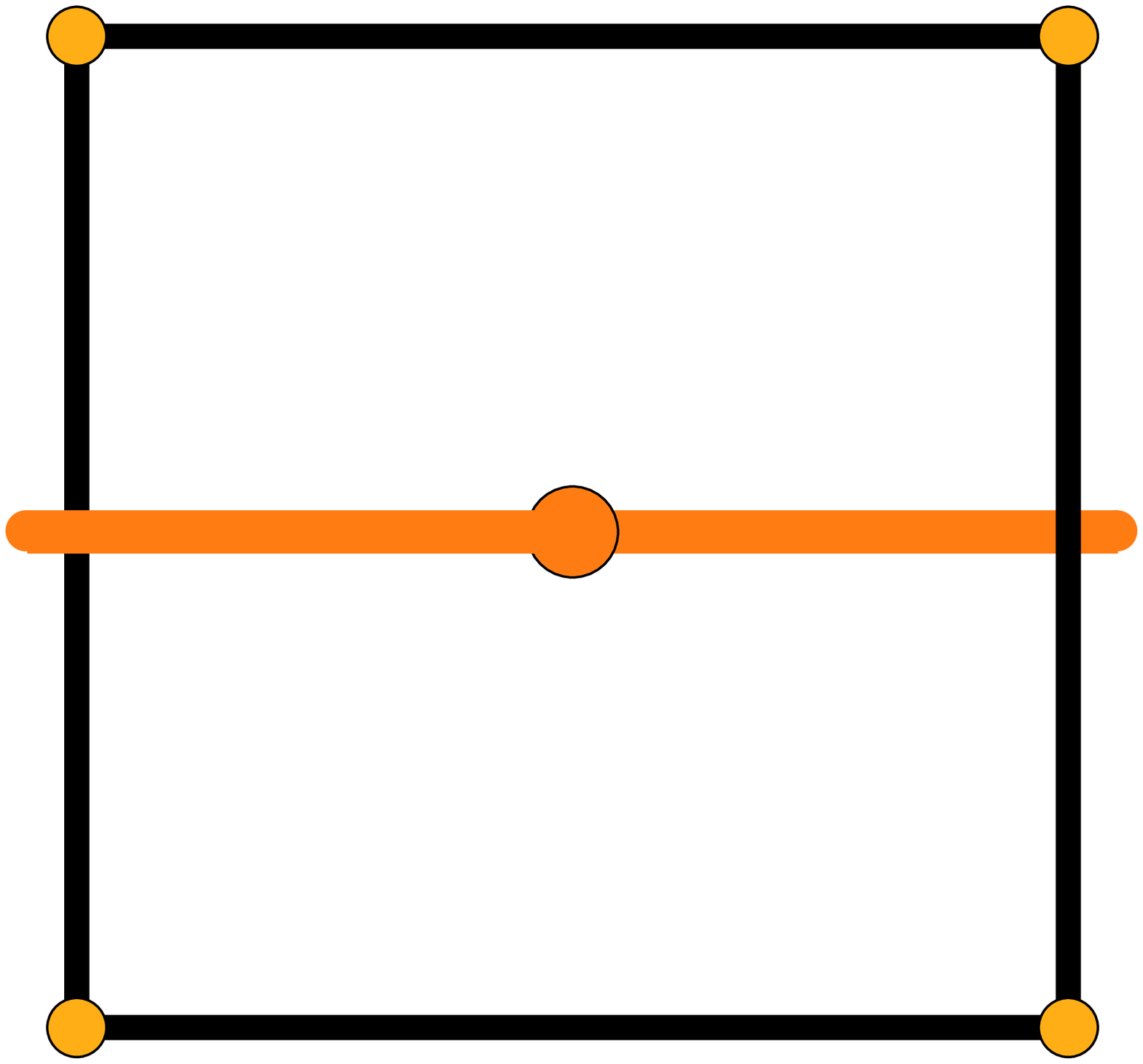}}\qquad
       \subfig{single2b}{\includegraphics[width=.1\textwidth]{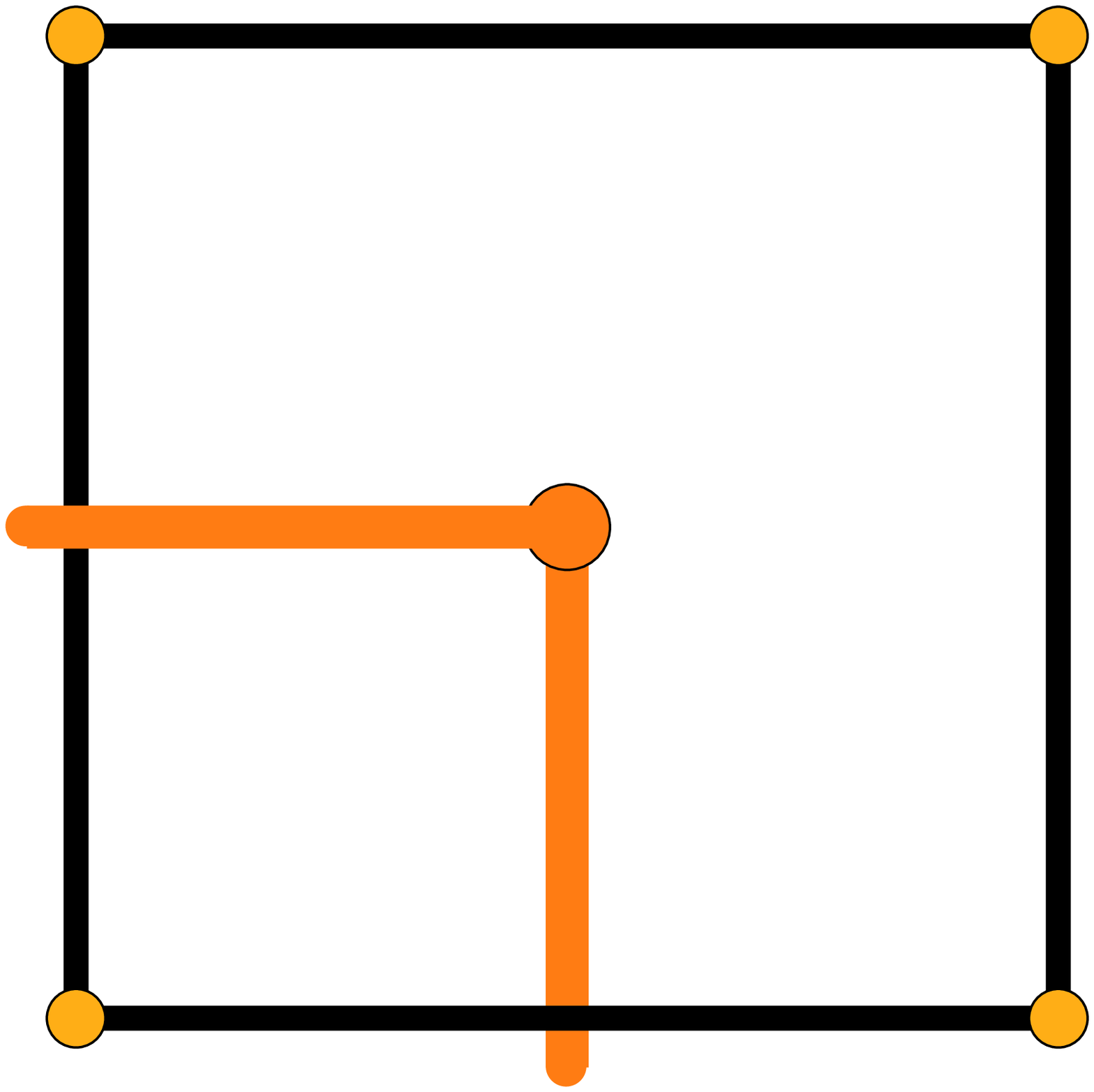}}\qquad
       \subfig{double}{\includegraphics[width=.1\textwidth]{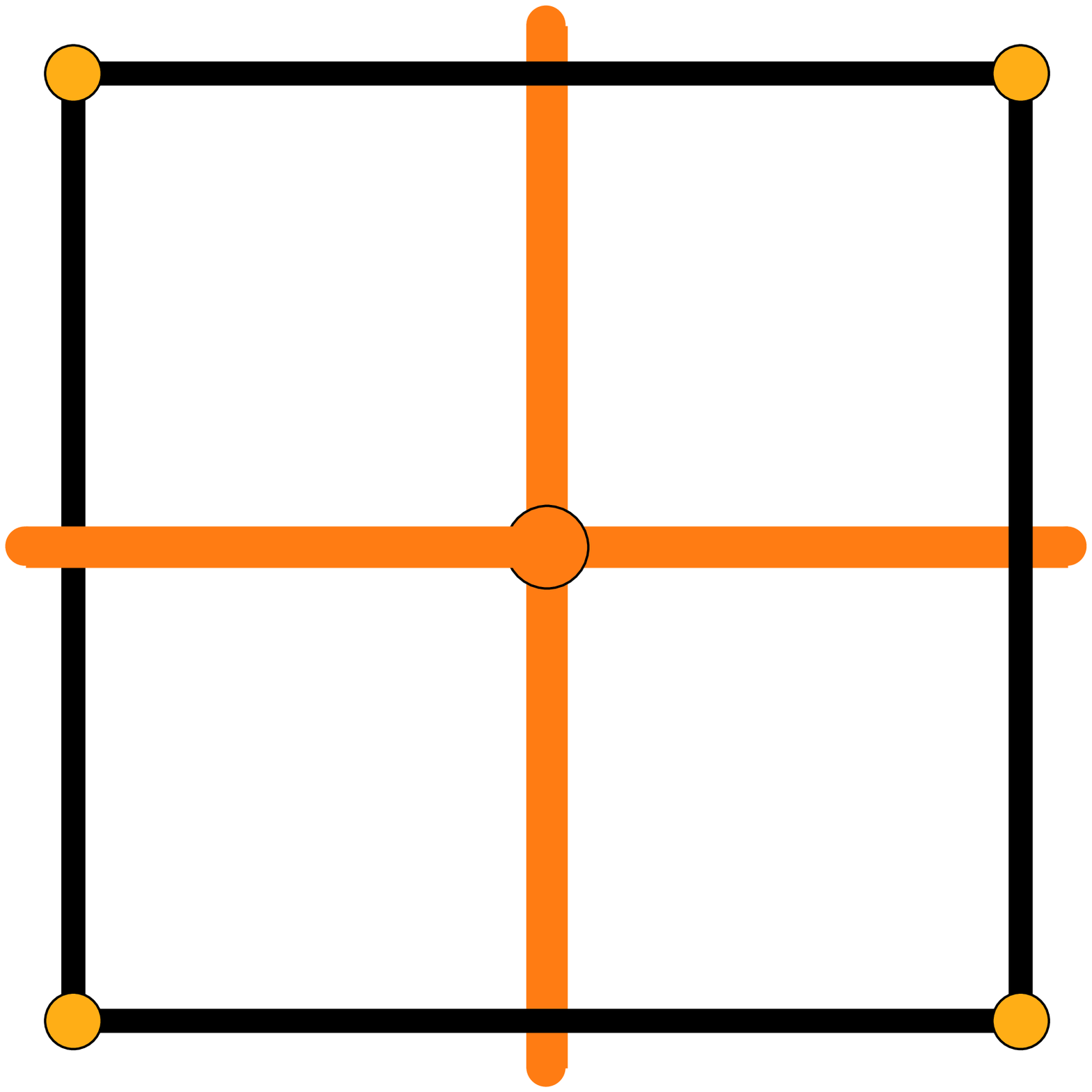}}
       \end{Figure}    
     
       \begin{enumerate}[(a)]
       \item In the cases ``single2a'' and ``double'' there is a
         coordinate hyperplane~$H$ such that $(Q,j(\mathcal{M})\cap
         Q)$ is symmetric with respect to $H$, and a vertex
         $\boldsymbol{v}$ of $\mathcal{M}$ such that
         $j(\boldsymbol{v})=\boldsymbol{w}$ and the image of the
         vertex star is contained in~$H$. Let $F$ be a facet of~$Q$
         that does not intersect~$H$.
       Then $\mathcal{B}^2_Q$ is taken to be a copy of the
       product $(\mathcal{B}^1\cap F)\times \mathcal{U}_3$:     
                 \begin{Figure}
\raisebox{1cm}{\includegraphics[height=2.5cm]{img/grid/2d_patterns.org/T3/immersion_fat}}
                   \qquad\transformsToArrow{20mm}{10mm}\qquad
                   \psfrag{0}{$\scriptstyle 0$}
                   \psfrag{1}{$\scriptstyle 1$}
                   \psfrag{1/4}{$\scriptstyle
                       \frac{1}{4}$}
                   \psfrag{5/12}{$\scriptstyle
                       \frac{5}{12}$}
                   \psfrag{1/6}{$\scriptstyle
                       \frac{1}{6}$}
 \includegraphics[height=4.5cm]{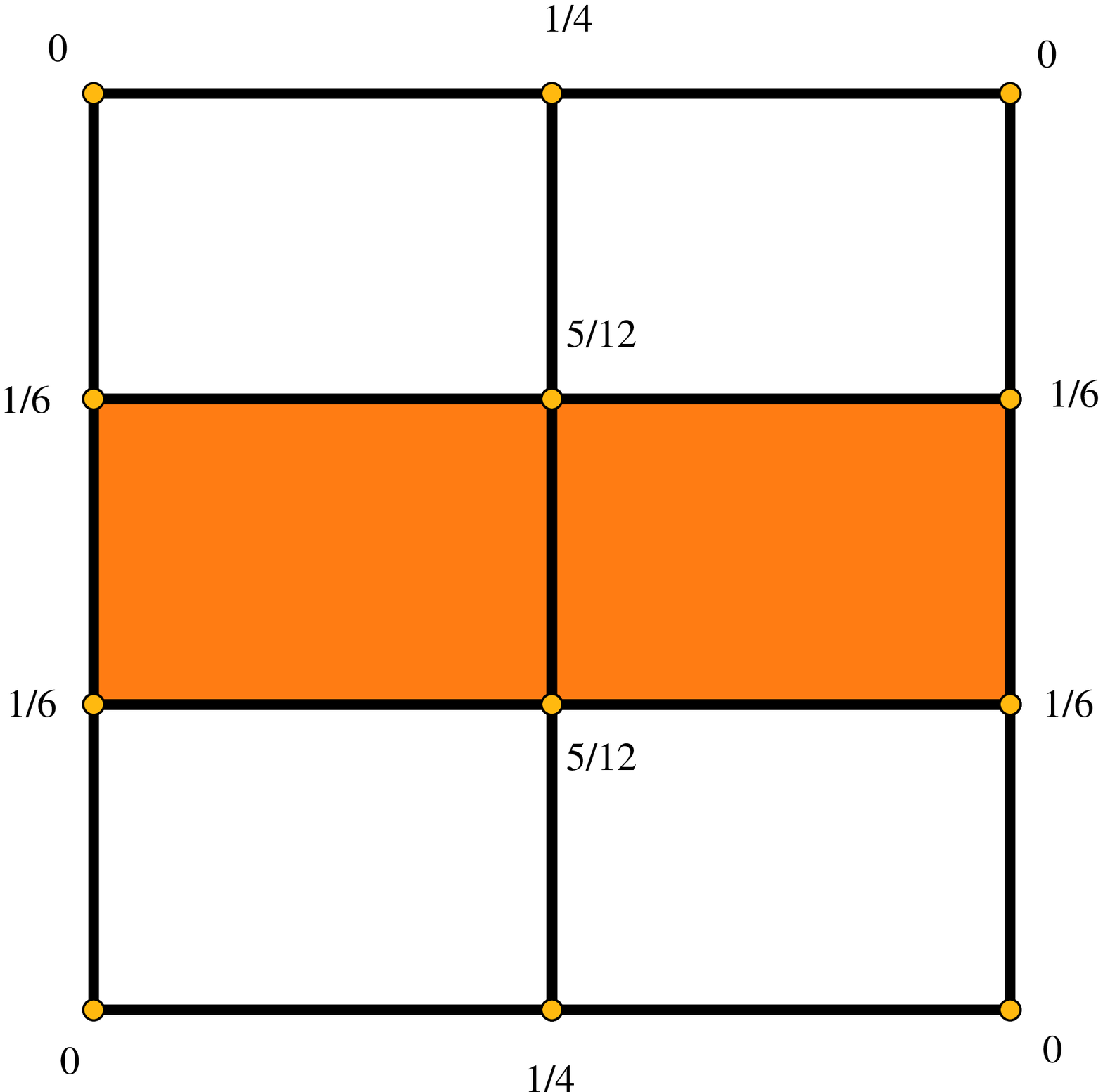}
                 \end{Figure}
\vskip-3mm
                 \begin{Figure}
\raisebox{1cm}{\includegraphics[height=2.5cm]{img/grid/2d_patterns.org/T2/immersion_fat}}
                   \qquad\transformsToArrow{20mm}{10mm}\qquad
                  \psfrag{0}{$\scriptstyle 0$}
                   \psfrag{1}{$\scriptstyle 1$}
                   \psfrag{1/3}{$\scriptstyle
                       \frac{1}{3}$}
                   \psfrag{1/6}{$\scriptstyle
                       \frac{1}{6}$}
                    \includegraphics[height=4.5cm]{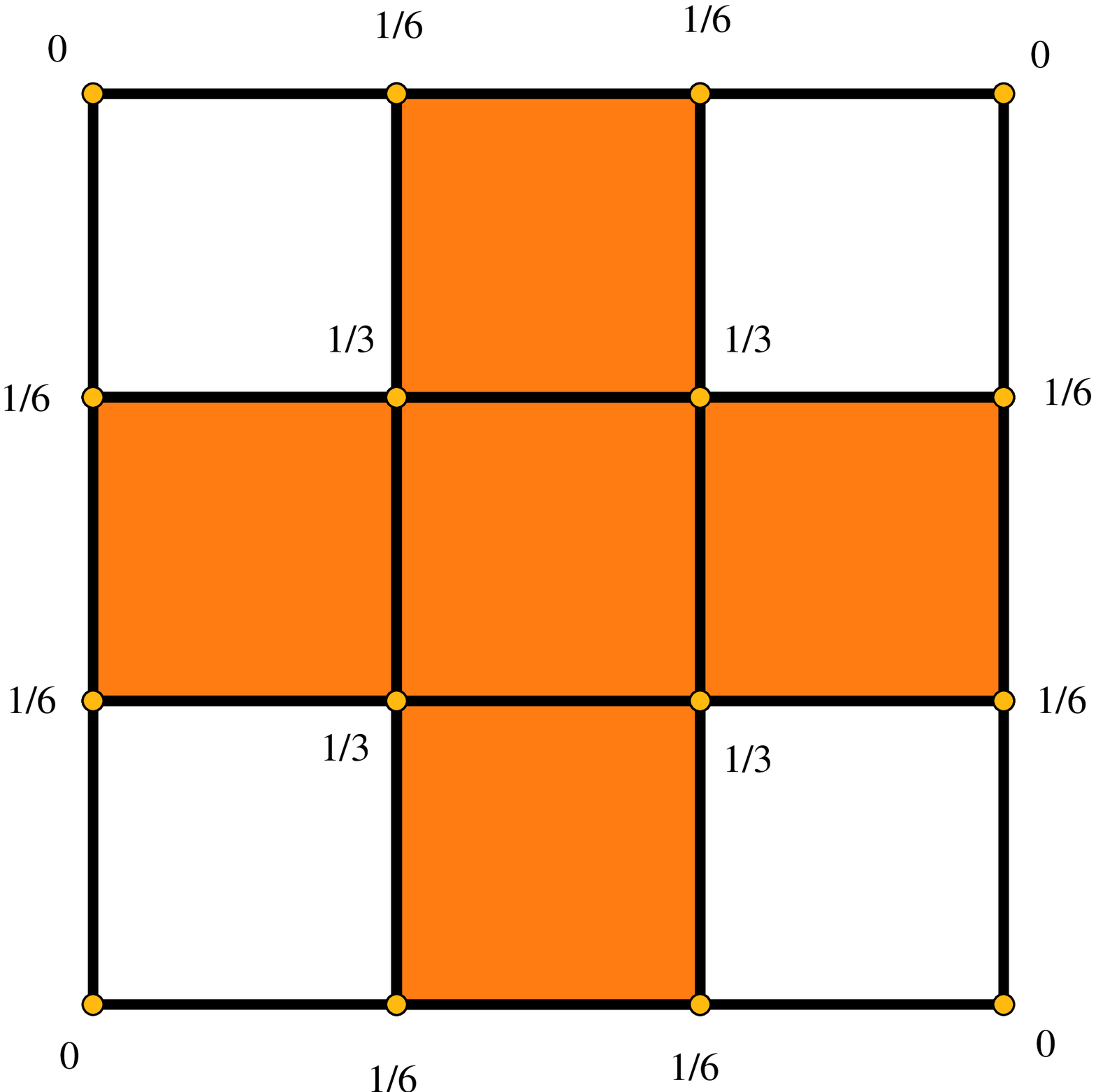}
                   \end{Figure}
                
       
                   \item If the immersion does not intersect $Q$, then
                   $\mathcal{B}^2_Q$ is a copy of the 
                   lifted cubical $2$-complex
                   $\mathcal{V}$ which
                   arises as the cubical
                   barycentric subdivision of the stellar subdivision 
                   of~$[-\tfrac{1}{2}, \tfrac{1}{2}]^2$:
\vskip-5mm
                 \begin{Figure}
                     \raisebox{1cm}{\includegraphics[height=2.5cm]{img/grid/2d_patterns.org/T1/immersion_fat}}
                     \qquad\transformsToArrow{20mm}{10mm}\qquad
                     \psfrag{0}{{$\scriptstyle 0$}}
                     \psfrag{1}{{$\scriptstyle 1$}}
                     \psfrag{1/4}{{$\scriptstyle \frac{1}{4}$}}
                     \psfrag{2/3}{{$\scriptstyle \frac{2}{3}$}}
                     \psfrag{3/4}{{$\scriptstyle \frac{3}{4}$}}
                     \includegraphics[height=4.5cm]{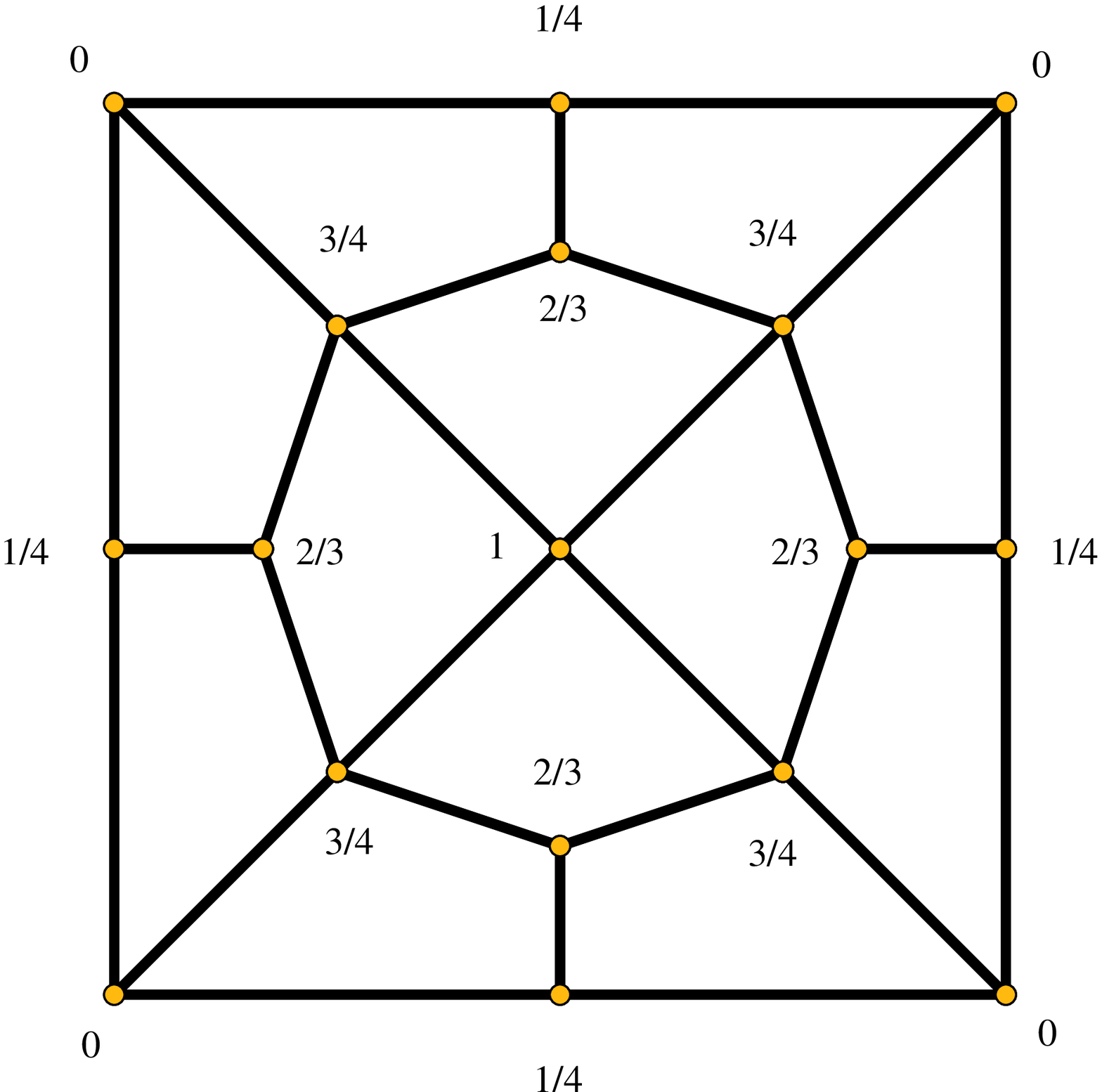}
                 \end{Figure}

               \item In the case ``single2b''  we define
                 $\mathcal{B}^2_Q$ as an affine copy of the lifted cubical
                 $2$-complex~$\mathcal{V}'$,
                 which is given by $\mathcal{V}$
                 truncated by four additional planes:
                 \begin{Figure}
                     \raisebox{1cm}{\includegraphics[height=2.5cm]{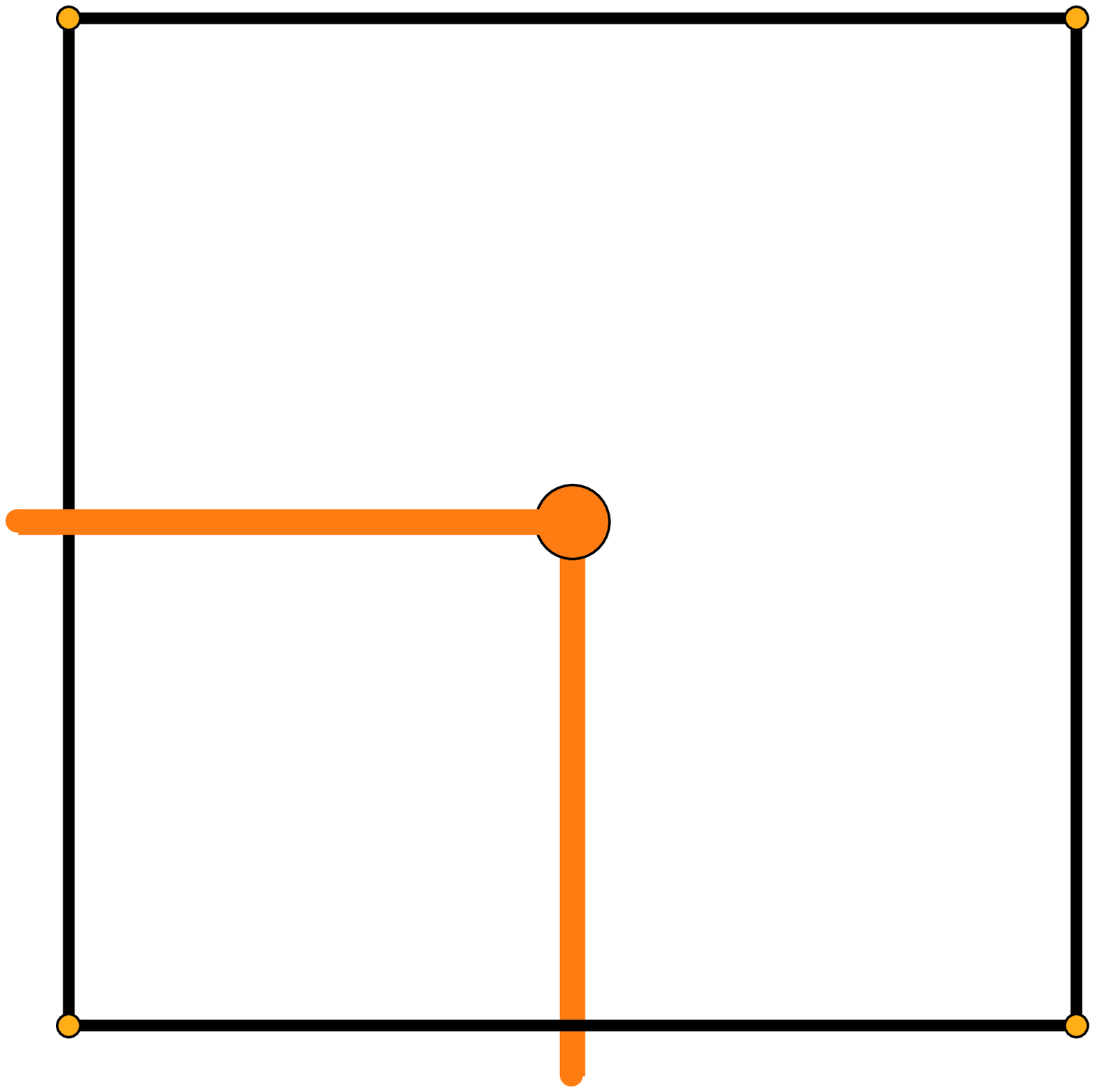}}
                     \quad\transformsToArrow{20mm}{10mm}\qquad
                     \psfrag{0}{{$\scriptstyle 0$}}
                     \psfrag{1}{{$\scriptstyle 1$}}
                     \psfrag{1/4}{{$\scriptstyle \frac{1}{4}$}}
                     \psfrag{2/3}{{$\scriptstyle \frac{2}{3}$}}
                     \psfrag{3/4}{{$\scriptstyle \frac{3}{4}$}}
                     \psfrag{1/6}{{$\scriptstyle \frac{1}{6}$}}
                     \psfrag{1/2}{{$\scriptstyle \frac{1}{2}$}} 
                     \psfrag{7/10}{{$\scriptstyle \frac{7}{10}$}}
                     \psfrag{5/6}{{$\scriptstyle \frac{5}{6}$}}
                     \includegraphics[height=4.5cm]{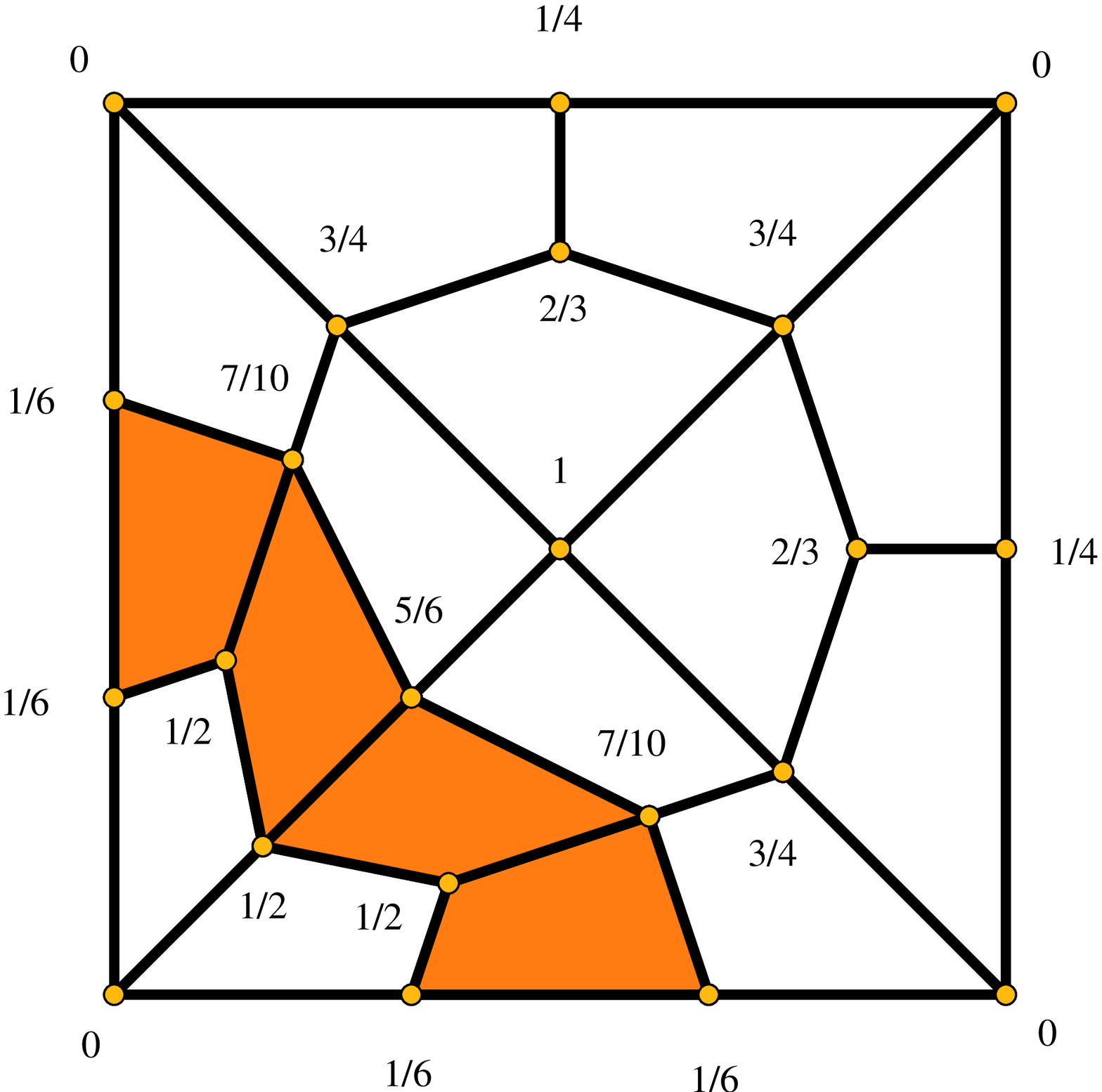}
                 \end{Figure}
                 \end{enumerate}
     Observe that the conditions
     \InvariantConsistency{2}--\InvariantDiagSubcomplex{2} are satisfied.

  \item \textbf{Subdivision of 3-cubes.}
      Let $Q$ be a facet of~$\mathcal{A}$ and $\boldsymbol{w}$ the
      unique vertex of $j(\mathcal{M})$ that is mapped to the
      barycenter of $Q$. Let $\mathcal{S}:=\mathcal{B}^2\cap Q$ be the
      induced lifted cubical boundary subdivision of~$Q$. 

   All templates for the lifted cubification $\mathcal{B}^3_Q$ of
      $\mathcal{S}$ arise either as a generalized regular Hexhoop, or as a
      product of~$\mathcal{U}_3$ with a lifted cubical subdivision of a
      facet of~$Q$.
      \begin{enumerate}[(a)]
     \item For the following four types of vertex stars we use
        a product with~$\mathcal{U}_3$:
\vskip-4mm
      \begin{Figure}\qquad\qquad
        \subfig{single4a}{
              \includegraphics[width=.15\textwidth]{img/grid/cases/single_4a}}
        \subfig{double8a}{
              \includegraphics[width=.15\textwidth]{img/grid/cases/double_8a}}
        \subfig{double8b}{
              \includegraphics[width=.15\textwidth]{img/grid/cases/double_8b}}
        \subfig{triple12}{
              \includegraphics[width=.15\textwidth]{img/grid/cases/triple}}
        \medskip
      \end{Figure}
\vskip-4mm
      In all these cases there is a coordinate symmetry plane $H$ 
      such that $H\cap Q$ is a sheet of~$j(\mathcal{M})\cap Q$.
      Hence all facets of $Q$ that
      intersect $H$ are subdivided by
      $\mathcal{U}_3\times\mathcal{U}_3$ or
      $\mathcal{U}_3\times\mathcal{U}_2$.
      Let $F$ be one of the two
      facets of $Q$ that do not intersect~$H$. Then the product
      $(\mathcal{B}^2\cap F)\times\mathcal{U}_3$ yields the lifted subdivision~$\mathcal{B}^3_Q$
      of~$Q$.
      Clearly $\mathcal{B}^3_Q$ is consistent with
     \InvariantConsistency{3} and \InvariantPLEquivalence{3}.

    \item In the remaining five cases we take a generalized regular
      Hexhoop with a diagonal plane of symmetry of~$Q$ to produce
      $\mathcal{B}^3_Q$:
\vskip-4mm
      \begin{Figure}\qquad\qquad
         \subfig{empty}{
              \includegraphics[width=.15\textwidth]{img/grid/cases/empty}}\hspace*{-2mm}
         \subfig{single3}{
              \includegraphics[width=.15\textwidth]{img/grid/cases/single_3}}
         \subfig{single4b}{
              \includegraphics[width=.15\textwidth]{img/grid/cases/single_4b}}
         \subfig{single5}{
              \includegraphics[width=.15\textwidth]{img/grid/cases/single_5}}
         \subfig{single6a}{
              \includegraphics[width=.15\textwidth]{img/grid/cases/single_6a}}
            \medskip
      \end{Figure}
\vskip-4mm      
      In each of these cases, $(Q,Q\cap j(\mathcal{M}))$ has a
      diagonal plane $H$ of symmetry. This plane
      intersects the relative interior of two facets of~$Q$. Since
      \InvariantDiagSubcomplex{2} holds, no facet
      of~$\mathcal{S}=\mathcal{B}^2\cap Q$ intersects $H$ in its
      relative interior.  By \InvariantSymmetry{2} the lifted boundary
      subdivision $\mathcal{S}$ is symmetric with respect to $H$.
      Hence all preconditions of the generalized regular Hexhoop construction
      are satisfied.  The resulting cubification $\mathcal{B}^3_Q$
      satisfies \InvariantConsistency{3} and 
      \InvariantPLEquivalence{3}.
   \end{enumerate}
  \end{enumerate}
\end{construction}

    \subsection{Correctness}

\begin{proposition}\label{prop:imm2ball}
  Let $\mathcal{M}$ be a $2$-dimensional closed
  cubical PL-surface, and
  $j\colon\mathcal{M}\immersed\mathbbm{R}^3$ a {locally symmetric}
  codimension one grid immersion.
    
  Then the cubical $3$-ball $\mathcal{B}$ given by
  Construction~\ref{constr:cubical-3ball_with_prescribed_dmf} has the
  following properties:
  \begin{thm_enumerate_i}
    \item \label{prop:imm2ball_i}
        $\mathcal{B}$ is regular, with a lifting function $\psi$.

    \item \label{prop:imm2ball_ii}
      There is a dual manifold $\mathcal{M}'$ of~$\mathcal{B}$
      and associated immersion
      $y:\mathcal{M}'\immersed\support{\mathcal{B}}$ such that
      $\mathcal{M}'$ is a cubical subdivision of $\mathcal{M}$, and
      the immersions $j$ and $y$ are PL-equivalent.
    \item \label{prop:imm2ball_iii}
      The number of facets of~$\mathcal{B}$ is congruent modulo
      two to the number~$t(j)$ of triple points of the immersion~$j$.
      
    \item \label{prop:imm2ball_iv} There is a lifted cubification
      $\mathcal{C}$ of
      $(\boundaryOP{\mathcal{B}},\psi|_{\boundaryOP\mathcal{B}})$ with
      an even number of facets.

  \end{thm_enumerate_i}
\end{proposition}

\begin{proof}
  (i) \emph{Regularity.}  By construction the lifting functions
  $\psi_F$, $F\in\facets{\mathcal{A}}$, satisfy the consistency
  precondition of the Patching Lemma (Lemma~\ref{lemma:patching_lemma}).  
  Since every pile of cube is regular
  the Patching Lemma implies that $\mathcal{B}$ is regular, too.

  (ii) \emph{PL-equivalence of manifolds} is guaranteed by 
  Property $(\text{I}_32)$.  

  (iii) \emph{Parity of the number of facets.}  For each $3$-cube $Q$
  of~$\mathcal{A}$, its cubification~$\mathcal{B}^3_Q$ is either a
  product $\mathcal{B}^2_F\times\mathcal{U}_3$ (where $\mathcal{B}^2_F$ is
  a cubification of a facet $F$ of~$Q$), or the outcome of a
  generalized regular Hexhoop construction.  In the latter case the the number of
  facets of~$\mathcal{B}^3_Q$ is even. In the first case the number of
  facets depends on the number of $2$-faces of $\mathcal{B}^2_F$. The
  number of quadrangles of $\mathcal{B}^2_F$ is odd only in the
  case ``double,'' if
  $j(\mathcal{M})\cap F$ has a double intersection point.  Hence,
  $f_3(\mathcal{B}^3_Q)$ is odd if and only if the immersion~$j$ has a
  triple point in $Q$.
  
  (iv) \emph{Alternative cubification.}  Applying
  Construction~\ref{constr:cubical-3ball_with_prescribed_dmf} to
  $\PileOfCubes{3}{\ell_1,\ell_2,\ell_3}$ without an immersed manifold 
  yields a regular cubification $\mathcal{C}$ of
  $\boundaryOP{\mathcal{B}}$ with the same lifting function as
  $\mathcal{B}$ on the boundary.
  Since the immersion $\emptyset\immersed\R^3$ has no triple points
  the number of facets of $\mathcal{C}$ is even.
\end{proof}

\subsection{Proof of the main theorem}

\begin{proof}[Proof of Theorem~\ref{thm:C4P_with_prescribed_dmf}]
  Let $j\colon\mathcal{M}\immersed\mathbbm{R}^3$ be a locally flat
  normal crossing immersion of a compact $(d-1)$-manifold
  $\mathcal{M}$ into $\R^{d}$.
  By Proposition~\ref{prop:PL_to_grid_imm} and
  Proposition~\ref{prop:local_symmetry} there is a cubical subdivision
  $\mathcal{M}'$ of $\mathcal{M}$ with a locally symmetric, codimension
  one grid immersion $j'\colon{\mathcal{M}}'\immersed\mathbbm{R}^3$
  that is PL-equivalent to~$j$.

  Construct a convex cubical $3$-ball $\mathcal{B}$ with
  prescribed dual manifold immersion $j'$ as described above. By
  Proposition~\ref{prop:imm2ball}(\ref{prop:imm2ball_i}) the ball
  $\mathcal{B}$ is regular, and by
  Proposition~\ref{prop:imm2ball}(\ref{prop:imm2ball_iv}) there is a
  cubification~$\mathcal{C}$ of~$\boundaryOP{\mathcal{B}}$ with an
  even number of facets and the same lifting function on the boundary.

  The lifted prism\index{lifted prism} over $\mathcal{B}$ and
  $\mathcal{C}$ (Construction~\ref{constr:Lifted_prism_over_two_balls}
  of Section~\ref{subsec:cubical_prisms_two_balls}).  This yields a
  cubical $4$-polytope~$P$ with
  \[
      f_3(P) = f_3(\mathcal{B}) + f_3(\mathcal{C}) +
      f_2(\boundaryOP\mathcal{B}),
  \]
  whose boundary contains $\mathcal{B}$ and thus has a dual manifold
  immersion PL-equivalent to~$j$.

  For (\ref{main_theorem:parity_of_f_3}) observe that
  for every cubical $3$-ball the number of facets of the boundary
  is even. Hence $f_2(\boundaryOP\mathcal{B})$ is even.  Since the number
  of facets of~$\mathcal{C}$ is even, we obtain
  \[
      f_3(P)\ \equiv\ f_3(\mathcal{B}) \ \equiv\ t(j)\ \ \bmod 2.
  \]

  Now consider (\ref{main_theorem:odd}).
  By a famous theorem of Banchoff \cite{Banchoff1974} the number of
  triple points of a normal crossing codimension one immersion of a
  surface has the same parity as the Euler characteristic.  Hence, if
  $\mathcal{M}$ is a non-orientable surface of odd genus the
  number of triple points of $j$ is odd, which implies that the cubical
  $4$-polytope $P$ has an odd number of facets.  
\end{proof}

\subsection{Symmetric templates}\label{subsec:sym_templates}

The three-dimension templates constructed above, which we call the
\emph{standard templates}, do not satisfy the conditions
\InvariantSymmetry{3} and \InvariantDiagSubcomplex{3}. In particular,
the symmetry requirement \InvariantSymmetry{3} is violated by the templates
corresponding to the cases ``empty'', ``single3'', and ``single6a,''
and it is satisfied by all others.
For example, the standard template for ``single5'' is illustrated in
Figure~\ref{fig:dmf_three_dim_gen_Hexhoop}; it
satisfies \InvariantSymmetry{3} since there is only one diagonal
symmetry hyperplane.

For the ``empty'' case an alternative template
may be obtained from the cubical barycentric subdivision. The resulting
cubification satisfies both conditions \InvariantSymmetry{3} and
\InvariantDiagSubcomplex{3}, and furthermore, it has less faces --- 96
facets, 149 vertices --- than the standard template.
\\
For the case ``single3'' an alternative cubification, of full symmetry,
can be constructed from $\mathcal{C}''$ by truncating the lifted
polytope corresponding to the lifted cubical ball $\mathcal{C}''$ by
some additional hyperplanes.   
\\
For the case ``single6a'' we do not know
how to get a cubification of full symmetry.
This is the main obstacle for an extension of our constructions
to higher dimensions (cf.~Section~\ref{sec:nextstep}).



\section{An odd cubical 4-polytope with a dual Boy's surface}\label{sec:boy}

  Cubical $4$-polytopes with odd numbers of facets exist
  by our Main Theorem~\ref{thm:C4P_with_prescribed_dmf}.
In this section we describe the construction of a
cubical $4$-polytope with an odd number of facets in more
detail. The data for the corresponding model will be submitted
to the \url{eg-models} archive. 

\begin{theorem}\label{thm:odd_C4P}  
  There is a cubical $4$-polytope $P^{}_{\mathrm{Boy}}$ with $f$-vector
  \[  
     f\ =\ (17\,718, 50\,784, 49\,599, 16\,533)
  \]
  that has has a Boy surface as a dual manifold.
\end{theorem}

\subsection*{A grid immersion of Boy's surface.} 
%
The construction starts with a grid immersion (cf.~\cite{Petit1995})
of Boy's surface, that is, an immersion of the real projective plane with
exactly one triple point and three double-intersection curves in a
pattern of three loops \cite{Boy1903} \cite{HilbertCohn-Vossen} \cite{Apery1987}. 
This immersion $j:\mathcal{M}\immersed\R^3$ in shown in Figure~\ref{grid_boy}.

\begin{Figure}\label{grid_boy}
  \subfigure[grid immersion]{\includegraphics[height=.45\textwidth]{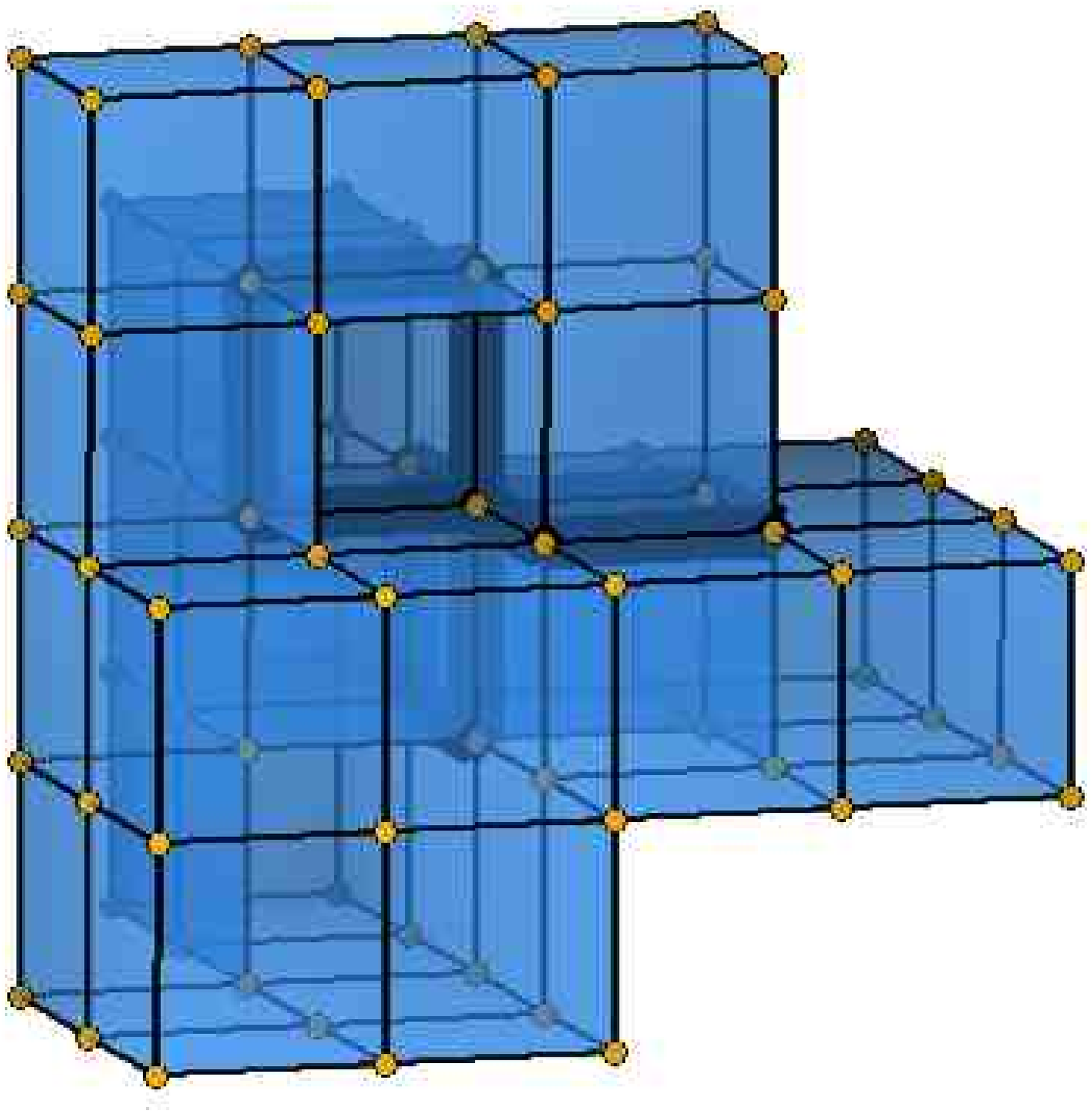}}\qquad
  \subfigure[]{\includegraphics[height=.45\textwidth,clip=true]{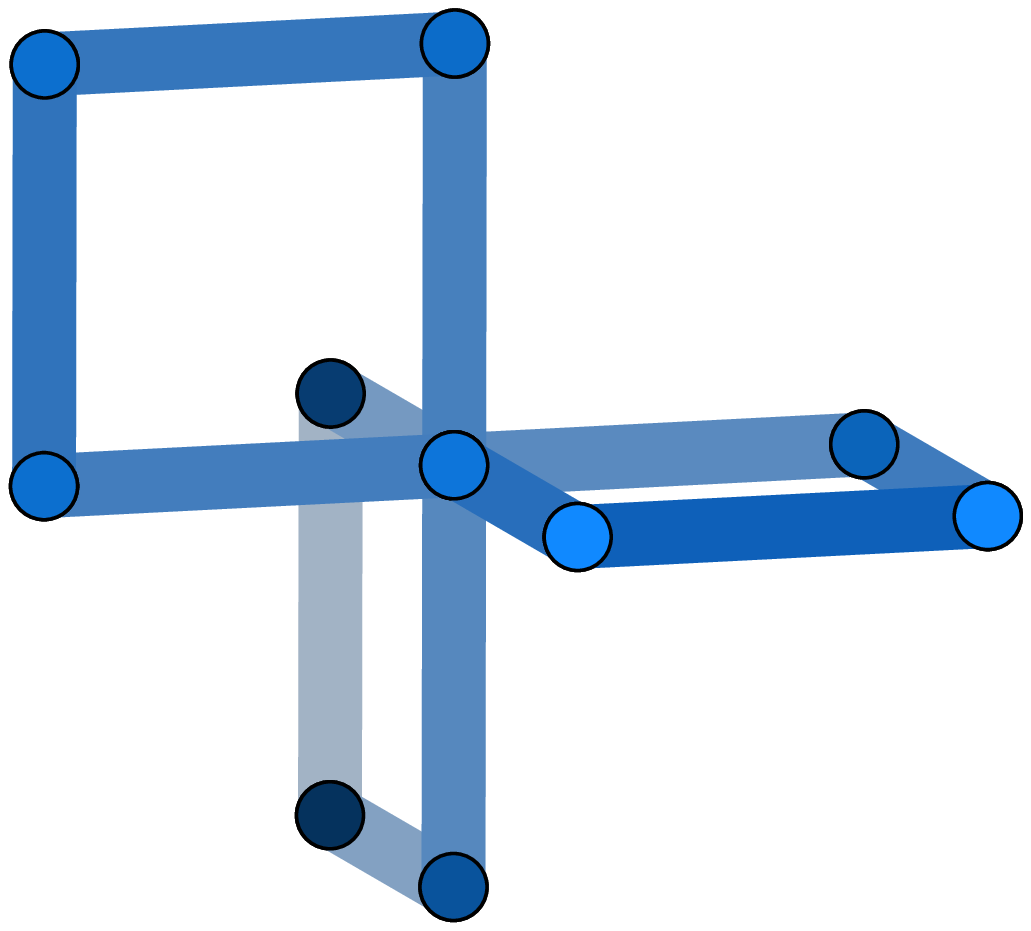}}
\vskip-5mm
  \caption[A grid immersion of Boy's surface]{
    A grid immersion the Boy's surface.
    Each double-intersection loop is of length four.}
  \label{fig:grid_boy}
\end{Figure}

The $2$-manifold $\mathcal{M}$ has the $f$-vector $f(\mathcal{M})\ =\ 
(85,\, 168,\, 84),$ whereas the image of the grid immersion has the
$f$-vector $f(j(\mathcal{M}))\ =\ (74,\,156,\,84)$. The vertex
coordinates can be chosen such that the image $j(\mathcal{M})$ is
contained in a pile of cubes $\PileOfCubes{3}{4,4,4}$.

\subsection*{A cubical 3-ball with a dual Boy's surface.}

We apply Construction~\ref{constr:cubical-3ball_with_prescribed_dmf} to the
grid immersion $j:\mathcal{M}\immersed\R^3$ to obtain a cubical
$3$-ball with a dual Boy's surface, and with an odd number of facets.

\begin{figure}[!h]\centering
  \includegraphics[height=.8\textwidth]{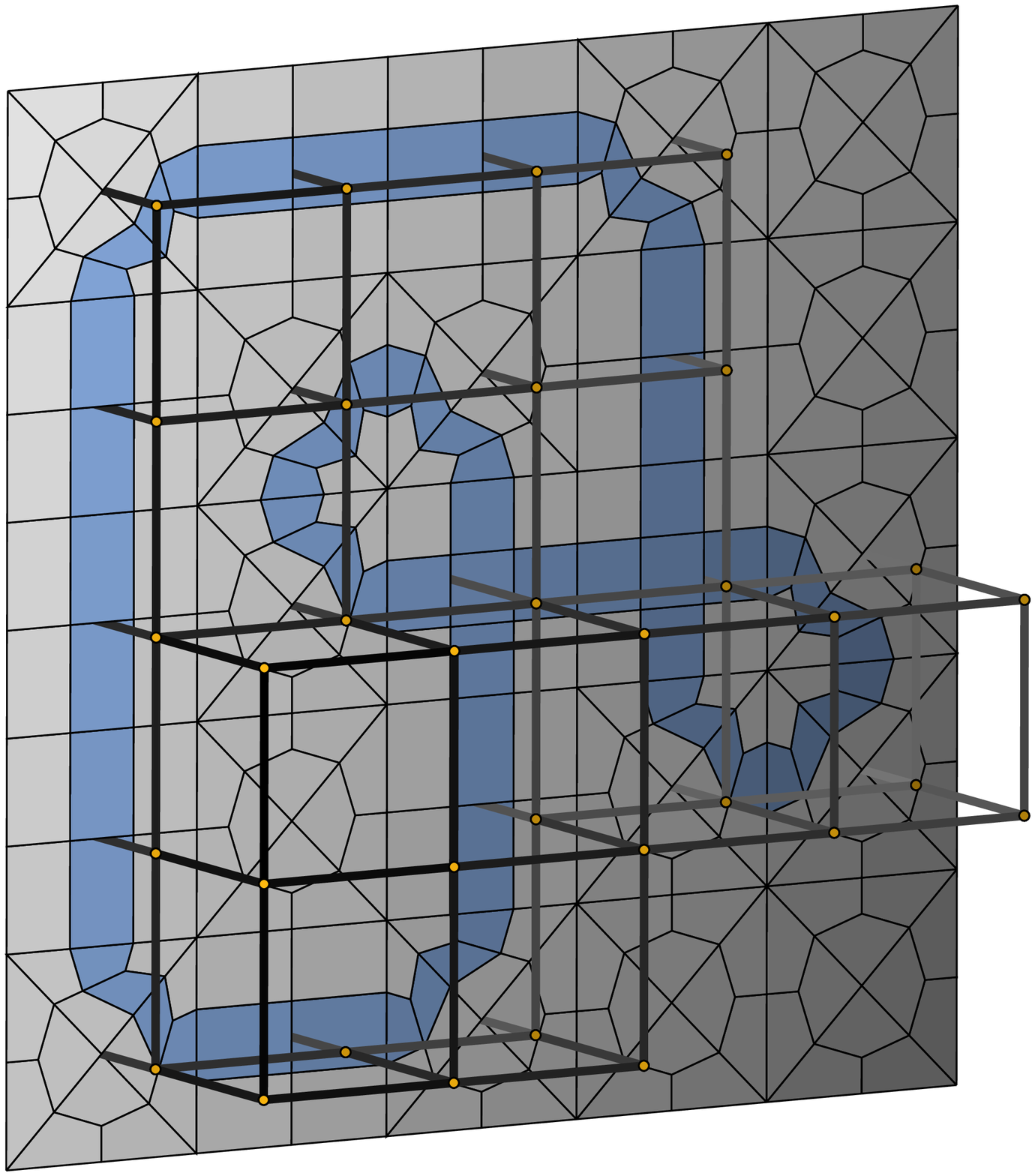}
  \caption{A sketch of the cubification of the $2$-skeleton of $\mathcal{A}$.}
  \label{fig:grid_boy_2_skel}
\end{figure}

Since the image $j(\mathcal{M})$ is contained in a pile of cubes
$\PileOfCubes{3}{4,4,4}$, the raw complex $\mathcal{A}$ given by
Construction~\ref{constr:cubical-3ball_with_prescribed_dmf} is
isomorphic to $\PileOfCubes{3}{5,5,5}$. 
Hence we have $5^3 - 74 = 51$ vertices of
$\mathcal{A}$ that are not vertices of~$j(\mathcal{M})$.
We try to give an impression of the subdivision $\mathcal{C}^2$ of the
$2$-skeleton of~$\mathcal{A}$ in Figure~\ref{fig:grid_boy_2_skel}.
The $f$-vector of~$\mathcal{C}^2$ is $f=(4\,662,\,  9\,876,\, 5\,340)$.

The subdivision of the boundary of~$\mathcal{A}$ consists of
$150=6\cdot 5\cdot 5$ copies of the two-dimensional ``empty pattern''
template.  Hence the subdivision of the boundary of~$\mathcal{A}$
(given by $\mathcal{C}^2 \cap \support{\boundaryOP{\mathcal{A}}}$) has
the $f$-vector $f=(1\,802,\, 3\,600,\, 1\,800)$.

The refinement $\mathcal{B}$ of~$\mathcal{A}$
depends on templates that are used for the $3$-cubes. We use
the ``symmetric'' templates of Section~\ref{subsec:sym_templates}. 
The $f$-vector
of~$\mathcal{B}$ is then $f=(15\,915,\, 45\,080,\, 43\,299,\, 14\,133)$.
(The ``standard set'' of templates yields a cubical ball 
with $18\,281$ facets.)

\begin{figure}[H]\centering
  \includegraphics[width=.8\textwidth]{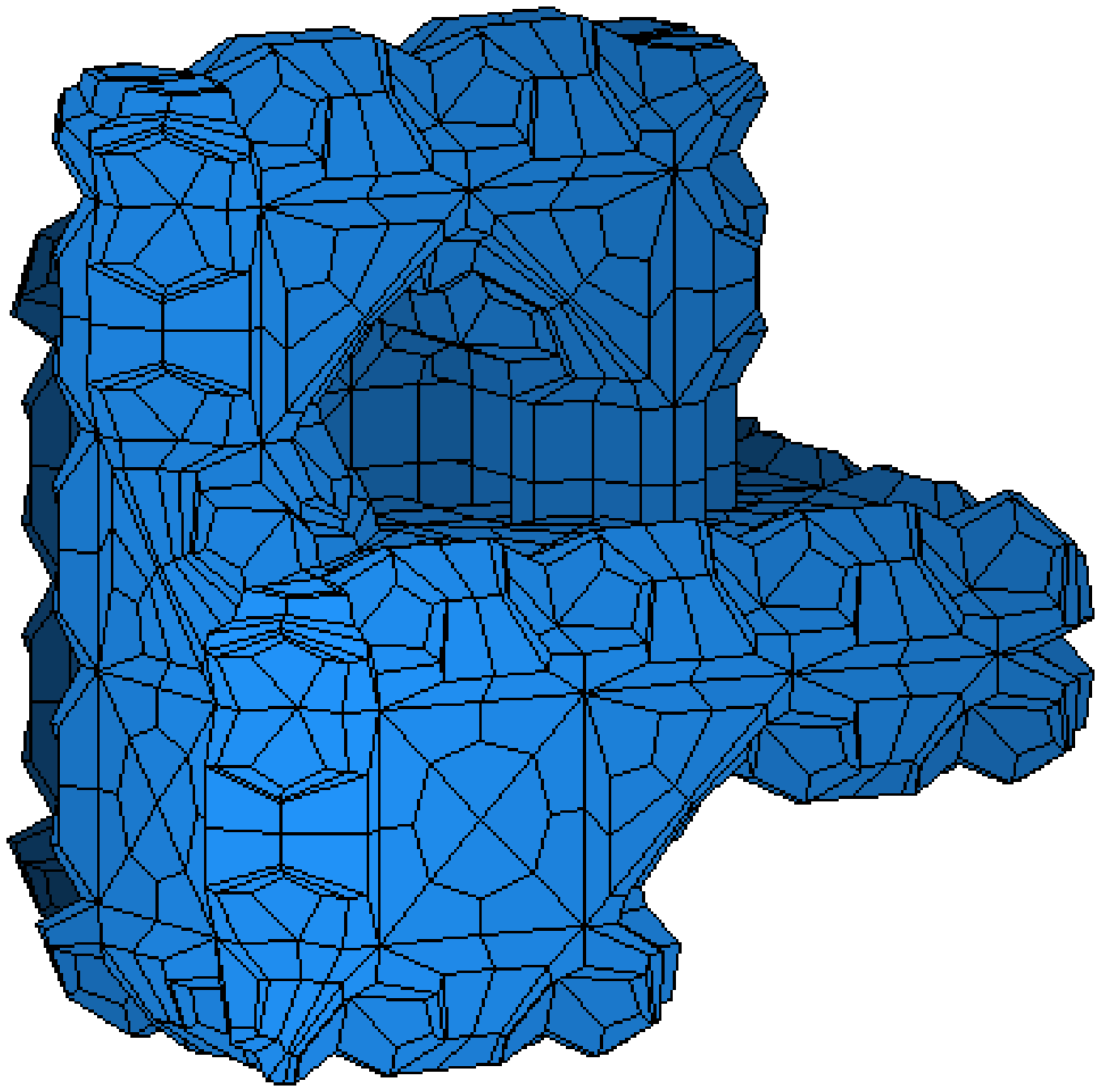}
  \caption[Dual Boy's surface]{The dual Boy's surface of $f$-vector $f=(1\,998,\,  3\,994,\, 1\,997)$ of the cubical $3$-ball $\mathcal{B}$.}
  \label{fig:dual_boys}
\end{figure}

Figure~\ref{fig:dual_boys} illustrates the dual Boy's surface of
the cubical $3$-ball $\mathcal{B}$. It has the $f$-vector
$f=(1\,998,\, 3\,994,\, 1\,997)$; its multiple-intersection
loops have length~$16$.
The ball $\mathcal{B}$ has  612 dual manifolds is total 
(339 of them without boundary).


\begin{figure}[H]\centering
  \includegraphics[width=.4\textwidth]{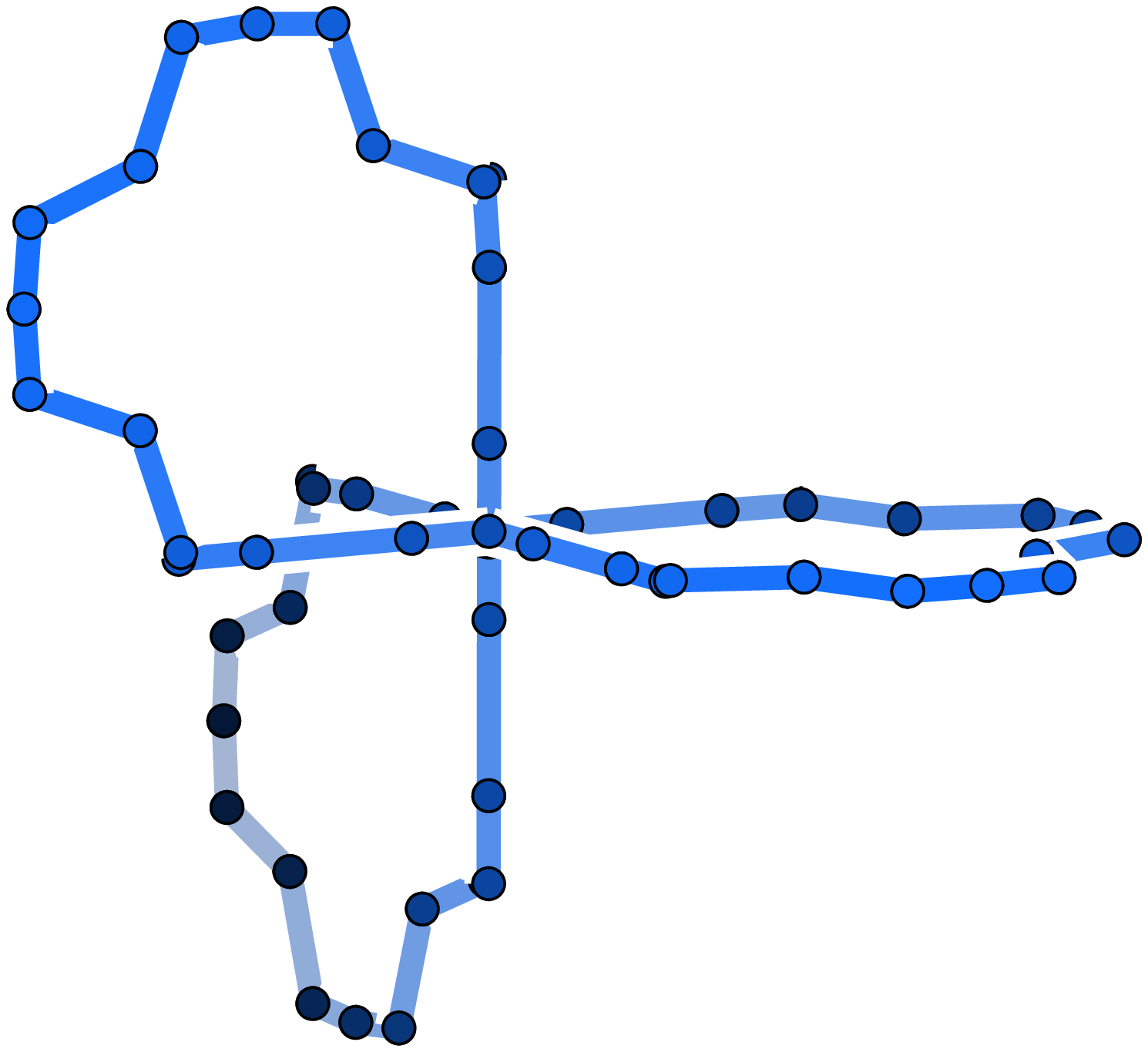}
  \caption[Intersection quadrangles and curve of dual Boy's surface]{The multiple-intersection curve of the dual Boy's surface of the cubical $3$-ball $\mathcal{B}$.}
  \label{fig:dual_boys_propellor}
\end{figure}
\pagebreak

\subsection*{A cubical 4-polytope with a dual Boy's surface}

A cubification $\mathcal{B}'$ of $\boundaryOP{\mathcal{B}}$ with an
even number of facets is given by subdividing each facet of the raw
ball $\mathcal{A}$ with a cubification for the empty pattern. Using
the symmetric cubification for the empty pattern yields a regular
cubical $3$-ball $\mathcal{B}'$ with $12\,000$ facets.  The lifted
prism over $\mathcal{B}$ and~$\mathcal{B}'$ yields a cubical
$4$-polytope with $27\,933$ facets.

%

However, using first stellar subdivisions, then a (simplicial) cone,
and then cubical subdivisions on
${\boundaryOP\PileOfCubes{3}{5,5,5}}$, it is possible to produce a
significantly smaller alternative cubification $\mathcal{B}''$ of
$\boundaryOP\mathcal{B}$ with an even number of facets.
Moreover, for this one can form the cone based directly on the
boundary complex of $\mathcal{B}$ and thus ``save the vertical part''
of the prism.  The resulting cubical $4$-polytope
$P^{}_{\mathrm{Boy}}$ has $f_0=17\,718$ vertices and $f_3=16\,533$
facets.  A further analysis of the dual manifolds
of~$P^{}_{\mathrm{Boy}}$ shows that there are 613 dual manifolds in
total: One dual Boy's surface of $f$-vector $f=(1\,998,\, 3\,994,\,
1\,997)$, one immersed surface of genus~20 (immersed with 104 triple
points) with $f$-vector $(11\,470,\ 23\,016, 11\,508)$, and 611
embedded $2$-spheres with various distinct $f$-vectors.

\subsection*{Verification of the instances}

All the instances of the cubical $4$-polytopes described above were
constructed and verified as electronic geometry models in the
\textsf{polymake} system by \textsc{Gawrilow \& Joswig}
\cite{polymake}, which is designed for the construction and analysis
of convex polytopes.  A number of our own tools for handling cubical
complexes are involved as well.  These cover creation, verification,
and visualization of cubical complexes (for $d\in\{2,3\}$).

The instances are available from
   \url{http://www.math.tu-berlin.de/~schwartz/c4p}.
   
   Whereas the construction of the instances involves new tools that
   were writted specifically for this purpose, the verification
   procedure uses only standard polymake tools.  All tools used in the
   verification procedure are parts of polymake system which have been
   used (and thereby verified) by various users over the past years
   (using a rich variety of classes of polytopes).

The topology of the dual manifolds of our instances was examined
using all the following tools:
\begin{compactitem}[~$\bullet$]
\item A homology calculation code based  written by
  Heckenbach \cite{Heckenbach1997}.
\item The \textsf{topaz} module of
  the \textsf{polymake} project, which covers the
  construction and analysis of simplicial complexes.
  \item Our own tool for the calculation of the Euler characteristics.
\end{compactitem}


\section{Consequences}\label{sec:conseq}\label{sec:consequences}
\label{sec:C4P_with_prescribed_nor_dmfs_even}

In this section we list a few immediate corollaries and consequences
of our main theorem and of the constructions that lead to it.
The proofs are quite immediate, so we do not give extended explanations
here, but refer to \cite{Schwartz2} for details.

\subsection{Lattice of f-vectors of cubical 4-polytopes}
  
Babson \& Chan \cite{BabsonChan3} have obtained a 
characterizazion of the
$\mathbbm{Z}$-affine span of the $f$-vectors of cubical $3$-spheres:
With the existence of cubical $4$-polytopes with an odd number of facets
this extends to cubical $4$-polytopes.

\begin{corollary}\label{cor:span_of_f_vectors_of_C4P}
  The $\mathbbm{Z}$-affine span of the $f$-vectors
  ($f_0,f_1,f_2,f_3)$ of the cubical $4$-polytopes is
  characterized by
\begin{compactenum}[\rm(i)]
\item integrality ($f_i\in\mathbbm{Z}$ for all~$i$), 
\item the cubical Dehn-Sommerville equations $f_0-f_1+f_2-f_3=0$ 
and $f_2=3f_3$,
and
\item the extra condition $f_0\equiv0\bmod2$.
 \end{compactenum}
\end{corollary}

Note that 
this includes modular conditions such as $f_2\equiv 0\bmod 3$,
which are not ``modulo~$2$.'' 
The main result of Babson \& Chan \cite{BabsonChan3} says that
for cubical $d$-spheres and $(d+1)$-polytopes, $d\ge2$,
``all congruence conditions are modulo~$2$.'' However, this
refers only to the modular conditions 
\emph{which are not implied by integrality
and the cubical Dehn-Sommerville equations}.
The first example of such a condition is, for $d=4$,
the congruence (iii) due to Blind \& Blind~\cite{BlindBlind94}.

\subsection
{Cubical 4-polytopes with dual manifolds of prescribed genus}

By our main theorem, from any
embedding $M_g\rightarrow\R^3$
we obtain a cubical $4$-polytope that has the
orientable connected $2$-manifold $M_g$ of genus
$g$ as an embedded dual manifold.
Indeed, this may for example 
be derived from a grid embedding of 
  $M_g$ into the pile of cubes $P(1,3,1+2g)$.

  However, cubical $4$-polytopes with an orientable dual
     manifold of prescribed genus can 
much more efficiently, and with more control on the 
topological data, be produced by means of
     connected sums of copies of the ``neighborly cubical''
     $4$-polytope~$C^5_4$ with the graph of a $5$-cube (compare
     Section~\ref{subsec:NCP}).

\begin{proposition}\label{lemma:C4P_with_dmf_of_prescr_genus}
     For each $g>0$, there is a cubical $4$-polytope $P_g$ with the following properties.
     \begin{compactenum}[\ \rm(i)]
     \item The polytope $P_g$ has exactly one embedded orientable dual
       $2$-manifold $\mathcal{M}$ of genus~$g$ with $f$-vector
       $f(\mathcal{M})=(12g+4,\, 28g+4,\,14g+2)$.  
     \item There is a facet $F$ of $P$ which is not intersected by the
       image of the dual manifold~$\mathcal{M}$, and which is
       \emph{affinely regular}, that is, there is an affine transformation between $F$ and
       the standard cube $[-1,+1]^3$.
        \item All other dual manifolds of $P_g$ are embedded $2$-spheres.
        \item $f(P_g) = (24 g + 8,\, 116g+12,\, 138g+6,\,46 g + 2)$.
      \end{compactenum}
\end{proposition}

Taking now the connected sum of one example of 
a $4$-polytope with a non-orientable dual $2$-manifold,
we obtain $4$-polytopes with a non-orientable dual manifold
of prescribed genus.

  \begin{corollary} 
     For each even $g>0$, there is a cubical $4$-polytope that has a
      cubation of the non-orientable connected $2$-manifold $M_g'$ of
      genus $g$ as a dual manifold (immersed without triple points and with 
      one double-intersection curve). 
   \end{corollary} 
   
\begin{figure}[ht]\centering
          \includegraphics[height=50mm]{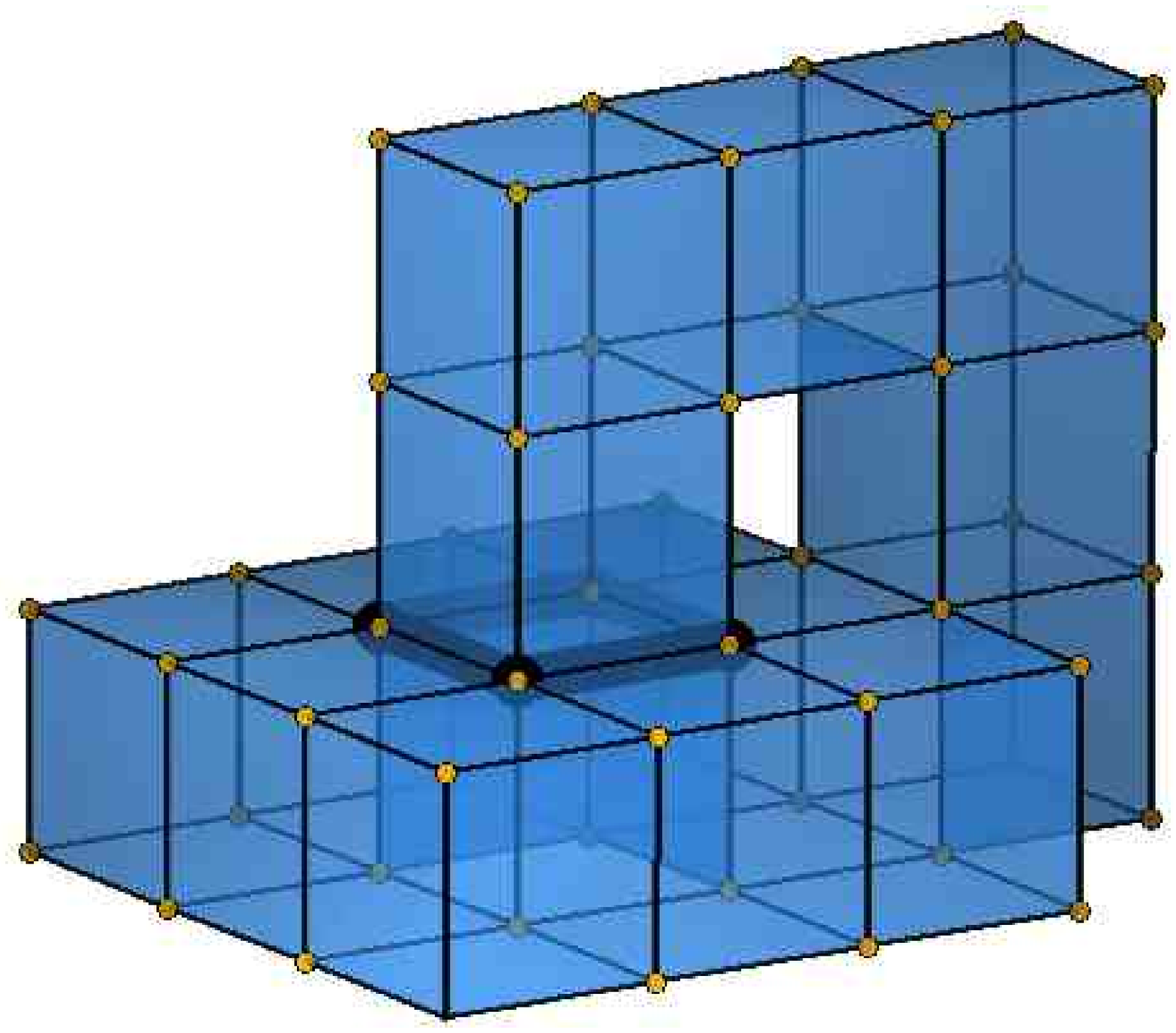}\quad 
           \caption[A grid immersion of the Klein bottle]{
             A grid immersion of the
             Klein bottle with one double-intersection
             curve and without triple points.}
           \label{fig:grid_klein}
     \end{figure}

For this, one can for example construct the $4$-polytope associated 
with the  grid immersion of the Klein bottle of
  $f$-vector~$f=(52,108,56)$
  as depicted in Figure~\ref{fig:grid_klein}.

Smaller cubical $4$-polytopes with non-orientable cubical
$4$-polytopes can be produced by means of connected sums of the
cubical $4$-polytope $P_{62}$ of Section~\ref{sec:klein} with a dual
Klein bottle, and several copies of the neighborly cubical
$4$-polytope $C^5_4$. (Some ``connector cubes'' of
Lemma~\ref{lem:connector_cube} have to be used.) 
The resulting cubical $4$-polytope has rather small $f$-vector entries, but
the set of multiple-intersection points consists of five double-intersection curves.
 
Applying the same proof as above to the grid immersion of Boy's
surface of the previous section yields the following result.

  \begin{corollary} 
    For each odd $g>0$, there is a cubical $4$-polytope that has a
    cubation of the non-orientable connected $2$-manifold $M_g'$ of
    genus $g$ as a dual manifold (immersed with one triple point and
    three double-intersection curves of length $14$).
   \end{corollary}

\subsection[Higher-dimensional cubical polytopes with non-orientable dual manifolds]{%
            Higher-dimensional cubical polytopes with non-orientable dual manifolds}

\begin{corollary}
For each $d\ge4$ there are cubical $d$-polytopes
with non-orientable dual manifolds.
\end{corollary}

\begin{proof}
  By construction, the $4$-dimensional instance $P_{62}$ of
  Section~\ref{sec:klein} comes with a subdivision into a regular
  cubical $4$-ball.  Since one of its dual manifolds is not
  orientable, its $2$-skeleton is not edge orientable, i.\,e.\ it
  contains a cubical M\"obius strip with parallel inner edges.  So if
  we now iterate the lifted prism construction\index{lifted prism} of
  Section~\ref{subsec:lifted_prisms}, then the resulting cubical
  $d$-polytopes ($d\ge4$) will contain the $2$-skeleton of~$P_{62}$.
  By Proposition~\ref{prop:NEO} they must also have non-orientable
  dual manifolds.
\end{proof}


\section{Applications to hexa meshing}\label{sec:hexameshing}  

In the context of computer aided design (CAD) the surface of a
workpiece (for instance a part of a car, ship or plane) is often
modeled by a \emph{surface mesh}. In order to analyze physical and
technical properties of the mesh (and of the workpiece),
finite element methods (FEM) are widely used.

Such a surface mesh is either a \emph{topological mesh}, that is,
$2$-dimensional regular CW complex, or a
\emph{geometric mesh}, that is, a (pure) $2$-dimensional polytopal
complex cells.  
Common cell types of a surface mesh are triangles
($2$-simplices) and quadrangles. Thus a 
\emph{geometric quad mesh} is a
$2$-dimensional cubical complex, and a 
\emph{topological} one is a cubical $2$-dimensional regular CW complex.

In recent years there has been growing interest in volume meshing.
Tetrahedral volume meshes (simplicial $3$-complexes) are
well-understood, whereas there are interesting and challenging open
questions both in theory and practice for hexahedral volume meshes,
\emph{hexa meshes} for short. 
That is, a \emph{geometric hexa mesh} is a $3$-dimensional cubical
complex, and a \emph{topological hexa mesh} is a cubical $3$-dimensional
regular CW complex.

A challenging open question in this context is whether each cubical
quadrilateral geometric surface mesh with an even number of
quadrangles admits a geometric hexa mesh.  In our terminology this
problem asks whether each cubical PL $2$-sphere with an even
number of facets admits a cubification.  Thurston \cite{Thurston1993}
and Mitchell \cite{Mitchell1996} proved independently that every
topological quad mesh with an even number of quadrangles admits a
topological hexa mesh. Furthermore, Eppstein showed in
\cite{Eppstein1999} that a linear number of topological cubes
suffices, and Bern, Eppstein \& Erickson proved the existence of a
(pseudo-)shellable topological hexa mesh~\cite{BernEppsteinErickson}.

\subsection{Parity change}
Another interesting question deals with the parity of the number of
facets of a mesh.  For quad meshes there are several known parity
changing operations, that is, operations that change the numbers of
facets without changing the boundary. In \cite{BernEppsteinErickson},
Bern, Eppstein \& Erickson raised the following questions:
  \begin{compactenum}[(i)]
  \item Are there geometric quad meshes with geometric hexa meshes of
    both parities?
    
  \item Is there a \emph{parity changing operation} for geometric hexa
    meshes, which would change the parity of the
    number of facets of a cubical $3$-ball without changing the
    boundary?
  \end{compactenum}
%
From the existence of a cubical $4$-polytope with odd number of facets
we obtain positive answers to these questions.

\begin{corollary}\ 
   \begin{compactenum}[\ \rm (i)]
   \item Every combinatorial $3$-cube has a cubification with an
     {even} number of facets.  Furthermore, this cubification is
     regular and even Schlegel.
    \item  Every combinatorial $3$-cube is a facet of a cubical
         $4$-polytope with an {odd} number of facets.
       \item There is a parity changing operation for geometric hexa
         meshes.
   \end{compactenum}
\end{corollary}

\begin{proof}
  For (ii) let $F$ be a combinatorial $3$-cube and $P$ a cubical
  $4$-polytope with an odd number of facets.  By
  Lemma~\ref{lem:connector_cube} there is a combinatorial $4$-cube $C$
  that has both $F$ and a projectively regular $3$-cube $G$ as facets.
  Let $F'$ be an arbitrary facet of $P$. Then there is a combinatorial
  $4$-cube $C'$ that has both $F'$ and a projectively regular $3$-cube
  $G'$ as facets. Then the connected sum of $P$ and $C$ based on the
  facet $F'$ yields a cubical $4$-polytope $P'$ with an odd number of
  facets, and with a projectively regular $3$-cube $G''$ as a facet.
  The connected sum of $P'$ and $C$ glueing the facets $G$ and $G''$
  yields a cubical $4$-polytope with an odd number of facets, and with
  a projective copy of $F$ as a facet.
  
  The statements (i) and (iii) follow from (ii) via Schlegel diagrams.
\end{proof}

\subsection{Flip graph connectivity}
In analogy to the concept of flips for simplicial (pseudo-)manifolds
one can define \emph{cubical flips} for quad or hexa meshes; compare
\cite{BernEppsteinErickson}.  In the meshing terminology the flip
graph is defined as follows.  For any domain with boundary mesh, and a
type of mesh to use for that domain, define the {\em flip graph} to be
a graph with (infinitely many) vertices corresponding to possible
meshes of the domain, and an edge connecting two vertices whenever the
corresponding two meshes can be transformed into each other by a
single flip.

In this framework, the question concerning a parity changing operation
can be phrased as asking for a description of the connected components
of the flip graph.  As an immediate consequence of the corollary above
we obtain the following result.
 
\begin{corollary}\label{thm:flip_graph}
  For every geometric hexa mesh the cubical flip graph has at least
  two connected components.
\end{corollary}




\section{The next step}\label{sec:nextstep}

In this paper we are primarily concerned with the realization of
2-manifold immersions in terms of cubical 4-polytopes, but the
higher-dimensional cases are interesting as well.  For example, one
would like to know whether there are cubical $5$-polytopes with an odd
number of facets.  (There are \emph{no} such $d$-polytopes for $d=6$,
or for $8\le d\le 13$; see the \cite[Sect.~7]{BabsonChan3}.)  For this
we have to realize a normal crossing immersion of $3$-manifold into
$S^4$ by a cubical $5$-polytope with an odd number of quadruple
points. Such immersions exist by an abstract result of Freedman
\cite{freedman78:_quadr_s} \cite{akhmetev96:_freed}, but more
concretely by John Sullivan's observation (personal communication)
that there are regular sphere eversions of the $2$-sphere with exactly
one quadruple point \cite{Smale1958} \cite{Francis1998} and from any
such one obtains a normal-crossing immersion $S^3\immersed S^4$ with a
single quadruple point.

\subsubsection*{Acknowledgements}
We are indebted to Eric Babson, Mark de Longueville, Nikolai Mn\"ev,
Matthias M\"uller-\allowbreak Hanne\-mann, John Sullivan, and Arnold
Wa\ss mer for conversations and comments that turned out to be
essential for this project.

The illustrations of polytopes and balls were produced in the
\textsf{polymake} system of Gawrilow \& Joswig \cite{polymake}, via
\textsf{javaview} of Polthier
et.\,al.\,\cite{PolthierKhademPreussReitebuch01}.  
(A small number of figures was drawn with \textsf{xfig}.)


\begin{small}
\bibliographystyle{siam}
\bibliography{main}

\end{small}


\end{document}